\documentclass{article}
\usepackage[utf8]{inputenc}
\usepackage[english]{babel}
\usepackage{xypic}
\usepackage[all]{xy}
\usepackage{eucal}
\usepackage{latexsym, amscd,eucal,xr, makeidx,stmaryrd,color,mathrsfs,amsmath, amsfonts,amssymb,bbm,stmaryrd,amsthm,graphicx,pstricks,xypic,comment,graphics}
\usepackage{wrapfig}
\usepackage{color,epsfig,enumerate}
\usepackage{xcolor}
\usepackage{caption}
\usepackage{hyperref}
\xyoption{2cell}

\newtheorem{theorem}{Theorem}[section]
\newtheorem{thm}[theorem]{Theorem}
\newtheorem{lemma}[theorem]{Lemma}

\newtheorem{cor}[theorem]{Corollary}
\newtheorem{prop}[theorem]{Proposition}

\newtheorem{defn}[theorem]{Definition}

\theoremstyle{definition}
\newtheorem{remark}[theorem]{Remark}
\newtheorem{construction}[theorem]{Construction}
\newtheorem{example}[theorem]{Example}
\newtheorem{notation}[theorem]{Notation}
\newtheorem{warning}[theorem]{Warning}

\numberwithin{equation}{section}

\newcommand{\C}{\mathcal{C}}

\newcommand{\Sp}{\mathit{Sp}}

\newcommand{\V}{\mathcal{V}}

\newcommand{\Fin}{\mathit Fin_*}

\newcommand{\onefin}{\langle 1 \rangle}

\newcommand{\aff}{\mathit{N(AffSm}^{ft}(k))}

\newcommand{\M}{\mathcal{M}}

\newcommand{\Op}{\mathcal{O}}

\newcommand{\Ass}{\mathcal{A}ss}

\newcommand{\Opmonoidal}{\mathcal{O}^{\otimes}}

\newcommand{\Dmonoidal}{\mathcal{D}^{\otimes}}

\newcommand{\Spaces}{\mathcal{S}}

\newcommand{\stk}{\mathcal{S}\mathcal{H}(k)}

\newcommand{\stmonoidal}{\mathcal{S}\mathcal{H}(S)^{\otimes}}
\newcommand{\stmonoidalk}{\mathcal{S}\mathcal{H}(k)^{\otimes}}
\newcommand{\stnck}{\mathcal{S}\mathcal{H}_{nc}(k)}
\newcommand{\stncmonoidalk}{\mathcal{S}\mathcal{H}_{nc}(k)^{\otimes}}
\newcommand{\stncmonoidal}{\mathcal{S}\mathcal{H}_{nc}(S)^{\otimes}}
\newcommand{\stnc}{\mathcal{S}\mathcal{H}_{nc}(S)}

\newcommand{\nck}{\mathcal{N}cS(k)}

\newcommand{\dg}{\mathcal{D}g(k)}

\newcommand{\refnci}[1]{[I,\ref{I-#1}]}

\newcommand{\A}{\mathcal A}

\newcommand{\W}{\mathcal W}
\newcommand{\X}{\mathcal X}
\newcommand{\Y}{\mathcal Y}
\newcommand{\UU}{\mathcal U}

\begin{document}

\title{ Noncommutative Motives II:\\ \LARGE{$K$-Theory and Noncommutative Motives}}

\author{Marco Robalo\footnote{Supported by the Portuguese Foundation for Science and Technology - FCT Project SFRH / BD / 68868 / 2010 } \\ \small{I3M} \\  \small{Universit\'e de Montpellier2}}

\maketitle

\hfill \textit{À minha avó Diamantina...}\\

\hfill \textit{ Um Eterno Obrigado.}\\

\begin{abstract}

We continue the work initiated in \cite{nc1}, where we introduced  a new stable symmetric monoidal $(\infty,1)$-category $\stncmonoidalk$ encoding a motivic stable homotopy theory for the noncommutative spaces of Kontsevich \cite{kontsevich3, kontsevich1, kontsevich2} and obtained a canonical  monoidal colimit-preserving functor $\mathcal{L}^{\otimes}:\stmonoidalk\to \stncmonoidalk$ relating this new theory to the $(\infty,1)$-category $\stk$ encoding the theory of Morel-Voevodsky \cite{voevodsky-morel, Voevodsky-icm}.  For a scheme $X$ this map recovers the dg-derived category of perfect complexes $L_{pe}(X)$.

In this sequel we address the study of the different flavours of algebraic $K$-theory of dg-categories. As in the commutative case, these can be understood as spectral valued $\infty$-presheaves over the category of noncommutative smooth spaces and therefore provide objects in $\stnck$ once properly localized. Our first main result is the description of non-connective $K$-theory of dg-categories \cite{schlichting-negative, MR2822869} as the noncommutative Nisnevich sheafification of connective $K$-theory. In particular it follows that its further $\mathbb{A}^1$-localization is an object in $\stnck$. As a corollary of the recent developments in \cite{Anthony-thesis}, we prove that this object is a unit for the monoidal structure.

As a first immediate corollary we obtain a precise proof for an original conjecture of Kontsevich claiming that $K$-theory gives the correct mapping spaces in noncommutative motives. A second major consequence is the discovery of a canonical factorization of our functor $\mathcal{L}^{\otimes}:\stmonoidalk\to \stncmonoidalk$ through $Mod_{\mathcal{K}H}(\stk)$ - the $(\infty,1)$-category of modules over the commutative algebra object $\mathcal{K}H $ in $\stk$ representing homotopy invariant algebraic $K$-theory of schemes. To conclude and as a corollary of the results in \cite{riou-spanierwhitehead} we show that if  $k$  is a field admitting resolutions of singularities, this factorization is fully faithful, so that, at the motivic level, no information (below $K$-theory) is lost by passing to the noncommutative world.

\end{abstract}

\setcounter{tocdepth}{2}
\tableofcontents

\section{Introduction}

\subsection{Notations}
This paper is a sequel to our previous work \cite{nc1}. We continue to use the theory of quasi-categories \cite{joyal-article, lurie-htt} as a model for $(\infty,1)$-categories and will  follow the same notations, conventions and universe considerations. In particular, whenever we say "commutative diagram in a quasi-category $\C$", we actually mean the existence of an higher cell in $\C$ providing the commutativity. For the reader's convenience,  references to \cite{nc1} will appear as [I, x.xx]. Moreover,  we will continue to make extensive use of the tools developed by J. Lurie in \cite{lurie-htt, lurie-ha} and assume the reader is familiar with those.\\

 For the rest of this paper we will fix $k$ a base commutative ring. 

\subsection{Previous Work}
To start with, let us briefly and suitably summarize our previous results :

\begin{enumerate}[1)]

\item Let $N(Sm^{ft}(k))^{\times}$ be the category of smooth schemes of finite type over $k$, endowed with the monoidal structure induced by the cartesian product. We provided a universal characterization for the symmetric monoidal $ (\infty,1)$-category $\stmonoidalk$ encoding the motivic stable homotopy theory of Morel-Voevodsky \cite{voevodsky-morel, Voevodsky-icm}: namely,  for any pointed presentable symmetric monoidal $(\infty,1)$-category $\Dmonoidal$, the composition with the canonical monoidal map $N(Sm^{ft}(k))^{\times}\to \stmonoidal$ 

\begin{equation}
Fun^{\otimes,L}(\stmonoidalk, \Dmonoidal)\to Fun^{\otimes}(N(Sm^{ft}(k))^{\times}, \Dmonoidal)
\end{equation}

\noindent is fully-faithful and its image consist of those monoidal functors satisfying Nisnevich descent, $\mathbb{A}^1$-invariance and such that the cofiber of the image of the point at $\infty$ in $\mathbb{P}^1$ is an tensor-invertible object  (see \refnci{universalpropertymotives}). 

Throughout this paper, it will be convenient to use the equivalent description of $\stmonoidalk$  given in \refnci{usingpresheavesofspectra1}, namely, as the composition of monoidal functors

\begin{equation}
\xymatrix{
\aff^{\times} \ar@{^{(}->}[r]^{(\Sigma^{\infty}_{+}\circ j)^{\otimes}} &Fun(\aff^{op},  \widehat{\Sp})^{\otimes} \ar@/^/[d]^{l_{Nis}^{\otimes}}&\\
& \ar@{^{(}->}[u]   Fun_{Nis}(\aff^{op},  \widehat{\Sp})^{\otimes}  \ar@/^/[d]^{l_{\mathbb{A}^1}^{\otimes}}&\\
& \ar@{^{(}->}[u] Fun_{Nis,\mathbb{A}^1}(\aff^{op},  \widehat{\Sp})^{\otimes} \ar[d]^{\Sigma_{\mathbb{G}_m}^{\otimes}}&\\
&Fun_{Nis,\mathbb{A}^1}(\aff^{op},  \widehat{\Sp})^{\otimes}[\mathbb{G}_m^{-1}]\simeq \stmonoidalk&
}
\end{equation}

\noindent with $\aff^{\times}$ the nerve of the standard category of smooth affine schemes of finite type over $k$ endowed with the cartesian product and $\widehat{\Sp}$ the big $(\infty,1)$-category of spectra, obtained from the stabilization of the big $(\infty,1)$-category of spaces $\widehat{\Spaces}$. The first  map $\Sigma^{\infty}_{+}\circ j$ is the composition of the Yoneda's map with the stablization map \footnote{ More precisely, it is the composition

\begin{equation}
\xymatrix{
\aff\ar@{^{(}->}[r]^j &\mathcal{P}^{big}(\aff)\ar[r]^{(-)_{+}} &\mathcal{P}^{big}(\aff)_{\ast} \ar[r]^{\Sigma^{\infty}}& Fun(\aff^{op}, \widehat{\Sp})
}
\end{equation}

\noindent  with $j$ the Yoneda's embedding , $(-)_{+}$ the pointing map , $\Sigma^{\infty}$ the stabilization and we recall that $Stab(\mathcal{P}^{big}(\aff))\simeq Fun(\aff^{op}, \widehat{\Sp})$ .} and the second and third maps are, respectively, the Nisnevich and $\mathbb{A}^1$ reflexive monoidal localizations. The last map is the formal inversion of $\mathbb{G}_m$ with respect to the monoidal structure (see \refnci{defformal}). By our previous results, the last step can be identified with the stabilization with respect to the multiplication by $\mathbb{G}_m$ (see \refnci{main5} for the complete details). Because of the Adjoint Functor Theorem \cite[Corollary 5.5.2.9]{lurie-htt}, this last functor also has a right adjoint which we will denote as $\Omega^{\infty}_{\mathbb{G}_m}$.

\begin{remark}
\label{allarestable}
Notice that all the $(\infty,1)$-categories in the vertical column of the previous diagram are stable. The first, $Fun(\aff^{op},\widehat{\Sp})$, is stable because limits and colimits are computed objectwise in $\widehat{\Sp}$ (see \cite[5.1.2.3]{lurie-htt}), which is stable.  The others require a small discussion. In general, a reflexive localization of a stable $(\infty,1)$-category is stable if and only if the localization functor is left exact \cite[Lemma 1.4.4.7]{lurie-ha}. Here we will use a different argument rather than the discussion of left-exactness. Namely, we use the description of stability for the pointed presentable setting given by our method of inverting the topological circle $S^1$, as explained in \refnci{monoidalstabilization}. In this case, the fact that the two reflexive localizations are stable follows from the combination of  $i)$ the fact that  $Fun(\aff^{op},\widehat{\Sp})$ is stable, $ii)$ the fact that each functor in the diagram is monoidal and colimit preserving and finally $iii)$ the fact that both the Nisnevich and $\mathbb{A}^1$ localizations are pointed, which follows from the fact that the zero object in $Fun(\aff^{op},\widehat{\Sp})$ is of course Nisnevich local and $\mathbb{A}^1$ invariant. The last is stable because $\mathbb{G}_m$ is a symmetric object (see \refnci{symmetric} and \refnci{remarksymmetric},  together with \refnci{corolariodacaca}.
\end{remark}

\item
Let $\dg^{idem, \otimes}$ be the symmetric monoidal $(\infty,1)$-category underlying the homotopy theory of small dg-categories equipped with the Morita model structure of \cite{tabuada-invariantsadditifs} and let $\dg^{ft}\subseteq \dg^{idem}$ be the full subcategory spanned by the small dg-categories of finite type of Toën-Vaquié \cite{toen-vaquie}. In our previous work we defined the $(\infty,1)$-category $\nck$ of smooth noncommutative spaces as the opposite of $\dg^{ft}$ and defined the tensor product of two noncommutative spaces as the tensor product between their underlying dg-categories of finite type. Finally, we explain how to write the formula $X\mapsto L_{pe}(X)$ sending a smooth affine scheme $X$ to its associated dg-category of perfect complexes $L_{pe}(X)$ as a monoidal functor $L_{pe}:\aff^{\times}\to \nck^{\otimes}$ (see \refnci{section633}).

\item
We introduced a noncommutative version of the Nisnevich topology \refnci{defncnisnevich} and used it to construct a new symmetric monoidal $(\infty,1)$-category $\stncmonoidalk$ encoding a noncommutative analogue of the motivic stable homotopy theory of Morel-Voevodsky, together with a monoidal functor $\nck^{\otimes}\to \stncmonoidalk$  satisfying a universal property similar to the one given for schemes. More precisely, and following \refnci{usingpresheavesofspectra2}, it can be written as the composition of monoidal functors

\begin{equation}
\label{nc2picapi}
\xymatrix{
\nck^{\otimes}=(\dg^{ft})^{op,\otimes} \ar@{^{(}->}[r]^{(\Sigma^{\infty}_{+}\circ j_{nc})^{\otimes}} &Fun(\dg^{ft}, \widehat{\Sp})^{\otimes} \ar@/^/[d]^{l_{Nis}^{nc,\otimes}}&\\
& \ar@{^{(}->}[u]   Fun_{Nis}(\dg^{ft}, \widehat{\Sp})^{\otimes}  \ar@/^/[d]^{l_{\mathbb{A}^1}^{nc,\otimes}}&\\
& \ar@{^{(}->}[u] Fun_{Nis,L_{pe}(\mathbb{A}^1)}(\dg^{ft}, \widehat{\Sp})^{\otimes} \ar[d]^{\Sigma_{L_{pe}(\mathbb{G}_m)}^{\otimes}}_{\sim}&\\
&Fun_{Nis,L_{pe}(\mathbb{A}^1)}(\dg^{ft}, \widehat{\Sp})^{\otimes}[\mathbb{G}_m^{-1}]\simeq \stncmonoidalk&
}
\end{equation}

\noindent where this time the last map corresponding to the inversion of $L_{pe}(\mathbb{G}_m)$ is already an equivalence (see \refnci{alreadystable}). In particular, an object $F:\dg^{ft}\to \widehat{\Sp}$ is in $\stnck$ if and only if it satisfies the following properties: $a)$ sends Nisnevich Squares of dg-categories to pullback-pushout squares in $\widehat{\Sp}$, $b)$ sends the zero-dg-category to zero and $c)$ for every dg-category $T$, the canonical map $F(T)\to F(T\otimes L_{pe}(\mathbb{A}^1))$  is an equivalence of spectra.

\begin{remark}
\label{allarestable2} 
By the same arguments as in the Remark \ref{allarestable}, all the  $(\infty,1)$-categories in the vertical column of the diagram (\ref{nc2picapi}) are stable.
\end{remark}

\begin{notation}
\label{notationunit}
Let $1_k$ denote the dg-category with a single object and having the ring $k$ (considered as a complex concentrated in degree zero) as its complex of endomorphisms. Its idempotent completion $\widehat{(1_k)}_c$ is a dg-category of finite type canonically equivalent to $L_{pe}(k)$.  It is the unit for the monoidal structure in $\dg^{idem}$ so that $L_{pe}(k)$ seen as a noncommutative space is a unit in $\nck$. Since the construction of the symmetric monoidal $(\infty,1)$-category $\stncmonoidal$  is obtained as sequence of monoidal steps, its unit object is the image of $L_{pe}(k)$ through this sequence. Since $\Sigma^{\infty}_{+}\circ j(L_{pe}(k))$ is Nisnevich  local, it is given by $l_{\mathbb{A}^1}( \Sigma^{\infty}_{+}\circ j(L_{pe}(k)))$. Throughout this paper we will denote it as $1_{nc}$.
\end{notation}

\begin{remark}
\label{nc2strictmodel}(Strictification) 
It is important to remark that an object $F$ in $Fun(\dg^{ft}, \widehat{\Sp})$ can always be identified up to equivalence with an actual strict functor $F_s$ from the category of dg-categories endowed with the Morita model structure  of  \cite{tabuada-invariantsadditifs} to some combinatorial model category whose underlying $(\infty,1)$-category is $\widehat{\Sp}$ (for instance, the big model category of symmetric spectra $Sp^{\Sigma}$ of \cite{MR1695653}), with $F_s$ sending Morita equivalences to weak-equivalences and commuting with filtered homotopy colimits. Indeed, as explained in \refnci{dgideminddgft}, $\dg^{ft}$ generates $\dg^{idem}$ under filtered colimits. Since $\Sp$ admits all small filtered colimits, using \cite[Thm 5.3.5.10]{lurie-htt} we find  an equivalence of $(\infty,1)$-categories between $Fun(\dg^{ft}, \widehat{\Sp})$ and $Fun_{\omega}(\dg^{idem}, \widehat{\Sp})$ - the full subcategory of $Fun(\dg^{idem},\widehat{\Sp})$ spanned by the functors that preserve filtered colimits. Moreover, we have also seen that  $\dg^{idem}$ is the underlying $(\infty,1)$-category of the Morita model structure for small dg-categories (see the discussion in \refnci{morita}). Finally, with the appropriate universe considerations, we can use the strictification result of \cite[1.3.4.25]{lurie-ha} and the characterization of homotopy limits and colimits in a model category as limits and colimits in its underlying $(\infty,1)$-category \cite[1.3.4.24]{lurie-ha} to deduce the existence of a canonical equivalence between $Fun_{\omega}(\dg^{idem}, \widehat{\Sp})$ and the localization along the levelwise equivalences of the category of strict functors from the category of dg-categories to the strict model for spectra $Sp^{\Sigma}$, 
which commute with filtered homotopy colimits and send Morita weak-equivalences to weak-equivalences in $Sp^{\Sigma}$.
\end{remark}

 \item
Finally,  using all the universal properties involved, we were able to build up a homotopy commutative diagram of colimit preserving monoidal functors extending the initial $L_{pe}$

\begin{equation}
\label{diagramaleft}
\xymatrix{
\aff^{\times}\ar@{^{(}->}[d]^{(\Sigma^{\infty}_{+}\circ j)^{\otimes}}\ar[r]^{L_{pe}^{\otimes}}& \nck^{\otimes}\ar@{^{(}->}[d]^{(\Sigma^{\infty}_{+}\circ j_{nc})^{\otimes}}\\
Fun(\aff^{op}, \widehat{\Sp})^{\otimes}\ar[d]^{l_{Nis}^{\otimes}}\ar@{-->}[r]& Fun(\dg^{ft}, \widehat{\Sp})^{\otimes}\ar@<1ex>[d]_{l_{Nis}^{nc, \otimes}}\\
Fun_{Nis}(\aff^{op}, \widehat{\Sp})^{\otimes}  \ar[d]^{l_{\mathbb{A}^1}^{\otimes}}\ar@{-->}[r]& Fun_{Nis}(\dg^{ft}, \widehat{\Sp})^{\otimes}\ar[d]_{l_{\mathbb{A}^1}^{nc,\otimes}} \\
Fun_{Nis, \mathbb{A}^1}(\aff^{op}, \widehat{\Sp})^{\otimes}  \ar[d]^{\Sigma_{\mathbb{G}_m}^{\otimes}}\ar@{-->}[r]& Fun_{Nis, L_{pe}(\mathbb{A}^1)}(\dg^{ft}, \widehat{\Sp})^{\otimes}\ar[d]^{\sim} \\
\stmonoidalk\ar@{-->}[r]^{\mathcal{L}^{\otimes}}&\stncmonoidalk
}
\end{equation}
\end{enumerate}

\noindent providing a canonical mechanism to compare the theory of Morel-Voevodsky with our new approach.

\subsection{Present Work}

Our ultimate goal in this paper is to explore how the comparison mechanism in (\ref{diagramaleft}) can be used to give a canonical interpretation to the various flavours of algebraic $K$-theory of schemes. In order to state our results, we observe first that, due to the Adjoint Functor Theorem (\cite[Corollary 5.5.2.9]{lurie-htt}),  each of the dotted monoidal functors  in (\ref{diagramaleft}) has a right adjoint.  This is because at each level, the source and target $(\infty,1)$-categories are presentable and each dotted map is, by construction, colimit-preserving. Furthermore, since each dotted map is monoidal, these right adjoints are lax-monoidal (see \refnci{monoidaladjoint}). In this case,  together with the lax-monoidal inclusions associated to the reflexive monoidal localizations,  we have a new commutative diagram of lax-monoidal functors

\begin{equation}
\label{diagramaright}
\xymatrix{
Fun(\aff^{op}, \widehat{\Sp})^{\otimes}& \ar[l]_{\M_1^{\otimes}} Fun(\dg^{ft}, \widehat{\Sp})^{\otimes}\\
\ar@{^{(}->}[u]  Fun_{Nis}(\aff^{op}, \widehat{\Sp})^{\otimes}  & \ar@{^{(}->}[u] \ar[l]_{\M_2^{\otimes}} Fun_{Nis}(\dg^{ft}, \widehat{\Sp})^{\otimes}\\
\ar@{^{(}->}[u]  Fun_{Nis, \mathbb{A}^1}(\aff^{op}, \widehat{\Sp})^{\otimes} & \ar@{^{(}->}[u] \ar[l]_{\M_3^{\otimes}} Fun_{Nis, L_{pe}(\mathbb{A}^1)}(\dg^{ft}, \widehat{\Sp})^{\otimes} \\
\stmonoidalk \ar[u]^{\Omega^{\infty,\otimes}_{\mathbb{G}_m}}& \ar[l]_{\M^{\otimes}}\stncmonoidalk \ar[u]^{\sim}
}
\end{equation}

Let us present some remarks that will be useful all along the text. 

\begin{remark}
\label{Mpreservescolimits}
The first functor $\M_1$ commutes with small colimits. We can deduce this either from the fact that colimits in $Fun(\aff^{op}, \widehat{\Sp})$ and in $Fun(\dg^{ft}, \widehat{\Sp})$ are computed objectwise (see \cite[5.1.2.3]{lurie-htt}) or from the spectral enriched version of Yoneda's lemma (\refnci{enrichedyoneda}).
\end{remark}

\begin{remark}
\label{Mcompatibleinternalhom}
All the symmetric monoidal $(\infty,1)$-categories appearing in the previous diagram are closed monoidal (they are presentable). In particular, recall that if $\C_0\subseteq \C$ is a monoidal reflexive localization and if $\C$ admits internal-homs $\underline{Hom}_{\C}$ then $\C_0$ admits internal-homs: given $X$ and $Y$ local, we can easily see that $\underline{Hom}_{\C}(X,Y)$ is also local and works as an internal-hom in $\C_0$.
 
We observe that each functor $\M_\ast$ is compatible with the respective internal-homs, in the sense that at each level, for every object $X$ on the left and $F$ on the right, we have

\begin{equation}
\label{nc2internalhomispreserved}
\M_\ast(\underline{Hom}_\ast(\mathcal{L}_\ast(X), F))\simeq \underline{Hom}_{\ast}(X, \M_\ast(F))
\end{equation}

\noindent where $\mathcal{L}_\ast$ denotes the respective monoidal left adjoint appearing in the diagram (\ref{diagramaleft})
\end{remark}

\begin{remark}
\label{nc2enrichedyoneda}
Recall that if $\C$ is a stable $(\infty,1)$-category, it is canonically enriched over $\Sp$. More precisely, the universal property of the stabilization tells us that the composition with $\Omega^{\infty}:\Sp\to \Spaces$ induces an equivalence of $(\infty,1)$-categories $Exc_{\ast}(\C, \Sp)\simeq Exc_{\ast}(\C, \Spaces)$ (see \cite[1.4.2.22]{lurie-ha}). In particular, this provides for  any object $X$  an essentially unique factorization of the functor $Map_{\C}(X,-):\C\to \Spaces$ as

\begin{equation}
\xymatrix{
\C\ar[rr]^{Map_{\C}(X,-)}\ar@{-->}[d]_{Map_{\C}^{Sp}(X,-)}&& \Spaces\\
\Sp \ar[urr]_{\Omega^{\infty}}&
}
\end{equation}

\noindent such that for any object $Y$,  the spectrum $Map_{\C}^{Sp}(X,Y)$ can be identified with the collection of spaces  $\{Map_{\C}(X, \Sigma^n Y)\}_{n\in \mathbb{Z}}$. This holds for any universe.

In particular, and thanks to the enriched version of  Yoneda's lemma for spectral presheaves (see the Remark \refnci{enrichedyoneda}), given an object $F\in Fun(\dg^{ft}, \widehat{\Sp})$, we have for each scheme $X$ an equivalence of spectra 

\begin{equation}
Map_{Fun(\aff^{op}, \widehat{\Sp})}^{Sp}(\Sigma^{\infty}_+\circ j(X), \M_1(F))\simeq Map^{Sp}_{ Fun(\dg^{ft}, \widehat{\Sp})}(\Sigma^{\infty}_+\circ j_{nc}(L_{pe}(X)), F)\simeq F(L_{pe}(X))
\end{equation}

\noindent so that $\M_1(F)$ can be thought of as a restriction of $F$ to the commutative world. The same is valid for $\M_2$ and $\M_3$ because the upper vertical arrows are inclusions.
\end{remark}

This mechanism allows us to restrict noncommutative invariants to the commutative world. 

\begin{example}
An important noncommutative invariant is the Hochschild homology of dg-categories. Thanks to the works of B. Keller in \cite{keller-exact} and as explained in the Remark \ref{nc2strictmodel} this invariant can be completely encoded by means  of an $\infty$-functor $HH: \dg^{ft}\to \widehat{\Sp}$. Another important example is the so called \emph{periodic cyclic homology of dg-categories} $HP$. It follows from the famous $HKR$ theorem that the restriction of $HP$ to the commutative world recovers the classical de Rham cohomology of schemes. For more details see the discussion in \cite[Section 3.1]{Anthony-thesis}.
\end{example}

In this paper we will be interested in the restriction of the various algebraic $K$-theories of dg-categories. As we shall explain below, all of them live as objects in $Fun(\dg^{ft}, \widehat{\Sp})$. There are two of primary relevance: 

\begin{itemize}
\item $K^c$, encoding the \emph{connective} $K$-theory given by Waldhausen's $S$-construction. See the discussion in Section \ref{connectivektheorydg} below. 
\item $K^S$, encoding the \emph{non-connective} $K$-theory of dg-categories as defined in \cite{MR2822869} using the adaptation of the Schlichting's framework of \cite{schlichting-negative} to the context of dg-categories. (see the discussion in Section \ref{nonconnectivektheorydg}). By construction, this functor comes naturally equipped with a canonical natural transformation $K^c\to K^S$ which is an equivalence in the connective part.
\end{itemize}

For the first one, it follows immediately from the spectral version of  Yoneda's lemma mentioned above and from the definition in \cite[Section 3]{thomasonalgebraic} that $\M_1(K^c)$ recovers the connective algebraic $K$-theory of schemes. The second one, by the comparison result \cite[Theorem 7.1]{schlichting-negative}, recovers the non-connective K-theory of schemes of Bass-Thomason-Trobaugh of \cite{thomasonalgebraic}. The constructon of $K^S$ in \cite{MR2822869} using the methods of \cite{schlichting-negative} is somehow ad-hoc. Our first main result explains how the non-connective version of $K$-theory $K^S$ can be canonically obtained from the connective version $K^c$ as a result of forcing our noncommutative-world version of Nisnevich descent.

\begin{thm}
\label{teorema1}
The canonical morphism $K^c\to K^S$ presents non-connective $K$-theory of dg-categories as the (noncommutative) Nisnevich localization of connective $K$-theory.
\end{thm}

To prove this result we will first check that $K^S$ is Nisnevich local. This follows from the well-known localization theorem for non-connective $K$-theory (see the Corollary \ref{nc2negativektheoryisnisnevichlocal} bellow). The rest of the proof will require a careful discussion concerning the behavior of the noncommutative Nisnevich localization.  There are two main ingredients:

\begin{enumerate}[Step 1)]
\item \emph{Every Nisnevich local $F:\dg^{ft}\to \widehat{\Sp}$ is determined by its connective part by means of the Bass exact sequences}. More precisely, we show that every Nisnevich local functor $F:\dg^{ft}\to \widehat{\Sp}$ satisfies the familiar Bass exact sequences for any integer $n$. We will see that the proof in \cite{thomasonalgebraic} can be easily adapted to our setting. Namely, we start by showing that every Nisnevich local $F$ satisfies the Projective Bundle theorem. This result is central and appears as a consequence of one of the most important features of the noncommutative world, namely, the fact that Nisnevich coverings of non-geometrical origin are allowed, in particular, those appearing from semi-orthogonal decompositions and exceptional collections. The projective bundle theorem is a direct consequence of the existence of an exceptional collection on  $L_{pe}(\mathbb{P}^1) $ generated by the sheaves  $\mathcal{O}_{\mathbb{P}^1}$ and  $\mathcal{O}_{\mathbb{P}^1}(-1)$ (see \cite{beilinsonprojective}).  Its existence forces  the image of $L_{pe}(\mathbb{P}^1)$ in $Fun_{Nis}(\dg^{ft}, \widehat{\Sp})$  to become equivalent to the direct sum $L_{pe}(k)\oplus L_{pe}(k)$. To complete the proof we proceed as in  \cite[Theorem 6.1]{thomasonalgebraic} and explain how this direct sum decomposition can be suitably  adapted in order to extract the familiar Bass exact sequences out of the classical Nisnevich covering of $\mathbb{P}^1$ by two affine lines.

\item \emph{The connective truncation of the localization map $K^c\to l_{Nis}^{nc}(K^c)$ is an equivalence} \footnote{Recall that $\widehat{\Sp}$ has a natural t-structure $(\widehat{\Sp}_{\geq 0}, \widehat{\Sp}_{\leq -1})$ with $\widehat{\Sp}_{\geq 0}$ the full subcategory spanned by connective spectra. As a consequence, the inclusion $\widehat{\Sp}_{\geq 0}\subseteq \widehat{\Sp}$ (resp. $\widehat{\Sp}_{\leq -1}\subseteq \widehat{\Sp}$) admits a right adjoint $\tau_{\geq 0}$ (resp. left adjoint $\tau_{\leq -1}$).  In particular, we have an induced adjunction

$$
\xymatrix{Fun(\dg^{ft}, \widehat{\Sp}_{\geq 0})\ar@{^{(}->}[r]& \ar@/^1pc/[l]_{\tau_{\geq 0}} Fun(\dg^{ft},\widehat{\Sp}) }
$$

\noindent with $\tau_{\geq 0}$ a right adjoint to the inclusion.}. In other words, the information stored in the connective part of $l_{Nis}(K^c)$ remains the information of connective $K$-theory. We will prove something a bit more general, namely, that this property holds not only for $K^c$ but  for the whole class of functors $F:\dg^{ft}\to \widehat{\Sp}$ satisfying the formal properties of $K^c$, namely, having values in connective spectra and sending Nisnevich squares of dg-categories to pullback squares of connective spectra (for $K^c$ this follows from the fibration theorem of  Waldhausen \cite[1.6.4]{waldhausen-ktheoryofspaces} - see Prop. \ref{kcsemilocal} below.). These will be called \emph{connectively-Nisnevich local}. We prove that the connective truncation functor induces a canonical equivalence between the theory of connective-Nisnevich functors and that of Nisnevich functors (see \ref{nc2proposicao1}). For this we will show that if $F$ is  connectively-Nisnevich local, its noncommutative Nisnevich localization $l_{Nis}^{nc}(F)$ is equivalent to $F^{B}$ -  the more familiar $B$-construction of Thomason  of \cite[Def. 6.4]{thomasonalgebraic}. 

\end{enumerate}

\begin{remark}
Since the functor $\M_2$ in the diagram (\ref{diagramaright}) sends Nisnevich local objects to Nisnevich local objects, our Theorem \ref{teorema1} provides a new proof that the spectral presheaf giving the Bass-Thomason-Trobaugh $K$-theory of schemes satisfies Nisnevich descent.
\end{remark}

We can now go one step further and consider the $\mathbb{A}^1$-localization of $K^S$. We will prove that

\begin{thm}
\label{12}
$\M_3(l_{\mathbb{A}^1}^{nc}(K^S))$ is the Nisnevich local $\mathbb{A}^1$-invariant spectral presheaf giving Weibel's homotopy invariant $K$-theory of schemes of \cite{weibel-homotopyinvariantktheory}. In particular, $\M(l_{\mathbb{A}^1}^{nc}K^S)$ is canonically equivalent to the object $\mathcal{K}H$ in $\stk$ studied in \cite{Voevodsky-icm} and in \cite{cisinski-descentpar} representing homotopy invariant algebraic $K$-theory of schemes.
\end{thm}

The proof of this result follows immediately from the results in \cite{cisinski-descentpar} and from our Theorem \ref{teorema1} using a nice description of the $\mathbb{A}^1$-localization functors. This will be done in Section \ref{comparisonA1}.\\

Our second main result in this paper is a new representability theorem for $K$-theory. 

\begin{thm}
\label{teorema2}
The further localization  $l_{\mathbb{A}^1}^{nc}(K^S)$ is a unit for the monoidal structure in $\stncmonoidalk$.
\end{thm}

In \cite{Anthony-thesis}, the author constructs an $\mathbb{A}^1$-equivalence between the split and the standard versions of Waldhausen's $S$-construction. In Section \ref{proofteorema2} we will explain how this $\mathbb{A}^1$-equivalence appears in our context and how the theorem follows as a consequence.\\

We deduce the following immediate corollaries

\begin{cor}(Kontsevich)
Let $\X$ and $\Y$ be two noncommutative spaces with $\Y$ smooth and proper \refnci{nc1smoothandproper}. Then, there is a natural equivalence of spectra

\begin{equation}
Map^{Sp}_{\stnck}(\X, \Y)\simeq (l^{nc}_{\mathbb{A}^1}K^S)(T_{\X}\otimes \check{T_{\Y}})
\end{equation}

\noindent where we identify $\X$ and $\Y$ with their images in $\stnck$ and where $T_{\X}$ (resp. $\check{T_{\Y}}$) denotes the dg-category of finite type associated to $\X$ (resp. the dual of the dg-category associated to $\Y$).
\begin{proof}
This follows directly from the spectral version of the Yoneda's lemma and from our Theorems \ref{teorema1} and \ref{teorema2}, together with the fact that a smooth and proper noncommutative space is dualizable (and vice-versa).
\end{proof}
\end{cor}

\begin{cor}
The object $\mathcal{K}H\in \stk$ representing homotopy algebraic $K$-theory is equivalent to $\M(1_{nc})$. In particular, for each scheme $X$ we have an equivalence of spectra

\begin{equation}
KH(X)\simeq Map_{\stk}^{Sp}(\Sigma^{\infty}_+\circ j(X), \mathcal{K}H)\simeq Map^{Sp}_{\stnck}(\Sigma^{\infty}_+\circ j_{nc}(L_{pe}(X)),1_{nc})
\end{equation}

\end{cor}

At this point we should emphasize that a different representability result for connective $K$-theory is already known from the thesis of G. Tabuada \cite{tabuada-higherktheory} and for non-connective $K$-theory from his  later works with D.C. Cisinski \cite{MR2822869}. Our setting and proofs are independent of theirs. In the appendix of this paper we describe the relation between the two approaches.\footnote{The contents of the appendix are independent of the rest of paper.} The main advantage of our theory is the existence of a canonical comparison with the original approach of Morel-Voevodsky and our new representabilty theorem  brings some immediate consequences to the nature of this comparison. Namely, since $\M$ is lax-monoidal, the object $\mathcal{K}H\simeq \M(1_{nc})$  acquires a canonical structure of commutative algebra-object in $\stk$ induced by the trivial algebra structure on the unit object $1_{nc}$. In this case, the comparison functor
$\mathcal{L}^{\otimes}:\stmonoidalk\to \stncmonoidalk$ admits a canonical colimit preserving monoidal factorization

\begin{equation}
\xymatrix{
\stmonoidalk\ar[r]^{\mathcal{L}^{\otimes}}\ar[d]^{-\otimes \mathcal{K}H}& \stncmonoidalk\ar@/^/[drr]^{(-\otimes 1_{nc})\simeq Id} \ar[d]_{-\otimes \mathcal{L}(\mathcal{K}H)}&&\\
Mod_{ \mathcal{K}H}(\stk)^{\otimes}\ar[r] &Mod_{\mathcal{L}(\mathcal{K}H)}(\stnck)^{\otimes} \ar[rr]_{-\otimes_{\mathcal{L}(\mathcal{K}H)} 1_{nc}}&& Mod_{1_{nc}}(\stnck)^{\otimes}
}
\end{equation}

\noindent where the first lower map is the monoidal functor induced by $\mathcal{L}$ at the level of modules (see  \refnci{changeofmodulesundermonoidalfunctor}) and the last map is base-change with respect to the canonical morphisms of algebra objects given by the counit of the adjunction $\mathcal{L}(\mathcal{K}H)\simeq \mathcal{L}\circ \M(1_{nc})\to 1_{nc}$ \footnote{Notice that the adjunction $(\mathcal{L}, \M)$ extends to an adjunction between the $(\infty,1)$-categories of commutative algebra-objects, so that this counit map is a morphism of algebras. In particular, we can perform base-change with respect to it.}.  We will write $\mathcal{L}_{\mathcal{K}H}$ for this factorization.\\

\begin{warning}
We will not prove here that the commutative algebra structure in $\mathcal{K}H$ obtained from our arguments is the same as the one already appearing in the literature and deduced from different methods (for instance, see \cite{gepner-algebraiccobordismalgebraicKtheory,1010.3944}). However, we believe that the arguments used in \cite{1010.3944} also work in the $\infty$-categorical setting, so that our algebra structure should match the standard one.
\end{warning}

Our representability result has the following corollary showing that under the existence of resolutions of singularities the passage to the noncommutative world produces no loss of information from the $K$-theoretic viewpoint.

\begin{cor}
Let $k$ be a field admitting resolutions of singularities. Then the canonical map 

\begin{equation}
\mathcal{L}_{\mathcal{K}H}: Mod_{\mathcal{K}H}(\stk)\to \stnck
\end{equation}

\noindent is fully faithful.

\begin{proof}
Thanks to the results of \cite{riou-spanierwhitehead} the family of dualizable objects in $\stk$ is a family of $\omega$-compact generators for the stable $(\infty,1)$-category $\stk$ in the sense of the Proposition \refnci{cg}. Thanks to the Proposition \refnci{compactgeneratorscategoriesofmodules}, the collection of all objects in the stable $(\infty,1)$-category $Mod_{\mathcal{K}H}(\stk)$ of the form $X\otimes \mathcal{K}H$ with $X$ dualizable in $\stk$ is a family of $\omega$-compact generators in the sense of \refnci{cg}. Since the functor $(-\otimes \mathcal{K}H)$ is monoidal, the objects $X\otimes \mathcal{K}H$ are dualizable in $Mod_{\mathcal{K}H}(\stk)$ and as $\mathcal{L}_{\mathcal{K}H}$ is monoidal, their image in $\stnck$ is dualizable and therefore compact (using the fact the monoidal structure is compatible with colimits in each variable). By the Proposition \refnci{equivalencecompactgenerators} we are now reduced to show that  $\mathcal{L}_{\mathcal{K}H}$ is fully faithful when restricted to the full subcategory spanned by all the objects of the form  $X\otimes \mathcal{K}H$ with $X$ dualizable in $\stk$. This follows from the canonical chain of equivalences

\begin{eqnarray}
Map_{Mod_{\mathcal{K}H}(\stk)}(X\otimes \mathcal{K}H, Y\otimes \mathcal{K}H)\simeq Map_{\stk}(X, Y\otimes \mathcal{K}H)\simeq Map_{\stk}(X\otimes \check{Y},\mathcal{K}H)\simeq \\
\simeq Map_{\stnck}(\mathcal{L}(X\otimes \check{Y}),1_{nc})\simeq Map_{\stnck}(\mathcal{L}(X)\otimes\mathcal{L}(\check{Y}),1_{nc})\simeq Map_{\stk}(\mathcal{L}(X)\otimes\check{\mathcal{L}(Y)},1_{nc})\simeq\\
\simeq Map_{\stnck}(\mathcal{L}(X)\otimes,\mathcal{L}(Y))
\end{eqnarray}

\noindent where we use the adjunction properties, the fact that $\mathcal{K}H\simeq \M(1_{nc})$ and the fact that $\mathcal{L}$ is monoidal and therefore preserves dualizable objects. This concludes the proof.
\end{proof}
\end{cor}

The content of this result as been known after a quite some time already. I think particularly of B.Toen, M. Vaquie and G. Vezzosi and also of D-C. Cisinski and G. Tabuada. Moreover, it is believe to be true over any base. I will address this question in a forthcoming part of my thesis by exploring the existence of a six-functors formalism for $\stncmonoidalk$.

\subsection{Acknowledgments and Credits}

This work is a second part of my doctoral thesis under the direction of Bertrand Toen in Montpellier. I want to expresss my profound gratitude and mathematical debt to him. The process of writing a thesis with him has been and continues to be one of the most enriching and motivating experiences of my life.

I have also a huge gratitude to Gabriele Vezzosi and Michel Vaquié. The idea of this possible approach to noncommutative motives is due to them, together with Bertrand and I'm thankful to them for letting me work in the subject. I should also express the huge influence of the works and writings of J.Lurie. This influence is omnipresent. Also, I thank Denis-Charles Cisinski for an email explanation about the explicit formula for the $\mathbb{A}^1$-localization.

I also want to express my deepest gratitute to Anthony Blanc.  I have learned a lot about $K$-theory from him, through many long discussions, beers, burgers and shadokian conversations. More importantly,  the key step leading to the second main result in this paper is due to him \cite{Anthony-thesis}. Finally, I'm also grateful to Brad Drew for his multiple valuable comments and for the many discussions we had during these last months concerning the intersection of our works and more generally, the subject of motives.

\section{$K$-theory Preliminaries}
\label{prelim}

\subsection{Connective $K$-theory - an historical overview}
\label{prelimkconnective}

$K$-theory was discovered by A. Grothendieck during his attempts to generalize the classical Riemman-Roch theorem to a more comprehensible form (see \cite{MR0116022,MR1644323}).  Given an abelian category $E$ he was led to consider an abelian group $K_{0}(E)$ together with a map $\theta: Obj(E)\to K_{0}(E)$ universal with respect to the following property: for any exact sequence $a\to b\to c$ in $E$ we have $\theta(b)=\theta(a)+ \theta(c)$.

The essential insight leading to the introduction of higher $K$-theory groups is the observation by Quillen \cite{MR0338129} that the group law on $K_{0}(E)$ can be understood as the $\pi_0$-reminiscent part of a grouplike homotopy commutative law on a certain space $K(E)$. Following his ideas, for any "exact category" $E$ we are able to define such a $K$-theory space $K(E)$ whose homotopy groups $\pi_n(K(E))$ we interpret as level $n$ $K$-theoretic information. In particular, this methodology allows us to attach a $K$-theory space to every scheme $X$ using the canonical structure of exact category on $E=Vect(X)$. 

An important step in this historical account is a theorem by Segal \cite[3.4]{MR0353298} (and its later formulation in terms of model categories in \cite{MR513569}) establishing an equivalence between the homotopy theory of grouplike homotopy commutative algebras in spaces and the homotopy theory of connective spectra. This is the reason why connective spectra is commonly used as the natural target for $K$-theory and the origin of the term "connective". In the modern days this equivalence can be stated by means of an equivalence of $(\infty,1)$-categories, namely, between the $(\infty,1)$-category $CAlg^{grplike}(S)$ and the $(\infty,1)$-category $\Sp_{\geq 0}$ (see \cite[Theorem 5.1.3.16 and Remark 5.1.3.17]{lurie-ha}).

Technical reasons and possible further applications led Waldhausen \cite{waldhausen-ktheoryofspaces} to extend the domain of $K$-theory from exact categories to what we nowadays call "Waldhausen categories".  Grosso modo, these are triples $(\C, W, Cof(\C))$ where $\C$ is a classical category having a zero object and both $W$ and $Cof(\C)$ are  classes of morphisms in $\C$, respectively called "weak-equivalences" and "cofibrations".  These triples are  subject to certain conditions which we will not specify here. The core of Waldhausen's method to construct a $K$-theory space out of this data is the algorithm known as the "$S$-Construction" which we review here very briefly: 

\begin{construction}($S$-Construction)
\label{Sconstruction0}
Let $Ar[n]$ be the category of arrows in the linear category $[n]$. In more explicit terms it can be described as the category where objects are pairs $(i,j)$ with $i\leq j$ and there is one morphism $(i,j)\to (l,k)$ everytime  $i\leq l$ and $j\leq k$. Let now $(\C, W, Cof(\C))$ be a Waldhausen category. We let $S_n(\C)$ denote the full subcategory of all functors $ Fun(Ar[n], \C)$ spanned by those functors $A$ verifying:

\begin{enumerate}
\item $A(i,i)$ is a zero object of $\C$ for all $0\leq i\leq n$;
\item for any $i$ the maps $A(i,j)\to A(i,k)$ with $j\leq k$ are cofibrations in $\C$; 
\item for any $i\leq j\leq k$ the induced diagram

\begin{equation}
\xymatrix{
A(i,j)\ar[r]\ar[d]&A(i,k)\ar[d]\\
0=A(j,j)\ar[r]& A(j,k)
}
\end{equation}

is a pushout $\C$.

\end{enumerate}

In other words, the objects in $S_n(\C)$ can be identified with sequences of cofibrations of length $n-1$ plus the datum of the sucessive quotients. In particular, $S_0(\C)$ is the category with a single object and $S_1(\C)$ is equivalent to $\C$. Moreover, the collection of categories $\{S_n(\C)\}_{n\in \mathbb{N}}$ assembles together to form a simplicial category $S_{\bullet}(\C)$ carrying at each level a canonical structure of Waldhausen category whose weak-equivalences $W_n$ are the levelwise weak-equivalences in $\C$. By definition, the $K$-theory space of $\C$ is the simplicial set $K^c(\C):= \Omega\, colim_{\Delta^{op}} N(S_n(\C)^{W_n})$ where $S_n(\C)^{W_n}$ denotes the subcategory of $S_n(\C)$ contaning all the objects and only those morphisms which are weak-equivalences and $N$ is the standard nerve functor. By iterating this procedure we can produce a spectrum. For the complete details see \cite{waldhausen-ktheoryofspaces}.
\end{construction}

There is a natural notion of exact functor between Waldhausen categories providing a category $Wald_{Classic}$ and the $K$-theory assignement can be understood as a functor

\begin{equation}
\xymatrix{K^c_{Wald}: N(Wald_{Classic})\ar[r]& Sp^{\Sigma}}
\end{equation}

\noindent where $Sp^{\Sigma}$ is a model category for the $(\infty,1)$-category $\Sp$.\\

Many Waldhausen categories used in practice appear as subcategories of a Quillen model category \cite{quillen} with the cofibrations and weak-equivalences therein. We will denote by $Wald_{Classic}^{Model}$ the full subcategory of $Wald_{Classic}$ spanned by those Waldhausen categories falling into this list of examples. These Waldhausen categories have a special advantage - the factorization axioms for the model category allow us to change the Construction \ref{Sconstruction0} to consider all morphisms in $\C$, not only the cofibrations.

What we risk to call the first era of connective $K$-theory finishes with the works of Thomason-Trobaugh in \cite{thomasonalgebraic} where the machinery of Waldhausen is applied to schemes and it is proven that the connective $K$-theory of a scheme $X$ introduced by Quillen can be recovered 
from the $K$-theory attached to the Waldhausen structure on the category of perfect complexes on the scheme.\\

The modern times for connective $K$-theory begin with the observation that the $K$-theory of a Waldhausen datum $(\C, W, Cof(\C))$ is not an invariant of the classical categorical localization $\C[W^{-1}]$ : there are examples of pairs of Waldhausen categories with the same homotopy categories but with different $K$-theory spaces (see \cite{MR1930883}). The crutial results of Toën-Vezzosi in \cite{toenvezzosi-remarkonKtheory} allow us to identify the world of $(\infty,1)$-categories as the natural ultimate domain for $K$-theory. They prove that if the underlying $(\infty,1)$-categories associated to a pair of Waldhausen categories (via the $\infty$-localization) are equivalent then the associated $K$-theory spaces are equivalent. Moreover, in the same paper, the authors remark that the classical $S$-construction of Waldhausen can be lifted to the setting of $(\infty,1)$-categories. Following this insight, in \cite{1204.3607} the author introduces the notion of a \emph{Waldhausen $(\infty,1)$-category} (which, grosso modo are pairs of $(\infty,1)$-categories $(\A_0, \A)$ with $\A_0$ a full subcategory of $\A$ containing its maximal $\infty$-groupoid, together with extra conditions on this pair)  and develops this $\infty$-version of the $S$-construction. The collection of Waldhausen $(\infty,1)$-categories forms itself an $(\infty,1)$-category $\mathcal{W}ald_{\infty}$ and the result of this new $\infty$-version of the $S$-construction can be encoded as an $\infty$-functor $K^c_{Barwick}: \mathcal{W}ald_{\infty}\to \Sp_{\geq 0}$. Moreover,  there is a canonical $\infty$-functor linking the classical theory to this new approach 

\begin{equation}
\xymatrix{N(Wald_{Classic}^{Model})\ar[r]& \mathcal{W}ald_{\infty}}
\end{equation}

\noindent sending  a classical Waldhausen data $(\C, W, Cof(\C))$ to the $\infty$-localization $N(\C)[W^{-1}]$ together with its smallest subcategory containing the equivalences and the images of the cofibrations under the localization functor (see \cite[Example 2.12]{1204.3607}). The author then proves that the two $S$-constructions, respectively, the classical and the new $\infty$-version agree by means of this assignement and therefore produce the same $K$-theory (\cite[10.6.2]{1204.3607}). Up to our days this framework seems to be the most natural and general domain for connective $K$-theory. However, we should remark that a different $\infty$-categorical domain  has been established in the paper \cite{tabuada-gepner} where the authors study $K$-theory spaces associated to pointed $(\infty,1)$-categories having all finite colimits, whose collection forms an $(\infty,1)$-category $Cat_{\infty}(\omega)_{\ast}$. They generalize the classical $S$-construction to this new domain obtaining a new $\infty$-functor $K_{BGT}^{c}: Cat_{\infty}(\omega)_{\ast}\to \Sp_{\geq 0}$ 
,and prove that for any Waldhausen category $\C$  with equivalences $W$ (appearing as a subcategory of a model category), the $K$-theory space which their method assigns to the $\infty$-localization $N(\C)[W^{-1}]$  is equivalent to the classical $K$-theory space attached to $\C$ through the classical methods of Waldhausen. This framework is of course related to the wider framework of \cite{1204.3607}: following the Example  \cite[2.9]{1204.3607}, every pointed $(\infty,1)$-category with finite colimits has a naturally associated Waldhausen $(\infty,1)$-category. Again, this assignement can be properly understood as an $\infty$-functor

\begin{equation}
\xymatrix{\Psi: Cat_{\infty}(\omega)_{\ast}\ar[r]& \mathcal{W}ald_{\infty}}
\end{equation}

We summarize this fast historical briefing with the existence of a diagram of $(\infty,1)$-categories

\begin{equation}
\xymatrix{
&Cat_{\infty}(\omega)_{\ast}\ar@{..>}[dd]^(0.3){K^c_{B.G.T.}}\ar[dr]^{\Psi}&\\
N(Wald_{Classic}^{Model})\ar[ur]\ar[rr]\ar[dr]_{K^c_{Wald}}&&\ar[dl]^{K^c_{Barwick}}\mathcal{W}ald_{\infty}\\
&\Sp_{\geq 0}&
}
\end{equation}

\noindent whose commutativity follows from the results in \cite{1204.3607} and in \cite{tabuada-gepner} and from the agreement of the two $\infty$-categorical versions of the $S$-construction via $\Psi$. This agreement follows from the very definition of the two procedures. Consult  \cite[Section 5]{1204.3607} and \cite[Section 7.1]{tabuada-gepner} for the complete details.

\subsection{Connective $K$-theory of dg-categories}
\label{connectivektheorydg}

Our goal in this section is to explain how to define the connective $K$-theory of a dg-category and how to present this assignement as an $\infty$-functor $K^c:\dg^{idem}\to \widehat{\Sp}_{\geq 0}$ commuting with filtered colimits. One possible way is to use the classical theory of Waldhausen categories. As discussed in the Remark \ref{nc2strictmodel}, the data of an object $F\in Fun_{\omega}(\dg^{idem},\widehat{\Sp})$ corresponds in a essentially unique way to the data of an actual strict functor $F_s$ from the category of small dg-categories $Cat_{Ch(k)}$ with the Morita model structure to some combinatorial model category  whose underlying $(\infty,1)$-category is $\widehat{\Sp}$ such that $1)$ $F_s$ sends Morita equivalences to weak-equivalences and $2)$ $F_s$ preserves filtered homotopy colimits. 
In the case of connective $K$-theory such a functor can be obtained by composing the strict functor $K^c_{Wald}: Wald_{Classic}\to Sp^{\Sigma}$ of the previous section with the functor $Cat_{Ch(k)}\to Wald_{Classic}^{Model}$ defined by sending a small dg-category $T$ to the strict category of perfect cofibrant dg-modules (obtained by forgetting the dg-enrichement), with its natural structure of Waldhausen category given by the weak-equivalences of $T$-dg-modules and the cofibrations of the module structure therein. This is well-defined because perfect modules are stable under homotopy pushouts and satisfy the "cube lemma" \cite[5.2.6]{hovey-modelcategories}. The conditions $1)$ and $2)$ are also well-known to be satisfied (for instance see \cite[Section 2.2]{Anthony-thesis}). For the most part of this paper it will be enough to work with the $\infty$-functor $K^c: \dg^{ft}\to \widehat{\Sp}_{\geq 0}$ associated to this composition via the Remark \ref{nc2strictmodel} or its canonical $\omega$-continuous extension $\dg^{idem}\to \widehat{\Sp}$. However, some of our purposes (namely the Theorem \ref{teorema2}) will require an alternative approach. More precisely, and in the same spirit of \cite[Section 7.1]{tabuada-gepner} for stable $\infty$-categories, we will need to have a description of the Waldhausen's $S$-construction within the setting of dg-categories. \\

\begin{construction}
\label{Sconstruction1}
Let $Ar[n]_k$ be the dg-category obtained as the $k$-linearization of the category $Ar[n]$ described in the Construction \ref{Sconstruction0}. More precisely, its objects are the objects in $Ar[n]$ and its complexes of morphisms are all given by the ring $k$ seen as a complex concentrated in degree zero. For each $n$ the dg-category $Ar[n]_k$ is locally cofibrant (meaning, enriched over cofibrant complexes - see \refnci{dg1}) so that for any other locally cofibrant dg-category $T$ we have $Ar[n]_k\otimes^{\mathbb{L}}T \simeq Ar[n]_k\otimes T$ (recall our discussion in \refnci{dg1} about the derived monoidal structure in the $(\infty,1)$-category $\dg$ underlying the homotopy theory of dg-categories).

Recall also from \cite{Toen-homotopytheorydgcatsandderivedmoritaequivalences} that the symmetric monoidal $(\infty,1)$-category $\dg^{\otimes}$ admits an internal-hom $\mathbb{R}\underline{Hom}(A,B)$ given by the full sub-dg-category of right-quasi-representable cofibrant $A\otimes^{\mathbb{L}}B^{op}$-dg-modules. Iif $T$ is a locally cofibrant dg-category and $\widehat{T}_{c}$ is its idempotent-completion  (which as explained in the Remark \refnci{omega3} we can always assume to be locally cofibrant), we find a canonical equivalence in $\dg$ between $\mathbb{R}\underline{Hom}(A,\widehat{T}_{c})$ and $\widehat{A^{op}\otimes^{\mathbb{L}}T}_{pspe}$ - the full sub-dg-category of cofibrant pseudo-perfect $A\otimes^{\mathbb{L}}T^{op}$-dg-modules (by definition these are cofibrant dg-modules $E$ such that for any object $a\in A$, the $T^{op}$- module $E(A,-)$ is perfect). After the discussion in the previous paragraph we have $\mathbb{R}\underline{Hom}(Ar[n]_k, \widehat{T}_{c})\simeq  \widehat{Ar[n]_k^{op}\otimes^{\mathbb{L}}T}_{pspe} \simeq  \widehat{Ar[n]_k^{op}\otimes T}_{pspe}$ so that the objects in this internal-hom can be identified with $Ar[n]$-indexed diagrams in the underlying strict category of perfect cofibrant $T^{op}$-modules (obtained by forgetting the dg-enrichement).  We now set $S^{dg}_n(T)$ to be the full sub-dg-category of $\mathbb{R}\underline{Hom}(Ar[n]_k ,\widehat{T}_{c})$ spanned by those diagrams satisfying the conditions in the construction \ref{Sconstruction0}. These conditions make sense for the same reasons the functor $Cat_{Ch(k)}\to Wald_{Classic}^{Model}$ of the previous section also makes sense (see \cite[Section 2.2]{Anthony-thesis}). Again, the collection of dg-categories $S_n^{dg}(T)$ for $n\geq 0$ forms a simplicial object in dg-categories and by contemplating each level as a category (omitting its dg-enrichement) we can recover the $K$-theory of $T$ as $\Omega\, colim_{\Delta^{op}}N(S_n^{dg}(T)^{W_n})$ where $W_n$ is the class of maps in  $S^{dg}_{n}(T)$ given by the levelwise weak-equivalences of dg-modules and $S_n^{dg}(T)^{W_n}$ is the full subcategory of $S_n^{dg}(T)$ spanned by all the objects and only those morphisms which are in $W_n$.

Let now $[n]_k$ be the dg-category obtained as the $k$-linearization of the ordered category $[n]=\{0\leq 1\leq .... \leq n\}$. This dg-category is again locally-cofibrant and for the same reasons as above the underlying category obtained from $\mathbb{R}\underline{Hom}([n]_k, \widehat{T}_{c})$ by forgetting the dg-enrichement is the category of sequences of perfect cofibrant $T^{op}$-dg-modules of lenght $n+1$. As cofibers of maps are essentially uniquely determined up to isomorphism, we have a canonical equivalence of categories between $S^{dg}_{n}(T)$ and $\mathbb{R}\underline{Hom}([n-1]_k, \widehat{T}_{c})$. Since the model structure on $T^{op}$-dg-modules verifies the "cube lemma"  \cite[5.2.6]{hovey-modelcategories} (because $Ch(k)$ verifies for the projective model structure) this equivalence becomes an equivalence of pairs $(S^{dg}_{n}(T), W_n)$ and $(\mathbb{R}\underline{Hom}([n-1]_k, \widehat{T}_{c}), W_n')$ where we contemplate both dg-categories as categories by forgetting the dg-enrichements and where $W_n'$ denotes the class of maps of sequences which are levelwise given by weak-equivalences of dg-modules. Thanks to this equivalence we find an homotopy equivalence of simplicial sets between $N(S^{dg}_{n}(T)^{W_n})$ and $N(\mathbb{R}\underline{Hom}([n-1]_k, \widehat{T}_{c})^{ W_n'})$. Finally, and thanks to the main theorem of \cite{Toen-homotopytheorydgcatsandderivedmoritaequivalences} the last is exactly the mapping space $Map_{\dg}([n-1]_k, \widehat{T}_{c})$ which by adjunction is equivalent to $Map_{\dg^{idem}}(\widehat{([n-1]_k)}_c, \widehat{T}_{c})$. Under this chain of equivalences this family of mapping spaces for $n\geq 0$ inherits the structure of a simplicial object in the $(\infty,1)$-category of spaces and the K-theory space of $T$ can finally be rewritten as 

\begin{equation}
\label{nc2formulaktheory}
\Omega\, colim_{[n]\in\Delta^{op}}Map_{\dg^{idem}}(\widehat{([n-1]_k)}_c, \widehat{T}_{c})
\end{equation}

This concludes the construction.
\end{construction}

To conclude this section we remark two important properties of $K^c$. The first should be well-known to the reader:

\begin{prop}(\cite{waldhausen-ktheoryofspaces})
\label{waldausensfibration}
The $\infty$-functor $K^c:\dg^{idem}\to \widehat{\Sp}$ sends exact sequences of dg-categories to fiber sequences in $\widehat{\Sp}_{\geq 0}$.
\begin{proof}
This follows from the so called Waldhausen's Fibration Theorem \cite[1.6.4]{waldhausen-ktheoryofspaces} and  \cite[1.8.2]{thomasonalgebraic}, together with the dictionary between homotopy limits and homotopy colimits in the model category of spectra and limits and colimits in the $(\infty,1)$-category $\widehat{\Sp}$ (see \cite[1.3.4.23 and 1.3.4.24]{lurie-ha}).
\end{proof}
\end{prop}

The second  is a consequence of this first and will be very important to us:

\begin{prop}
\label{kcsemilocal}
$K^c$ sends Nisnevich squares of noncommutative smooth spaces to pullback squares of connective spectra.
\begin{proof}
Let

\begin{equation}
\label{nc215april1}
\xymatrix{
T_{\X}\ar[d]\ar[r]& T_{\UU}\ar[d]\\
T_{\V}\ar[r]& T_{\W}
}
\end{equation}

\noindent be a Nisnevich square of dg-categories. By definition, there are dg-categories $K_{\X-\UU}$ and $K_{\V-\W}$ in $\dg^{idem}$,  having compact generators, and such that the maps $T_{\X}\to T_{\UU}$ and $T_{\V}\to T_{\W}$ fit into strict short exact sequences in $\dg^{idem}$ (see \refnci{defopenimmersion} and the Remark \refnci{exactsequencesarestrict2})

\begin{equation}
\xymatrix{
K_{\X-\UU}\ar[d]\ar[r]& T_{\X}\ar[d]&&K_{\V-\UU}\ar[d]\ar[r]& T_{\V}\ar[d]\\
0\ar[r]& T_{\UU}&&0\ar[r]& T_{\W}
}
\end{equation}

Again by the definition of an open immersion and because of  \ref{waldausensfibration} we have pullback squares of connective spectra

\begin{equation}
\xymatrix{
K^c(K_{\X-\UU})\ar[d]\ar[r]& K^c(T_{\X})\ar[d]&&K^c(K_{\V-\UU})\ar[d]\ar[r]& K^c(T_{\V})\ar[d]\\
0\ar[r]& K^c(T_{\UU})&&0\ar[r]& K^c(T_{\W})
}
\end{equation}

With these properties in mind, we aim to show that the diagram

\begin{equation}
\label{nc215april2}
\xymatrix{
K^c(T_{\X})\ar[d]\ar[r]& K^c(T_{\UU})\ar[d]\\
K^c(T_{\V})\ar[r]& K^c(T_{\W})
}
\end{equation}

\noindent is a pullback of connective spectra. For that purpose we consider the pullback squares

\begin{equation}
\label{nc215april2}
\xymatrix{
K^c(K_{\X-\UU})\ar[d]\ar[r]& K^c(K_{\V-\W})\ar[d]\ar[r]& 0\ar[d]\\
K^c(T_{\X})\ar@{-->}[r]&K^c(T_{\V})\times_{K^c(T_{\W})}K^c(T_{\UU})\ar[d]\ar[r]& K^c(T_{\UU})\ar[d]\\
&K^c(T_{\V})\ar[r]& K^c(T_{\W})
}
\end{equation}

\noindent from which we extract a morphism of fiber sequences

\begin{equation}
\label{nc215april3}
\xymatrix{
K^c(K_{\X-\UU})\ar[d]\ar[r]& K^c(K_{\V-\W})\ar[d]\\
K^c(T_{\X})\ar@{-->}[r]\ar[d]&K^c(T_{\V})\times_{K^c(T_{\W})}K^c(T_{\UU})\ar[d]\\
K^c(T_{\UU})\ar@{=}[r]& K^c(T_{\UU})
}
\end{equation}

To conclude, since the square (\ref{nc215april1}) is Nisnevich, by definition, the canonical morphism $K_{\X-\UU}\to K_{\V-\W}$ is an equivalence in $\dg^{idem}$ so that the top map is an equivalence $K^c( K_{\X-\UU})\simeq K^c( K_{\V-\W})$. Using the associated long exact sequences we conclude that the canonical morphism $\xymatrix{K^c(T_{\X})\ar@{-->}[r]&K^c(T_{\V})\times_{K^c(T_{\W})}K^c(T_{\UU})}$ is also an equivalence, thus concluding the proof.

\end{proof}
\end{prop}

\subsection{Non-connective $K$-Theory}
\label{nonconnectivektheorydg}

The first attempts to define negative $K$-theory groups dates back to the works of Bass in \cite{MR0249491} and Karoubi in \cite{MR0233871}. The motivation to look for these groups is very simple: the higher $K$-theory groups of an exact sequence of Waldhausen categories do not fit in a long exact sequence. The full solution to this problem appeared in the legendary paper of Thomason-Trobaugh \cite{thomasonalgebraic} where the author provides a mechanism to extend the connective spectrum $K^C$ of Waldhausen to a new non-connective spectrum $K^B$ whose connective part recovers the classical data. His attention focuses on the $K$-theory of schemes and recovers the negative groups of Bass (by passing to the homotopy groups). Moreover, it satisfies the property people were waiting for \cite[Thm 7.4]{thomasonalgebraic}: for any reasonable scheme $X$ with an open subscheme $U\subseteq X$ with complementar $Z$, there is a pullback-pushout sequence of spectra $K(X \text{ on } Z)\to K(X)\to K(U)$ where $K(X \text{ on } Z)$ is the $K$-theory spectrum associated to the category of perfect complexes on $X$ supported on $Z$. Moreover, he proves that his non-connected version of $K$-theory satisfies descent with respect to the classical Nisnevich topology for schemes (see \cite[Thm 10.8]{thomasonalgebraic}).\\

More recently, Schlichting \cite{schlichting-negative} introduced a mechanism that allows us to define non-connective versions of $K$-theory in a wide range of situations and in \cite[Section 6 and 7]{tabuada-cisinski} the authors applied this algorithm to the context of dg-categories. The result is a procedure that sends Morita equivalences of dg-categories to weak-equivalences of spectra and commutes with filtered homotopy colimits (for instance, see \cite[2.12]{Anthony-thesis}) and comes canonically equipped with a natural transformation from connective $K$-theory inducing an equivalence in the connective part. By applying the arguments of the Remark \ref{nc2strictmodel} their construction can be encoded  in a unique way in the form of an $\omega$-continuous $\infty$-functor $K^S: \dg^{idem}\to \widehat{\Sp}$ together with a natural transformation $K^c\to K^S$ with $\tau_{\geq 0}K^c\simeq \tau_{\geq 0}K^S$. The \emph{motto} of non-connective $K$-theory can now be stated as

\begin{prop}
\label{schlichtinglocalization}
$K^S$ sends exact sequences in $\dg^{idem}$ to cofiber/fiber sequences in $\widehat{\Sp}$.
\begin{proof}
This follows from  \cite[12.1 Thm 9]{schlichting-negative} and from the adaptation of the Schlichting's setup to dg-categories in \cite[Section 6]{tabuada-cisinski}, together with the fact that our notion of exact sequences in $\dg^{idem}$ agrees with the notion of exact sequences in \cite{tabuada-cisinski} (see \refnci{kellerexact}-(3)). To conclude use again the dictionary between homotopy limits and homotopy colimits in a model category and limits and colimits on the underlying $(\infty,1)$-category.
\end{proof}
\end{prop}

Using the same arguments as in Prop. \ref{kcsemilocal}, we find 

\begin{cor}
\label{nc2negativektheoryisnisnevichlocal}
$K^S$ is Nisnevich local
\end{cor}

The method of Thomason (the so called $B$-construction) and the methods of Schlichting to create non-connective extensions of $K$-thery are somehow ad-hoc. In this paper we will show how these two constructions can both be understood as explicit models for the same process, namely, the Nisnevich "sheafification" \footnote{The noncommutative Nisnevich topology is not a Grothendieck topology.} in the noncommutative world.

\section{Proofs}

\subsection{Proof of the Theorem \ref{teorema1}: Non-connective $K$-theory is the Nisnevich localization of connective $K$-theory}
\label{proofteorema1}

As explained in the introduction,  the proof proceeds in two steps. First, in \ref{bassexactsequencefornisnevichlocal}, we prove that every Nisnevich local functor $F:\dg^{ft}\to \widehat{\Sp}$ satisfies the familiar Bass exact sequences for any integer $n$.  
The second step requires a more careful discussion. In \ref{seminisnevich} we introduce the notion of \emph{connective-nisnevich descent} for functors $F:\dg^{ft}\to \widehat{\Sp}$ with values in $\widehat{\Sp}_{\geq 0}$. We will see (Prop.\ref{connectivenisnevichislocalization} below) that the  full subcategory  $Fun_{Nis_{\geq 0}}(\dg^{ft}, \widehat{\Sp}_{\geq 0})$  spanned by those functors satisfying this descent conditon, is an accessible reflexive localization of $Fun(\dg^{ft}, \widehat{\Sp}_{\geq 0})$

\begin{equation}
\xymatrix{
Fun_{Nis_{\geq 0}}(\dg^{ft}, \widehat{\Sp}_{\geq 0})\ar@{^{(}->}[rr]^{\alpha} && \ar@/_2pc/[ll]_{l_{nis_{\geq 0}}} Fun(\dg^{ft}, \widehat{\Sp}_{\geq 0})
}
\end{equation}

\noindent and that as a consequence of the definition the connective truncation of a Nisnevich local is  connectively-Nisnevich local, and we have a natural factorization $\overline{\tau_{\geq 0}}$

\begin{equation}
\label{nc2rightadjointstruncation}
\xymatrix{
Fun(\dg^{ft}, \widehat{\Sp}_{\geq 0}) &\ar[l]^{\tau_{\geq 0}} Fun(\dg^{ft}, \widehat{\Sp})\\
Fun_{Nis_{\geq 0}}(\dg^{ft}, \widehat{\Sp}_{\geq 0})\ar@{^{(}->}[u]^{\alpha}&\ar@{-->}[l]_{\overline{\tau_{\geq 0}}} \ar@{^{(}->}[u]^{\beta}Fun_{Nis}(\dg^{ft}, \widehat{\Sp}) 
}
\end{equation}

\noindent where $\alpha$ and $\beta$ denote the inclusions. By abstract-nonsense, the composition $i_{!}:=l_{nis}^{nc}\circ i \circ \alpha$ provides a left adjoint to $\overline{\tau_{\geq 0}}$ and because the diagram of right adjoints commute, the diagram of left adjoints 

\begin{equation}
\label{nc2leftadjointstruncation}
\xymatrix{
Fun(\dg^{ft}, \widehat{\Sp}_{\geq 0})\ar[d]^{l_{nis_{\geq 0}}}\ar@{^{(}->}[r]^i & Fun(\dg^{ft}, \widehat{\Sp}) \ar[d]^{l_{nis}^{nc}}\\
Fun_{Nis_{\geq 0}}(\dg^{ft}, \widehat{\Sp}_{\geq 0})\ar@{-->}[r]^{i_!} & Fun_{Nis}(\dg^{ft}, \widehat{\Sp}) 
}
\end{equation}

\noindent also commutes.

The second step in our strategy amounts to check that the adjunction $(i_!, \overline{\tau_{\geq 0}})$ is an equivalence of $(\infty,1)$-categories. At this point our task is greatly simplified by the first step: the fact that Nisnevich local objects satisfy the Bass exact sequences for any integer $n$, implies that $\overline{\tau_{\geq 0}}$ is conservative. Therefore,  we are reduced to prove that the counit of the adjunction $\overline{\tau_{\geq 0}}\circ i_!\to Id$ is an equivalence of functors. In other words,  if $F$ is already connectively-Nisnevich local, its Nisnevich localization preserves the connective part. In order to achieve this we will need a more explicit description of the noncommutative Nisnevich localization of a connectively-Nisnevich local $F$. Our main result is that the more familiar $(-)^B$ construction of Thomason-Trobaugh (which we reformulate in our setting) provides such an explicit model, namely, we prove that if $F$ is connectively-Nisnevich local, $\tau_{\geq 0}(F^B)$ is naturally equivalent to $F$ and $F^B$ is Nisnevich local and naturally equivalent to $l_{nis}^{nc}(F)$.

\subsubsection{ Nisnevich descent forces all the Bass Exact Sequences }
\label{bassexactsequencefornisnevichlocal}
In this section we prove that every Nisnevich local $F: \dg^{ft}\to \widehat{\Sp}$ satisfies the familiar Bass exact sequences for any integer $n$. Our proof follows the arguments of \cite[6.1]{thomasonalgebraic}. The first step is to show that every Nisnevich local $F$ satisfies the Projective Bundle theorem. As explained in the introduction, this follows from the existence of an exceptional collection in $L_{pe}(\mathbb{P}^1)$ generated by the twisting sheaves $\mathcal{O}_{\mathbb{P}^1}$ and $\mathcal{O}_{\mathbb{P}^1}(-1)$, which, following \refnci{exceptional2}, provides a split short exact sequence of dg-categories

\begin{equation}
\label{nc2lola}
\xymatrix{
L_{pe}(k)\ar[d]\ar[r]_{i_{\mathcal{O}_{\mathbb{P}^1}}} & \ar@/_/[l] L_{pe}(\mathbb{P}^1)\ar[d]\\
0\ar[r]& L_{pe}(k) \ar@/_/[u]_{i_{\mathcal{O}_{\mathbb{P}^1}(-1)}}
}
\end{equation} 

\noindent where the map $i_{\mathcal{O}_{\mathbb{P}^1}}$, resp. $i_{\mathcal{O}_{\mathbb{P}^1}(-1)}$, is the inclusion of the full triangulated subcategory generated by $\mathcal{O}_{\mathbb{P}^1}$, respectively  $\mathcal{O}_{\mathbb{P}^1}(-1)$.  In particular, since $\dg^{ft}$ has direct sums,  we extract  canonical maps of dg-categories

\begin{eqnarray}
\label{nc2canonicalmap1}
\xymatrix{
 L_{pe}(k)\oplus L_{pe}(k)\ar[r]^{\psi} &L_{pe}(\mathbb{P}^1)&  L_{pe}(\mathbb{P}^1)\ar[r]^{\phi} &L_{pe}(k)\oplus L_{pe}(k)
}
\end{eqnarray}

We observe now that these maps become mutually inverse once we regard them in $Fun_{Nis}(\dg^{ft}, \widehat{\Sp})$ via the Yoneda's embedding. Indeed, the split exact sequence in (\ref{nc2lola}), or more precisely, its opposite in $\nck$, induces a split exact sequence in $Fun_{Nis}(\dg^{ft}, \widehat{\Sp})$

\begin{equation}
\xymatrix{
\Sigma^{\infty}_{+}\circ j_{nc}(L_{pe}(k))\ar[d]\ar[r] & \ar@/_/[l] \Sigma^{\infty}_{+}\circ j_{nc}(L_{pe}(\mathbb{P}^1))\ar[d]\\
0\ar[r]& \Sigma^{\infty}_{+}\circ j_{nc}(L_{pe}(k)) \ar@/_/[u]
}
\end{equation} 

This is because $Fun_{Nis}(\dg^{ft}, \widehat{\Sp})$  is stable (see the Remark \ref{allarestable2}), together with the effects of the Nisnevich localization. Again, we have canonical maps

\begin{equation}
\Sigma^{\infty}_{+}\circ j(L_{pe}(\mathbb{P}^1))\to \Sigma^{\infty}_{+}\circ j_{nc}(L_{pe}(k))\oplus \Sigma^{\infty}_{+}\circ j_{nc}(L_{pe}(k))
\end{equation}

\begin{equation}
 \Sigma^{\infty}_{+}\circ j_{nc}(L_{pe}(k))\oplus \Sigma^{\infty}_{+}\circ j_{nc}(L_{pe}(k))\to \Sigma^{\infty}_{+}\circ j(L_{pe}(\mathbb{P}^1))
\end{equation}

\noindent which, because in $\nck$ finite sums are the same as finite products (see the end of our discussion in \refnci{morita}), because the Yoneda's embedding commutes with finite products and because the pointing map $\Spaces\to \Spaces_{\ast}$ and the suspension $\Sigma^{\infty}$ commute with all colimits, can be identified with the image under $\Sigma^{\infty}_{+}\circ j_{nc}$ of the opposites of the canonical maps of dg-categories in (\ref{nc2canonicalmap1}), respectively.

This time, and as explained in \refnci{splitstablesums} and \refnci{exceptionalsplittosums}, because  $Fun_{Nis}(\dg^{ft}, \widehat{\Sp})$ is stable, these canonical maps are inverses to each other. In other words, we have a direct sum decomposition  

\begin{equation}
\Sigma^{\infty}_{+}\circ j(L_{pe}(\mathbb{P}^1))\simeq \Sigma^{\infty}_{+}\circ j_{nc}(L_{pe}(k))\oplus \Sigma^{\infty}_{+}\circ j_{nc}(L_{pe}(k))\simeq  \Sigma^{\infty}_{+}\circ j_{nc}(L_{pe}(k)\oplus L_{pe}(k))
\end{equation}

\noindent where the first (resp. second) component  can be identified with the part of $L_{pe}(\mathbb{P}^1)$ generated by $\mathcal{O}_{\mathbb{P}^1}$ (resp. $\mathcal{O}_{\mathbb{P}^1}(-1)$).\\

In particular, if we denote by \underline{Hom} the internal-hom in $Fun_{Nis}(\dg^{ft}, \widehat{\Sp})$ we find 

\begin{cor}
\label{nisnevichlocalprojectivebundle}
Let $F$ be a Nisnevich local functor $\dg^{ft}\to \widehat{\Sp}$. Then $F$ satisfies the projective bundle theorem. In other words, and following the above discussion, we have

\begin{eqnarray}
\underline{Hom}(\Sigma^{\infty}_{+}\circ j(L_{pe}(\mathbb{P}^1)), F)\simeq \underline{Hom}(\Sigma^{\infty}_{+}\circ j(L_{pe}(k)), F)\oplus \underline{Hom}(\Sigma^{\infty}_{+}\circ j(L_{pe}(k)), F) \simeq F\oplus F
\end{eqnarray}
\end{cor}

As in \cite[6.1]{thomasonalgebraic} we can now re-adapt this direct sum decomposition to a new one, suitably choosen to extract the Bass exact sequences out of the classical Zariski (therefore Nisnevich) covering of $\mathbb{P}^1$ given by

\begin{equation}
\label{nc2classicalcoveringP1}
\xymatrix{
\mathbb{G}_m \ar@{^{(}->}[r]^i \ar@{^{(}->}[d]^j& \mathbb{A}^1\ar@{^{(}->}[d]^{\alpha}\\
\mathbb{A}^1 \ar@{^{(}->}[r]^{\beta}& \mathbb{P}^1
}
\end{equation}

The basic ingredient is the induced pullback diagram of dg-categories

\begin{equation}
\label{bass1}
\xymatrix{
L_{pe}(\mathbb{P}^1)\ar[r]^{\alpha^{\ast}}\ar[d]^{\beta^{\ast}}& L_{pe}(\mathbb{A}^1)\ar[d]^{j^{\ast}}\\
 L_{pe}(\mathbb{A}^1)\ar[r]^{i^{\ast}} &   L_{pe}(\mathbb{G}_m)
}
\end{equation}

\noindent together with the composition

\begin{equation}
\label{nc2kka}
\xymatrix{
L_{pe}(k)\oplus L_{pe}(k)\ar[dr]^{\psi} \ar@/^/[drr]^{\alpha^{\ast}\circ \psi } \ar@/_/[ddr]_{\beta^{\ast}\circ \psi } &&\\
&L_{pe}(\mathbb{P}^1)\ar[r]^{\alpha^{\ast}}\ar[d]^{\beta^{\ast}}& L_{pe}(\mathbb{A}^1)\ar[d]^{j^{\ast}}\\
 &L_{pe}(\mathbb{A}^1)\ar[r]^{i^{\ast}} &   L_{pe}(\mathbb{G}_m)
}
\end{equation}

More precisely, we will focus on the diagram in $Fun(\dg^{ft}, \widehat{\Sp})$ induced by the opposite of the above diagram, namely,

\begin{equation}
\label{nc22april1}
\xymatrix{
\Sigma^{\infty}_{+}\circ j_{nc}(L_{pe}(\mathbb{G}_m)) \ar[r]^{L_{pe}(i)} \ar[d]^{L_{pe}(j)}& \Sigma^{\infty}_{+}\circ j_{nc}(L_{pe}(\mathbb{A}^1))\ar[d]^{L_{pe}(\alpha)}\ar@/^2pc/[ddr]^{}&\\
\Sigma^{\infty}_{+}\circ j_{nc}(L_{pe}(\mathbb{A}^1)) \ar@/_2pc/[drr]_(0.4){}\ar[r]^{L_{pe}(\beta)}&\Sigma^{\infty}_{+}\circ j_{nc}(L_{pe}(\mathbb{P}^1))\ar[dr]^{\Sigma^{\infty}_{+}\circ j_{nc}(\psi^{op})}&\\ 
&& \Sigma^{\infty}_{+}\circ j_{nc}(L_{pe}(k) \oplus  L_{pe}(k))
}
\end{equation}

\begin{remark}
\label{bass21}
It follows from \refnci{classicalnisnevichgoestoncnisnevich}, from the effects of the Nisnevich localization and from the above discussion that the exterior commutative square in (\ref{nc22april1}) becomes a pushout-pullback square in  $Fun_{Nis}(\dg^{ft}, \widehat{\Sp})$.
\end{remark}

In order to extract the Bass exact sequences, we consider a different direct sum decomposition of $\Sigma^{\infty}_{+}\circ j_{nc}(L_{pe}(\mathbb{P}^1))$. For that purpose let us start by introducing a bit of notation. We let $i_1$, $i_2$  denote the canonical inclusions $L_{pe}(k)\to L_{pe}(k)\oplus L_{pe}(k)$  in $\dg^{ft}$,  and let $\pi_1, \pi_2$ denote the projections $L_{pe}(k)\oplus L_{pe}(k)\to L_{pe}(k)$. At the same,  let $i_1^{op}$ and $i_2^{op}$ denote the associated projections in $\nck$ and $\pi_1^{op}$ and $\pi_2^{op}$ the canonical inclusions.  
Since the Yoneda's map $\Sigma^{\infty}_{+}\circ j_{nc}$ commutes with direct sums, the maps $\Sigma^{\infty}_{+}\circ j_{nc}(i_1^{op})$ and  $\Sigma^{\infty}_{+}\circ j_{nc}(i_1^{op})$ can be identified with the canonical projections 

\begin{equation}
\Sigma^{\infty}_{+}\circ j_{nc}(L_{pe}(k))\oplus \Sigma^{\infty}_{+}\circ j_{nc}(L_{pe}(k))\to  \Sigma^{\infty}_{+}\circ j_{nc}(L_{pe}(k))
\end{equation}

\noindent and  $\Sigma^{\infty}_{+}\circ j_{nc}(\pi_1^{op})$ and $\Sigma^{\infty}_{+}\circ j_{nc}(\pi_2^{op})$ the canonical inclusions

\begin{equation}
\Sigma^{\infty}_{+}\circ j_{nc}(L_{pe}(k))\to \Sigma^{\infty}_{+}\circ j_{nc}(L_{pe}(k)) \oplus \Sigma^{\infty}_{+}\circ j_{nc}( L_{pe}(k))
\end{equation}

\noindent  in $Fun(\dg^{ft}, \widehat{\Sp})$.

Let us proceed. To achieve the new decomposition, we compose the decomposition we had before with an equivalence $\Theta$ in  $Fun(\dg^{ft}, \widehat{\Sp})$

\begin{equation}
\xymatrix{
\Sigma^{\infty}_{+}\circ j_{nc}(L_{pe}(k))\oplus \Sigma^{\infty}_{+}\circ j_{nc}(L_{pe}(k))\ar[r]^{\Theta} &\Sigma^{\infty}_{+}\circ j_{nc}(L_{pe}(k))\oplus \Sigma^{\infty}_{+}\circ j_{nc}(L_{pe}(k))
}
\end{equation}

\noindent  defined to be  the map  

\begin{equation}
\xymatrix{
&\Sigma^{\infty}_{+}\circ j_{nc}(L_{pe}(k))\\
\Sigma^{\infty}_{+}\circ j_{nc}(L_{pe}(k))\oplus\Sigma^{\infty}_{+}\circ j_{nc}(L_{pe}(k))\ar[ur]^{\delta_1}\ar[dr]_{\delta_2}\ar@{-->}[r]^{\Theta}&\Sigma^{\infty}_{+}\circ j_{nc}(L_{pe}(k))\oplus\Sigma^{\infty}_{+}\circ j_{nc}(L_{pe}(k))\ar[u]_{\Sigma^{\infty}_{+}\circ j_{nc}(i_1^{op})}\ar[d]^{\Sigma^{\infty}_{+}\circ j_{nc}(i_2^{op})}\\
&\Sigma^{\infty}_{+}\circ j_{nc}(L_{pe}(k))
}
\end{equation}

\noindent obtained  from the universal property of the direct sum, where:

\begin{itemize}
\item $\delta_1$ it is the canonical dotted map obtained from the diagram

\begin{equation}
\xymatrix{
\Sigma^{\infty}_{+}\circ j_{nc}(L_{pe}(k))\ar[d]_{\Sigma^{\infty}_{+}\circ j_{nc}(\pi_1^{op})}\ar[dr]^{id}&\\
\Sigma^{\infty}_{+}\circ j_{nc}(L_{pe}(k))\oplus\Sigma^{\infty}_{+}\circ j_{nc}(L_{pe}(k))\ar@{-->}[r]&\Sigma^{\infty}_{+}\circ j_{nc}(L_{pe}(k))\\
\Sigma^{\infty}_{+}\circ j_{nc}(L_{pe}(k))\ar[u]^{\Sigma^{\infty}_{+}\circ j_{nc}(\pi_2^{op})}\ar[ur]^{0}& 
}
\end{equation}

\item $\delta_2$ is the canonical map obtained from

\begin{equation}
\xymatrix{
\Sigma^{\infty}_{+}\circ j_{nc}(L_{pe}(k))\ar[d]_{\Sigma^{\infty}_{+}\circ j_{nc}(\pi_1^{op})}\ar[dr]^{id}&\\
\Sigma^{\infty}_{+}\circ j_{nc}(L_{pe}(k))\oplus\Sigma^{\infty}_{+}\circ j_{nc}(L_{pe}(k))\ar@{-->}[r]&\Sigma^{\infty}_{+}\circ j_{nc}(L_{pe}(k))\\
\Sigma^{\infty}_{+}\circ j_{nc}(L_{pe}(k))\ar[u]^{\Sigma^{\infty}_{+}\circ j_{nc}(\pi_2^{op})}\ar[ur]^{-id}& 
}
\end{equation}
\end{itemize}

Of course, it follows from this definition that $\Theta$ is an equivalence with inverse equal to itself. Finally, we consider the composition

\begin{equation}
\label{bass24}
\xymatrix{
\Sigma^{\infty}_{+}\circ j_{nc}(L_{pe}(\mathbb{G}_m)) \ar[r]^{L_{pe}(i)} \ar[d]^{L_{pe}(j)}& \Sigma^{\infty}_{+}\circ j_{nc}(L_{pe}(\mathbb{A}^1))\ar[d]^{L_{pe}(\alpha)}\ar@/^2pc/[ddr]^{}&\\
\Sigma^{\infty}_{+}\circ j_{nc}(L_{pe}(\mathbb{A}^1)) \ar@/_2pc/[drr]_(0.4){}\ar[r]^{L_{pe}(\beta)}&\Sigma^{\infty}_{+}\circ j_{nc}(L_{pe}(\mathbb{P}^1))\ar[dr]^{\Theta\circ \Sigma^{\infty}_{+}\circ j_{nc}(\psi^{op})}&\\ 
&& \Sigma^{\infty}_{+}\circ j_{nc}(L_{pe}(k)) \oplus \Sigma^{\infty}_{+}\circ j_{nc}(L_{pe}(k)) 
}
\end{equation}

\noindent which again, as in the Remark \ref{bass21}, provides a pushout-pullback square in $Fun_{Nis}(\dg^{ft}, \widehat{\Sp})$. The important point of this new decomposition is the fact that both maps $\Theta\circ (\Sigma^{\infty}_{+}\circ j_{nc}(\psi^{op}\circ L_{pe}(\alpha)))$ and $\Theta\circ (\Sigma^{\infty}_{+}\circ j_{nc}(\psi^{op}\circ L_{pe}(\beta)))$  become simpler. In fact, since $\alpha^{\ast}(\mathcal{O}_{\mathbb{P}^1})= \alpha^{\ast}(\mathcal{O}_{\mathbb{P}^1}(-1))= \beta^{\ast}(\mathcal{O}_{\mathbb{P}^1})= \alpha^{\ast}(\mathcal{O}_{\mathbb{P}^1}(-1))= \mathcal{O}_{\mathbb{A}^1}$, we find that

\begin{itemize}
\item The composition $\Sigma^{\infty}_{+}\circ j_{nc}(i_1^{op})\circ\Theta\circ (\Sigma^{\infty}_{+}\circ j_{nc}(\psi^{op}\circ L_{pe}(\alpha)))$ can be identified with the map $\Sigma^{\infty}_{+}\circ j_{nc}(L_{pe}(p))$ induced by pullback along the canonical projection $p:\mathbb{A}^1\to Spec(k)$. Indeed, we have 

\begin{eqnarray}
\Sigma^{\infty}_{+}\circ j_{nc}(i_1^{op})\circ\Theta\circ (\Sigma^{\infty}_{+}\circ j_{nc}(\psi^{op}\circ L_{pe}(\alpha)))\simeq \\
\simeq \delta_1\circ (\Sigma^{\infty}_{+}\circ j_{nc}(\psi^{op}\circ L_{pe}(\alpha)))\simeq \\
\simeq \delta_1\circ (\Sigma^{\infty}_{+}\circ j_{nc}(\pi_1^{op}\circ i_1^{op} + \pi_2^{op}\circ i_2^{op}))\circ (\Sigma^{\infty}_{+}\circ j_{nc}(\psi^{op}\circ L_{pe}(\alpha)))\simeq\\
\simeq \Sigma^{\infty}_{+}\circ j_{nc}(i_1^{op}\circ \psi^{op}\circ L_{pe}(\alpha)) + 0 \simeq \\
\simeq \Sigma^{\infty}_{+}\circ j_{nc}((\alpha^{\ast}\circ \psi\circ i_1)^{op})\simeq\\
\simeq  \Sigma^{\infty}_{+}\circ j_{nc}((p^{\ast})^{op})\simeq\\
\simeq  \Sigma^{\infty}_{+}\circ j_{nc}(L_{pe}(p))
\end{eqnarray}

The same holds for the composition $\Sigma^{\infty}_{+}\circ j_{nc}(i_1^{op})\circ\Theta\circ (\Sigma^{\infty}_{+}\circ j_{nc}(\psi^{op}\circ L_{pe}(\beta)))$;

\item The maps $\Sigma^{\infty}_{+}\circ j_{nc}(i_2^{op})\circ\Theta\circ (\Sigma^{\infty}_{+}\circ j_{nc}(\psi^{op}\circ L_{pe}(\alpha)))$ and $\Sigma^{\infty}_{+}\circ j_{nc}(i_2^{op})\circ\Theta\circ (\Sigma^{\infty}_{+}\circ j_{nc}(\psi^{op}\circ L_{pe}(\beta)))$ are zero. Indeed, we have

\begin{eqnarray}
\Sigma^{\infty}_{+}\circ j_{nc}(i_2^{op})\circ\Theta\circ (\Sigma^{\infty}_{+}\circ j_{nc}(\psi^{op}\circ L_{pe}(\alpha)))\simeq\\  
\simeq \delta_2\circ (\Sigma^{\infty}_{+}\circ j_{nc}(\psi^{op}\circ L_{pe}(\alpha)))\simeq \\
\simeq \delta_2\circ (\Sigma^{\infty}_{+}\circ j_{nc}(\pi_1^{op}\circ i_1^{op} + \pi_2^{op}\circ i_2^{op}))\circ (\Sigma^{\infty}_{+}\circ j_{nc}(\psi^{op}\circ L_{pe}(\alpha)))\simeq\\
\simeq Id\circ(\Sigma^{\infty}_{+}\circ j_{nc}(i_1^{op}\circ \psi^{op}\circ L_{pe}(\alpha))) + (-Id)\circ (\Sigma^{\infty}_{+}\circ j_{nc}(i_1^{op}\circ \psi^{op}\circ L_{pe}(\alpha)))  \simeq \\
\simeq \Sigma^{\infty}_{+}\circ j_{nc}((\alpha^{\ast}\circ \psi\circ i_1)^{op}) - \Sigma^{\infty}_{+}\circ j_{nc}((\alpha^{\ast}\circ \psi\circ i_2)^{op})  
\end{eqnarray}

But since $\alpha^{\ast}(\mathcal{O}_{\mathbb{P}^1})= \alpha^{\ast}(\mathcal{O}_{\mathbb{P}^1}(-1))= \mathcal{O}_{\mathbb{A}^1}$, we have $\alpha^{\ast}\circ \psi\circ i_1\simeq \alpha^{\ast}\circ \psi\circ i_2$ so that the last  difference is zero. The same argument holds for $\beta^{\ast}$.
\end{itemize}

From these two facts, we conclude that $\Theta\circ (\Sigma^{\infty}_{+}\circ j_{nc}(\psi^{op}\circ L_{pe}(\alpha)))$ is equivalent to the sum $\Sigma^{\infty}_{+}\circ j_{nc}(L_{pe}(p))\oplus 0$ so that the outer commutative square of the diagram (\ref{bass24}) can now be written as

\begin{equation}
\label{bass1}
\xymatrix{
\Sigma^{\infty}_{+}\circ j_{nc}(L_{pe}(\mathbb{G}_m)) \ar[rrr]^{L_{pe}(i)} \ar[d]^{L_{pe}(j)}&&& \Sigma^{\infty}_{+}\circ j_{nc}(L_{pe}(\mathbb{A}^1))\ar[d]^{\Sigma^{\infty}_{+}\circ j_{nc}(L_{pe}(p))\oplus 0}\\
\Sigma^{\infty}_{+}\circ j_{nc}(L_{pe}(\mathbb{A}^1)) \ar[rrr]^(.4){\Sigma^{\infty}_{+}\circ j_{nc}(L_{pe}(p))\oplus 0}&&&\Sigma^{\infty}_{+}\circ j_{nc}(L_{pe}(k)) \oplus \Sigma^{\infty}_{+}\circ j_{nc}(L_{pe}(k)) 
}
\end{equation}

We are almost done.  To proceed, we rewrite the diagram \ref{bass24} as  

\begin{equation}
\label{nc2cafecafe0}
\xymatrix{
\Sigma^{\infty}_{+}\circ j_{nc}(L_{pe}(\mathbb{G}_m)) \ar[r] \ar[d]& \Sigma^{\infty}_{+}\circ j_{nc}(L_{pe}(\mathbb{A}^1))\oplus \Sigma^{\infty}_{+}\circ j_{nc}(L_{pe}(\mathbb{A}^1))\ar[d]_{(\Sigma^{\infty}_{+}\circ j_{nc}(L_{pe}(\alpha)),- \Sigma^{\infty}_{+}\circ j_{nc}(L_{pe}(\beta)))}\ar@/^2pc/[ddr]^(.8){(\Sigma^{\infty}_{+}\circ j_{nc}(L_{pe}(p)), -\Sigma^{\infty}_{+}\circ j_{nc}( L_{pe}(p)))\oplus 0}&\\
0 \ar@/_2pc/[drr]\ar[r]&\Sigma^{\infty}_{+}\circ j_{nc}(L_{pe}(\mathbb{P}^1))\ar[dr]^{\Sigma^{\infty}_{+}\circ j_{nc}(\Theta^{op}\circ\psi^{op})}&\\ 
&& \Sigma^{\infty}_{+}\circ j_{nc}(L_{pe}(k) \oplus  L_{pe}(k))
}
\end{equation}

\noindent where of course, since the Yoneda's map $\Sigma^{\infty}_{+}\circ j_{nc}$ commutes with direct sums, we have 

\begin{equation}
\Sigma^{\infty}_{+}\circ j_{nc}(L_{pe}(\mathbb{A}^1))\oplus \Sigma^{\infty}_{+}\circ j_{nc}(L_{pe}(\mathbb{A}^1))\simeq \Sigma^{\infty}_{+}\circ j_{nc}(L_{pe}(\mathbb{A}^1)\oplus L_{pe}(\mathbb{A}^1))
\end{equation}

We observe that both the inner and the outer squares become pullback-pushouts once we pass to the Nisnevich localization. Moreover, the map $\Theta\circ(\Sigma^{\infty}_{+}\circ j_{nc}(\psi^{op}))$ becomes an equivalence. 

In a different direction, we also observe that the pullback map of dg-categories $p^{\ast}: L_{pe}(k)\to L_{pe}(\mathbb{A}^1)$  admits a left inverse $s^{\ast}:L_{pe}(\mathbb{A}^1)\to  L_{pe}(k)$ given by the pullback along the zero section $s:Spec(k)\to \mathbb{A}^1 $\footnote{which in terms of rings is given by the evaluation at zero $ev_0:k[T]\to k$}. In terms of noncommutative spaces, this can be rephrased by saying that $L_{pe}(p)$ as a right inverse $L_{pe}(s)$. We can use this right-inverse to construct a right inverse to the first projection of $(\Sigma^{\infty}_{+}\circ j_{nc}(L_{pe}(p)), -\Sigma^{\infty}_{+}\circ j_{nc}( L_{pe}(p)))\oplus 0$, namely, we consider the map $(\Sigma^{\infty}_{+}\circ j_{nc}(L_{pe}(s)), 0)$ induced by the universal property of the  direct sum in $Fun(\dg^{ft}, \widehat{\Sp})$

\begin{equation}
\xymatrix{
&&&\Sigma^{\infty}_{+}\circ j_{nc}(L_{pe}(\mathbb{A}^1))\\
\Sigma^{\infty}_{+}\circ j_{nc}(L_{pe}(k))\ar@/^2pc/[urrr]^{\Sigma^{\infty}_{+}\circ j_{nc}(L_{pe}(s))}\ar@/_2pc/[drrr]_{0}\ar@{-->}[rrr]^(.4){(\Sigma^{\infty}_{+}\circ j_{nc}(L_{pe}(s)),0)}&&&\Sigma^{\infty}_{+}\circ j_{nc}(L_{pe}(\mathbb{A}^1))\oplus\Sigma^{\infty}_{+}\circ j_{nc}(L_{pe}(\mathbb{A}^1))\ar[u]\ar[d]\\
&&& \Sigma^{\infty}_{+}\circ j_{nc}(L_{pe}(\mathbb{A}^1))
}
\end{equation}

It is immediate to check that the composition $\Sigma^{\infty}_{+}\circ j_{nc}(i_1^{op})\circ( (\Sigma^{\infty}_{+}\circ j_{nc}(L_{pe}(p)), - \Sigma^{\infty}_{+}\circ j_{nc}(L_{pe}(p)))\oplus 0)\circ (\Sigma^{\infty}_{+}\circ j_{nc}((L_{pe}(s)),0)$ is the identity, so that $\Sigma^{\infty}_{+}\circ j_{nc}(i_1^{op})\circ( (\Sigma^{\infty}_{+}\circ j_{nc}(L_{pe}(p)), - \Sigma^{\infty}_{+}\circ j_{nc}(L_{pe}(p)))\oplus 0)$ has a right inverse that we can picture as a dotted arrow

\begin{equation}
\label{nc2cafecafe}
\xymatrix{
\Sigma^{\infty}_{+}\circ j_{nc}(L_{pe}(\mathbb{G}_m)) \ar[r] \ar[d]& \Sigma^{\infty}_{+}\circ j_{nc}(L_{pe}(\mathbb{A}^1))\oplus \Sigma^{\infty}_{+}\circ j_{nc}(L_{pe}(\mathbb{A}^1))\ar[d]^{(L_{pe}(\alpha),- L_{pe}(\beta))}\ar@/^10pc/[dd]^(0.6){(L_{pe}(p), - L_{pe}(p))\oplus 0}\\
0 \ar[r]&\Sigma^{\infty}_{+}\circ j_{nc}(L_{pe}(\mathbb{P}^1))\ar[d]^{\Theta\circ (\Sigma^{\infty}_{+}\circ j_{nc}(\psi^{op}))}\\ 
& \Sigma^{\infty}_{+}\circ j_{nc}(L_{pe}(k) \oplus  L_{pe}(k))\ar[d]^{\Sigma^{\infty}_{p}\circ j_{nc}(i_1^{op})}\\
&\Sigma^{\infty}_{+}\circ j_{nc}(L_{pe}(k))\ar@/_19pc/@{-->}[uuu]^{}
}
\end{equation}

At the same time, the preceding discussion implies  that the second projection

\begin{equation}
\label{nc2cafecafe2}
\xymatrix{
\Sigma^{\infty}_{+}\circ j_{nc}(L_{pe}(\mathbb{G}_m)) \ar[r] \ar[d]& \Sigma^{\infty}_{+}\circ j_{nc}(L_{pe}(\mathbb{A}^1))\oplus \Sigma^{\infty}_{+}\circ j_{nc}(L_{pe}(\mathbb{A}^1))\ar[d]^{(L_{pe}(\alpha),- L_{pe}(\beta))}\ar@/^15pc/[ddd]^(0.6){0}\\
0 \ar[r]&\Sigma^{\infty}_{+}\circ j_{nc}(L_{pe}(\mathbb{P}^1))\ar[d]^{\Theta\circ (\Sigma^{\infty}_{+}\circ j_{nc}(\psi^{op}))}\\ 
& \Sigma^{\infty}_{+}\circ j_{nc}(L_{pe}(k) \oplus  L_{pe}(k))\ar[d]^{\Sigma^{\infty}_{+}\circ j_{nc}(i_2^{op})}\\
&\Sigma^{\infty}_{+}\circ j_{nc}(L_{pe}(k))
}
\end{equation}

\noindent is just the zero map.

We now explain how to extract the familiar Bass exact sequence out of  these two diagrams. Given any object $F\in Fun(\dg^{ft}, \widehat{\Sp})$ and a noncommutative space $\X$, we set the notation $F_{\X}:= \underline{Hom}(\Sigma^{\infty}_{+}\circ j_{nc}(\X), F)$.  To proceed, we consider the image of the  diagram (\ref{nc2cafecafe}) under the functor $\underline{Hom}(-, F)$, to find a diagram

\begin{equation}
\xymatrix{
F_{L_{pe}(k)}\simeq F\ar[r]^{i^F_1} &F\oplus F\simeq F_{L_{pe}(k)\oplus L_{pe}(k)}\ar[r] & F_{L_{pe}(\mathbb{P}^1)}\ar[d]\ar[r] & F_{L_{pe}(\mathbb{A}^1)\oplus L_{pe}(\mathbb{A}^1) }\simeq F_{L_{pe}(\mathbb{A}^1)}\oplus F_{L_{pe}(\mathbb{A}^1)} \ar[d]\ar@/_2pc/[lll]\\
&&0\ar[r]&F_{L_{pe}(\mathbb{G}_m)}
}
\end{equation}

\noindent where the first map $i^F_1:F\to F\oplus F$ can be identified with the canonical inclusion in the first coordinate and  the composition 
$\xymatrix{F\ar[r]& F_{L_{pe}(\mathbb{A}^1)}\oplus  F_{L_{pe}(\mathbb{A}^1)}\ar@{-->}[r]& F}$ is the identity.

From this we can procude a new commutative diagram by taking sucessive pushouts

\begin{equation}
\label{nc2h1}
\xymatrix{
\ar[d] F\ar[r]^{i^F_1} &F\oplus F \ar[d]\ar[r]^{} & F_{L_{pe}(\mathbb{P}^1)}\ar[d]\ar[r] & F_{L_{pe}(\mathbb{A}^1)}\oplus F_{L_{pe}(\mathbb{A}^1)} \ar@/^2pc/[ddr]\ar[d]\\
0 \ar[r]& F\ar[r]^{}& F\coprod_{F\oplus F}F_{L_{pe}(\mathbb{P}^1)}\ar[d]\ar[r]&F_{L_{pe}(\mathbb{A}^1)}\coprod_{F} F_{L_{pe}(\mathbb{A}^1)}\ar@{-->}[dr] \\
&&0 \ar[rr]&&F_{L_{pe}(\mathbb{G}_m)}
}
\end{equation}

\noindent and we notice that the vertical map $F\oplus F\to F$ can be identified with the projection in the second coordinate.\\

In particular, if we denote as $U(F)$ the pullback 

\begin{equation}
\label{nc2h2}
\xymatrix{
U(F)\ar[r]\ar[d]& \ar[d] F_{L_{pe}(\mathbb{A}^1)}\coprod_{F} F_{L_{pe}(\mathbb{A}^1)}\\
0\ar[r]&F_{L_{pe}(\mathbb{G}_m)}
}
\end{equation}

\noindent we find a canonical map

\begin{equation}
\xymatrix{F\coprod_{F\oplus F}F_{L_{pe}(\mathbb{P}^1)}\ar@{-->}[r]& U(F)}
\end{equation}

\noindent induced from the diagram (\ref{nc2h1}) using the universal property of the pullback .\\

At the same time,  if we apply $\underline{Hom}(-,F)$ to the diagram (\ref{nc2cafecafe2}) we find a new commutative diagram

\begin{equation}
\label{nc2h2}
\xymatrix{
\ar@/_5pc/[dd]^{id}F\ar[d]^{i^F_2}\ar@/^2pc/[drrr]^0&&&\\
F\oplus F \ar[d]\ar[r]& F_{L_{pe}(\mathbb{P}^1)}\ar[d]\ar[rr] && F_{L_{pe}(\mathbb{A}^1)}\oplus F_{L_{pe}(\mathbb{A}^1)} \ar[d]\\
F\ar[r]& F\coprod_{F\oplus F}F_{L_{pe}(\mathbb{P}^1)}\ar[dd]\ar@{-->}[dr]\ar[rr]&&F_{L_{pe}(\mathbb{A}^1)}\coprod_{F} F_{L_{pe}(\mathbb{A}^1)}\ar[dd] \\
&&U(F)\ar[ur]\ar[dl]&\\
&0 \ar[rr]&&F_{L_{pe}(\mathbb{G}_m)}
}
\end{equation}

\noindent and discover that the map $F\to F\coprod_{F\oplus F}F_{L_{pe}(\mathbb{P}^1)}\to U(F)$ admits a natural factorization

\begin{equation}
\label{nc2h3}
\xymatrix{
&\Omega F_{L_{pe}(\mathbb{G}_m)}\ar[r]\ar[d]&0\ar[d]\\
F\ar@{-->}[ur]^{\sigma_F}\ar[r]&U(F)\ar[r]\ar[d]& F_{L_{pe}(\mathbb{A}^1)}\coprod_{F} F_{L_{pe}(\mathbb{A}^1)}\ar[d]\\
&0\ar[r]&F_{L_{pe}(\mathbb{G}_m)}
}
\end{equation}

\noindent because $\Omega F_{L_{pe}(\mathbb{G}_m)}$ is the fiber of $U(F)\to F_{L_{pe}(\mathbb{A}^1)}\coprod_{F} F_{L_{pe}(\mathbb{A}^1)}$. This concludes the preliminary steps.\\

From now, we suppose that $F$ is Nisnevich local. In this case, by the Corollary \ref{nisnevichlocalprojectivebundle}, the map $F\oplus F\to F_{L_{pe}(\mathbb{P}^1)}$ is an equivalence and the commutative square

\begin{equation}
\label{nc2h4}
\xymatrix{
F_{L_{pe}(\mathbb{P}^1)}\ar[d]\ar[r] &  F_{L_{pe}(\mathbb{A}^1)}\oplus F_{L_{pe}(\mathbb{A}^1)} \ar[d]\\
0\ar[r]&F_{L_{pe}(\mathbb{G}_m)}
}
\end{equation}

\noindent is a pushout-pullback because the image of the square (\ref{nc2classicalcoveringP1}) under $L_{pe}$ is a Nisnevich square of noncommutative spaces.  Using these two facts we conclude that the canonical maps constructed above, $F\to F\coprod_{F\oplus F}F_{L_{pe}(\mathbb{P}^1)}$ and $F\coprod_{F\oplus F}F_{L_{pe}(\mathbb{P}^1)}\to U(F)$ are equivalences so that the diagram

\begin{equation}
\label{nc2h5}
\xymatrix{
F\ar[r]\ar[d]& \ar[d] F_{L_{pe}(\mathbb{A}^1)}\coprod_{F} F_{L_{pe}(\mathbb{A}^1)}\\
0\ar[r]&F_{L_{pe}(\mathbb{G}_m)}
}
\end{equation}

\noindent is a pullback-pushout. In particular, as in  the diagram (\ref{nc2h3}) we find the existence of a section

\begin{equation}
\label{nc2h6}
\xymatrix{
&\Omega F_{L_{pe}(\mathbb{G}_m)}\ar[r]\ar[d]&0\ar[d]\\
F\ar@{-->}[ur]^{\sigma_F}\ar[r]^{Id}&F\ar[r]\ar[d]& F_{L_{pe}(\mathbb{A}^1)}\coprod_{F} F_{L_{pe}(\mathbb{A}^1)}\ar[d]\\
&0\ar[r]&F_{L_{pe}(\mathbb{G}_m)}
}
\end{equation}

We are almost done. To conclude, we consider the induced pullback-pushout square

\begin{equation}
\label{nc2h7}
\xymatrix{
F_{L_{pe}(\mathbb{A}^1)}\coprod_{F} F_{L_{pe}(\mathbb{A}^1)}\ar[r]\ar[d]& 0\ar[d]\\
F_{L_{pe}(\mathbb{G}_m)}\ar[r]& \Sigma F
}
\end{equation}

\noindent where now, the suspension $\Sigma(\sigma_F)$ makes $\Sigma(F)$ a retract of $F_{L_{pe}(\mathbb{G}_m)}$. We are done now. Since the evaluation maps commute with colimits and by the definition of $F_{(-)}$, we have  for each $T_{\X} \in \dg^{ft}$, a pullback-pushout diagram in $\widehat{\Sp}$

\begin{equation}
\label{nc2h8}
\xymatrix{
F(L_{pe}(\mathbb{A}^1)\otimes T_{\X} )\coprod_{F(T_{\X} )} F(L_{pe}(\mathbb{A}^1)\otimes T_{\X} ) \ar[rr]\ar[d]&& 0\ar[d]\\
F(L_{pe}(\mathbb{G}_m)\otimes T_{\X} )\ar[rr]&& \Sigma F(T_{\X} )
}
\end{equation}

\noindent and therefore, a long exact sequence of abelian groups

\begin{equation}
\label{nc2h9}
...\to \pi_n(F(L_{pe}(\mathbb{A}^1)\otimes T_{\X} )\coprod_{F( T_{\X} )} F(L_{pe}(\mathbb{A}^1)\otimes T_{\X} ))\to \pi_n(F(L_{pe}(\mathbb{G}_m)\otimes T_{\X} ))\to \pi_n(\Sigma F(T_{\X} ))= \pi_{n-1}(F(T_{\X} ))\to ...
\end{equation}

\noindent and because of the existence of $\Sigma(\sigma_F)$, the maps $ \pi_n(F(L_{pe}(\mathbb{G}_m)\otimes T_{\X} ))\to \pi_n(\Sigma F(T_{\X} ))= \pi_{n-1}(F(T_{\X} ))$ are necessarily surjective, so that the long exact sequence breaks up into short exact sequences

\begin{equation}
\label{nc2h10}
0\to \pi_n(F(L_{pe}(\mathbb{A}^1)\otimes T_{\X} )\coprod_{F(T_{\X} )} F(L_{pe}(\mathbb{A}^1)\otimes T_{\X} ))\to \pi_n(F(L_{pe}(\mathbb{G}_m)\otimes T_{\X} ))\to \pi_n(\Sigma F( T_{\X} ))= \pi_{n-1}(F( T_{\X} ))\to 0
\end{equation}

\noindent $\forall n\in \mathbb{Z}$.

At the same time, since the square

\begin{equation}
\label{nc2h11}
\xymatrix{
F\ar[d]\ar[r]^{i^F_1} & F\oplus F\ar[r]&  F_{L_{pe}(\mathbb{A}^1)}\oplus F_{L_{pe}(\mathbb{A}^1)}\ar[d]\\
0 \ar[rr] && F_{L_{pe}(\mathbb{A}^1)}\coprod_{F} F_{L_{pe}(\mathbb{A}^1)}
}
\end{equation}

\noindent is also a pullback-pushout and the top map $F\to F_{L_{pe}(\mathbb{A}^1)}\oplus F_{L_{pe}(\mathbb{A}^1)}$ admits a left inverse, the associated long exact sequences

\begin{eqnarray}
\label{nc2h12}
...\to \pi_n(F(L_{pe}(\mathbb{A}^1)\otimes T_{\X} )\coprod_{F(T_{\X} )} F(L_{pe}(\mathbb{A}^1)\otimes T_{\X} ))\to  \pi_n(F(T_{\X} ))\to \pi_n(F(L_{pe}(\mathbb{A}^1)\otimes T_{\X} )\oplus F(L_{pe}(\mathbb{A}^1)\otimes T_{\X} ))\to \\
\to \pi_n(F(L_{pe}(\mathbb{A}^1)\otimes T_{\X} )\coprod_{F(\X)} F(L_{pe}(\mathbb{A}^1)\otimes T_{\X} ))\to ...
\end{eqnarray}

\noindent breaks up into short exact sequences

\begin{eqnarray}
\label{nc2h13}
0 \to  \pi_n(F(T_{\X} ))\to \pi_n(F(L_{pe}(\mathbb{A}^1)\otimes T_{\X} )\oplus F(L_{pe}(\mathbb{A}^1)\otimes T_{\X} ))\to \pi_n(F(L_{pe}(\mathbb{A}^1)\otimes T_{\X} )\coprod_{F(\X)} F(L_{pe}(\mathbb{A}^1)\otimes T_{\X} ))\to 0.
\end{eqnarray}

Combining the two short exact sequences (\ref{nc2h10}) and (\ref{nc2h13}) we find the familar exact sequences of Bass-Thomason-Trobaugh

\begin{eqnarray}
\label{nc2h13}
0 \to  \pi_n(F(T_{\X} ))\to \pi_n(F(L_{pe}(\mathbb{A}^1)\otimes T_{\X} )\oplus F(L_{pe}(\mathbb{A}^1)\otimes T_{\X} ))\to\\
\to  \pi_n(F(L_{pe}(\mathbb{G}_m)\otimes T_{\X} ))\to  \pi_{n-1}(F(T_{\X} ))\to 0
\end{eqnarray}

This concludes this section.

\subsubsection{Nisnevich vs Connective-Nisnevich descent and the Thomason-Trobaugh $(-)^B$-Construction }
\label{seminisnevich}

In this section we study the class of functors sharing the same formal properties of $K^c$, namely, the one of sending Nisnevich squares to pullback squares of connective spectra. This will take us through a small digression aiming to understand how the truncation functor $\tau_{\geq 0}$ interacts with the Nisnevich localization.

\begin{defn}
Let $F\in Fun(\dg^{ft}, \widehat{\Sp}_{\geq 0})$. We say that $F$ is \emph{connectively-Nisnevich local } if for any Nisnevich square of dg-categories

\begin{equation}
\label{nc2hojeestaquase}
\xymatrix{
T_{\X}\ar[r]\ar[d]&T_{\UU}\ar[d]\\
T_{\V}\ar[r] &T_{\W}
}
\end{equation}

\noindent the induced square 

\begin{equation}
\label{nc2hojeestaquase2}
\xymatrix{
F(T_{\X})\ar[r]\ar[d]&F(T_{\UU})\ar[d]\\
F(T_{\V})\ar[r] &F(T_{\W})
}
\end{equation}

\noindent is a pullback of connective spectra. 
\end{defn}

\begin{remark}
\label{nc2truncationnisnevichisseminisnevich}
It follows that if $F$ belongs to $Fun_{Nis}(Dg^{ft}, \widehat{\Sp})$, its connective truncation $\tau_{\geq 0}(F)$ is connectively-Nisnevich local. This is because $\tau_{\geq 0}$ acts objectwise and is a right adjoint to the inclusion of connective spectra into all spectra, thus preserving pullbacks.
\end{remark}

It is also convenient to isolate the following small technical remark:

\begin{remark}
\label{nc2remarksplit}
Let $\C$ be a stable $(\infty,1)$-category  and let  $\C_{0}\subseteq \C$ be a subcategory such that the inclusion  preserves direct sums. Then, if 

\begin{equation}
\label{nc2ficha}
\xymatrix{
A\ar[r]^i\ar[d]& B\ar[d]^p\\
0\ar[r] & C
}
\end{equation}

\noindent is a pullback square in $\C_{0}$ such that 

\begin{itemize}
\item the map $i$ admits a left inverse $v$;
\item the map $p$ admits a right inverse $u$;
\item the sum $i\circ v+u\circ p$ is homotopic to the identity
\end{itemize}

\noindent we conclude, by the same arguments given in the Remark \refnci{splitstablesums} that $B\simeq A\oplus C$. Moreover, under the hypothesis that the inclusion preserves direct sums, the square \label{nc2ficha} remains a pullback after the inclusion $\C_{0}\subseteq \C$ and therefore, a pushout. In particular, it becomes a split exact sequence in $\C$. This holds for any universe.

In particular,  for any pullback square of dg-categories associated to a Nisnevich square of noncommutative spaces (\ref{nc2hojeestaquase}) such that $T_{\UU}$ is zero and the sequence splits, the induced diagram of connective spectra (\ref{nc2hojeestaquase2}) makes  $F(T_{\V})$ canonically equivalent to the direct sum $F(T_{\X})\oplus F(T_{\W})$ in $\widehat{\Sp}$.   
\end{remark}

We let $Fun_{Nis_{\geq 0}}(\dg^{ft}, \widehat{\Sp}_{\geq 0})$ denote the full subcategory of $Fun(\dg^{ft}, \widehat{\Sp}_{\geq 0})$  spanned by the connectively-Nisnevich local functors. For technical reasons it is convenient to observe that the inclusion $Fun_{Nis_{\geq 0}}(\dg^{ft}, \widehat{\Sp}_{\geq 0})\subseteq Fun(\dg^{ft}, \widehat{\Sp}_{\geq 0})$ admits a left adjoint $l_{nis_{\geq 0}}$. More precisely

\begin{prop}
\label{connectivenisnevichislocalization}
$Fun_{Nis_{\geq 0}}(\dg^{ft}, \widehat{\Sp}_{\geq 0})$ is an accessible reflexive localization of $Fun(\dg^{ft}, \widehat{\Sp}_{\geq 0})$.
\begin{proof}
We evoke the Proposition 5.5.4.15 of \cite{lurie-htt} so that we are reduced to show the existence of a small class of maps $S$ in $Fun(\dg^{ft}, \widehat{\Sp}_{\geq 0})$ such that an object $F$ is connectively-Nisnevich local if and only if it is local with respect to the maps in $S$.

To define $S$, we ask the reader to bring back to his attention our discussion and notations in \refnci{usingpresheavesofspectra1} and in \refnci{usingpresheavesofspectra2}. Using the same notations, we define $S$ to be the collection of all maps 

\begin{equation}
 \delta_{\Sigma^{\infty}_{+}\circ j_{nc}(\UU)}(K)\coprod_{\delta_{ \Sigma^{\infty}_{+}\circ j_{nc}(\W)}(K)}\delta_{ \Sigma^{\infty}_{+}\circ j_{nc}(\V)}(K)\to  \delta_{ \Sigma^{\infty}_{+}\circ j_{nc}(\X)}(K)
\end{equation}

\noindent given by the universal property of the pushout, this time with $K$ in $\widehat{\Sp}_{\geq 0}\cap (\widehat{\Sp})^{\omega}$ \footnote{Here $(\widehat{\Sp})^{\omega}$ denotes the full subcategory of $\widehat{\Sp}$ spanned by the compact objects. Recall that $\widehat{\Sp}\simeq Ind((\widehat{\Sp})^{\omega})$.} and $\W$,$V$,$\UU$ and $\X$ part of a Nisnevich square of noncommutative smooth spaces. As before, the fact that $S$ satisfies the required property follows directly from the definition of the functors $\delta_{ \Sigma^{\infty}_{+}\circ j_{nc}(-)}$ as left adjoints to $Map^{Sp}$ and from the enriched version of Yoneda's lemma.

\end{proof}
\end{prop}

It follows directly from the definition of the class $S$ in the previous proof and from the description of the class of maps that generate the Nisnevich localization in $Fun(\dg^{ft}, \widehat{\Sp})$ (see \refnci{usingpresheavesofspectra2}) that the inclusion 

\begin{equation}i:Fun(\dg^{ft}, \widehat{\Sp}_{\geq 0})\hookrightarrow Fun(\dg^{ft}, \widehat{\Sp})\end{equation}

\noindent sends connective-Nisnevich local equivalences to Nisnevich local equivalences. In particular,  the universal property of the localization  provides us with a canonical colimit preserving map 

\begin{equation}
\xymatrix{
Fun(\dg^{ft}, \widehat{\Sp}_{\geq 0})\ar[d]^{l_{nis_{\geq 0}}}\ar@{^{(}->}[r]^i & Fun(\dg^{ft}, \widehat{\Sp}) \ar[d]^{l_{nis}^{nc}}\\
Fun_{Nis_{\geq 0}}(\dg^{ft}, \widehat{\Sp}_{\geq 0})\ar@{-->}[r] & Fun_{Nis}(\dg^{ft}, \widehat{\Sp}) 
}
\end{equation}

\noindent rendering the diagram commutative. Moreover, since the localizations are presentable, the Adjoint Functor Theorem implies the existence of a right adjoint which makes the associated diagram of right adjoints 

\begin{equation}
\label{nc2rightadjointstruncation2}
\xymatrix{
Fun(\dg^{ft}, \widehat{\Sp}_{\geq 0}) &\ar[l]^{\tau_{\geq 0}} Fun(\dg^{ft}, \widehat{\Sp})\\
Fun_{Nis_{\geq 0}}(\dg^{ft}, \widehat{\Sp}_{\geq 0})\ar@{^{(}->}[u]^{\alpha}&\ar@{-->}[l] \ar@{^{(}->}[u]^{\beta}Fun_{Nis}(\dg^{ft}, \widehat{\Sp}) 
}
\end{equation}

\noindent commute. At the same time, the Remark \ref{nc2truncationnisnevichisseminisnevich} implies the existence of the two commutative diagrams  (\ref{nc2rightadjointstruncation2}) and (\ref{nc2rightadjointstruncation2}).  By comparison with the new diagrams, we find that the canonical colimit preserving map  $\xymatrix{Fun_{Nis_{\geq 0}}(\dg^{ft}, \widehat{\Sp}_{\geq 0})\ar@{-->}[r]& Fun_{Nis}(\dg^{ft}, \widehat{\Sp})}$  can be identified with the composition $i_{!}:=l_{Nis}^{nc}\circ i\circ \alpha$ and that its right adjoint can be identified with  $\overline{\tau_{\geq 0}}$,  the restriction of the truncation functor $\tau_{\geq 0}$ to the Nisnevich local functors. 

Our goal is to prove that this adjunction

\begin{equation}
\xymatrix{
Fun_{Nis_{\geq 0}}(\dg^{ft}, \widehat{\Sp}_{\geq 0})\ar[r]_{i_!}&\ar@/^1pc/[l]^{\overline{\tau_{\geq 0}}}Fun_{Nis}(\dg^{ft}, \widehat{\Sp}) 
}
\end{equation}

\noindent is an equivalence. Our results from \ref{bassexactsequencefornisnevichlocal}  already provide one step towards this:

\begin{prop}
\label{nc2301}
The functor $\overline{\tau_{\geq 0}}$ is conservative.
\begin{proof}
Recall from  \ref{bassexactsequencefornisnevichlocal} that for any Nisnevich local $F$ we can construct a pullback-pushout square

\begin{equation}
\xymatrix{
F_{L_{pe}(\mathbb{A}^1)}\coprod_{F} F_{L_{pe}(\mathbb{A}^1)}\ar[r]\ar[d]& 0\ar[d]\\
F_{L_{pe}(\mathbb{G}_m)}\ar[r]& \Sigma F
}
\end{equation}

\noindent such that for any $T_{\X} \in \dg^{ft}$, the associated long exact sequence of homotopy groups breaks up into short exact sequences for any $n\in \mathbb{N}$

\begin{equation}
0\to \pi_n(F(L_{pe}(\mathbb{A}^1)\otimes T_{\X} )\coprod_{F(T_{\X} )} F(L_{pe}(\mathbb{A}^1)\otimes T_{\X} ))\to \pi_n(F(L_{pe}(\mathbb{G}_m)\otimes T_{\X} ))\to \pi_{n-1}(F(T_{\X} ))\to 0
\end{equation}

Therefore, given a morphism $F\to G$ in $Fun_{Nis}(\dg^{ft}, \widehat{\Sp})$, we have an induced diagram

\begin{equation}
\xymatrix{
&G_{L_{pe}(\mathbb{A}^1)}\coprod_{G} G_{L_{pe}(\mathbb{A}^1)}\ar[rr]\ar[dd]&&0\ar[dd]\\
\ar[ru]\ar[dd]\ar[rr]F_{L_{pe}(\mathbb{A}^1)}\coprod_{F} F_{L_{pe}(\mathbb{A}^1)}&&\ar[dd]0\ar[ur]&\\
&G_{L_{pe}(\mathbb{G}_m)}\ar[rr]&&\Sigma G\\
F_{L_{pe}(\mathbb{G}_m)}\ar[ru]\ar[rr]&&\Sigma F\ar[ru]&
}
\end{equation}

\noindent which induces natural maps of short exact sequences

\begin{equation}
\xymatrix{
0\ar[r] \ar[d]& \ar[d]\pi_n(F(L_{pe}(\mathbb{A}^1)\otimes T_{\X} )\coprod_{F(T_{\X} )} F(L_{pe}(\mathbb{A}^1)\otimes T_{\X} ))\ar[r]& \pi_n(F(L_{pe}(\mathbb{G}_m)\otimes T_{\X} ))\ar[r]\ar[d] &\pi_{n-1}(F( T_{\X} ))\ar[r] \ar[d]& 0\ar[d]\\
0\ar[r] & \pi_n(G(L_{pe}(\mathbb{A}^1)\otimes T_{\X} )\coprod_{G(T_{\X} )} G(L_{pe}(\mathbb{A}^1)\otimes T_{\X} ))\ar[r]& \pi_n(G(L_{pe}(\mathbb{G}_m)\otimes T_{\X} ))\ar[r]&\pi_{n-1}(G(T_{\X} ))\ar[r] & 0
}
\end{equation}

In particular, if $f$ is an equivalence in the connective part, by induction on $n=0,-1, -2,...$, we conclude that $f$ is an equivalence.
\end{proof}
\end{prop}

With this result, in order to prove that $i_!$ is an equivalence we are reduced to show the counit of the adjunction $\overline{\tau_{\geq 0}}\circ i_!\to Id$ is a natural equivalence of functors. Notice that since $\alpha$ and $i$ are fully-faithful, this amounts to show that for any $F$ connectively-Nisnevich local, the canonical map $i\circ \tau_{\geq 0}\circ l_{Nis}^{nc}\circ i\circ \alpha (F)\to i\circ \alpha (F)$ is an equivalence. Of course, to achieve this we will need a more explicit description of the noncommutative Nisnevich localization functor $l_{Nis}^{nc}$ restricted to connectively-Nisnevich local objects. There is a naive candidate, namely, the familiar $(-)^B$ construction of Thomason-Trobaugh \cite[6.4]{thomasonalgebraic}. Our goal to the end of this section is to prove the following proposition  confirming that this guess is correct:

\begin{prop}
\label{nc2proposicao1}
There is an accessible localization functor $(-)^B: Fun(\dg^{ft}, \widehat{\Sp})\to Fun(\dg^{ft}, \widehat{\Sp})$ encoding the $B$-construction of \cite[6.4]{thomasonalgebraic} such that for any $F\in Fun_{Nis_{\geq 0}}(\dg^{ft}, \widehat{\Sp})$ we have:

\begin{itemize}
\item  $\overline{\tau_{\geq 0}}(i\circ \alpha (F)^{B})\simeq F$.
\item the object  $(i\circ \alpha (F))^{B}$ is Nisnevich local;
\item there is a canonical equivalence $(i\circ \alpha (F))^{B}\simeq l_{Nis}^{nc}((i\circ \alpha (F))$;

\end{itemize}
In particular, the natural transformation $\overline{\tau_{\geq 0}}\circ i_!\to Id$ is an equivalence. Together with the Proposition \ref{nc2301} we have an equivalence of $(\infty,1)$-categories between the theory of connectively-Nisnevich local functors and the theory of Nisnevich local functors.
\end{prop}

With these results available we can already uncover the proof of our first main theorem:\\

\textit{Proof of the Theorem \ref{teorema1}:}

Thanks to the  Corollary \ref{nc2negativektheoryisnisnevichlocal} we already know that $K^S$ is Nisnevich local. In this case, and by the universal property of the localization, the canonical map $K^c\to K^S$ admits a canonical uniquely determined factorization

\begin{equation}
\label{nc2presta}
\xymatrix{
K^c\ar[d]\ar[dr]&\\
l_{Nis}^{nc}(K^c)\ar@{-->}[r]& K^S
}
\end{equation}

\noindent so that we are reduced to show that this canonical morphism $l^{nc}_{Nis}(K^c)\to K^S$ is an equivalence. But since these are Nisnevich local objects and since we now know by the Prop. \ref{nc2proposicao1} that the truncation functor $\tau_{\geq 0}$ is an equivalence when restricted to Nisnevich locals, it suffices to check that the induced map $\tau_{\geq 0} l^{nc}_{Nis}(K^c)\to \tau_{\geq 0} K^S$ is an equivalence. But this follows because all the morphisms in the image of the commutative diagram (\ref{nc2presta}) become equivalences after applying $\tau_{\geq 0}$. This follows from the construction of $K^S$ and again by the results in the Proposition \ref{nc2proposicao1}. 

\hfill $\qed$

We now start our small journey towards the proof of the Proposition \ref{nc2proposicao1}. To start with we need to specify how the $B$-Construction of \cite[6.4]{thomasonalgebraic} can be formulated in our setting:

\begin{construction}
\label{Bconstruction}
(Thomason-Trobaugh  $(-)^B$-Construction) We begin by asking the reader to bring back to his attention the diagrams constructed in \ref{bassexactsequencefornisnevichlocal}. More precisely, to recall that for any $F\in Fun(\dg^{ft}, \widehat{\Sp})$, we found a commutative diagram

\begin{equation}
\label{nc2hkilla}
\xymatrix{
&\Omega F_{L_{pe}(\mathbb{G}_m)}\ar[r]\ar[d]&0\ar[d]\\
F\ar@{-->}[ur]^{\sigma_F}\ar[r]^{\alpha_F}&U(F)\ar[r]\ar[d]& F_{L_{pe}(\mathbb{A}^1)}\coprod_{F} F_{L_{pe}(\mathbb{A}^1)}\ar[d]\\
&0\ar[r]&F_{L_{pe}(\mathbb{G}_m)}
}
\end{equation}

\noindent where both squares are pushout-pullbacks. Iterating this construction, we  find a sequence of canonical maps

\begin{equation}
\label{nc2chinesmacaco}
\xymatrix{
F\ar[rr]^{\alpha_F} && U(F)\ar[rr]^{\alpha_{U(F)}} &&U(U(F))\ar[rr]^{\alpha_{U^2(F)}}&& ....
}
\end{equation}

\noindent and we define $F^B$ to be the colimit for sequence (which is of course unique up to canonical equivalence). The assignement $F\mapsto F^B$ provides an endofunctor $(-)^B$ of the $(\infty,1)$-category $Fun(\dg^{ft}, \widehat{\Sp})$. To see this we can use the fact the monoidal structure in $Fun(\dg^{ft}, \widehat{\Sp})$ admits internal-homs $\underline{Hom}$. More precisely, we consider the diagram of natural transformations induced by the image of the diagram (\ref{nc2cafecafe}) under the first entry of $Hom(-,-)$. With this, and keeping the notations we have been using, we define $f_1$ to be the functor cofiber of $Id= (-)_{L_{pe}(k)}\to (-)_{L_{pe}(\mathbb{A}^1)}\oplus  (-)_{L_{pe}(\mathbb{A}^1)}$. The universal property of the cofiber gives us a canonical natural transformation $f_1\to (-)_{L_{pe}(\mathbb{G}_m)}$ and define a new functor $U$ as the fiber of  this map (recall that colimits and limits in the category of functors are determined objectwise). Finally, we consider $(-)^B$ as the colimit of the natural transformations 

\begin{equation}
\xymatrix{
Id \ar[dr]\ar[r]^{\alpha}& \ar[d]U=Id\circ U\ar[r]& \ar[dl]U^2=Id \circ U^2\ar[r]& ....\\
& (-)^B&&
}
\end{equation}

\end{construction}

We prove that for any $F$ the object $F^B$ satisfies the exact sequences of Bass-Thomason-Trobaugh for any $n\in \mathbb{Z}$. The proof requires some technical steps:

\begin{lemma}
\label{nc2lemma1ultimo}
The functor $U$ commutes with small colimits.
\begin{proof}
 Let $\{F_{i}\}_{i\in I}$ be a diagram in $\C$. Then, by definition we have a pullback diagram

\begin{equation}
\xymatrix{
U(colim_I \, F_{i})\ar[r] \ar[d]& (colim_I \, F_{i})_{L_{pe}(\mathbb{A}^1)} \coprod_{ (colim_I \, F_{i})} (colim_I \, F_{i})_{L_{pe}(\mathbb{A}^1)} \ar[d]\\
\ar[r]0&  (colim_I \, F_{i})_{L_{pe}(\mathbb{G}_m)}
}
\end{equation}

\noindent but since $(-)_{L_{pe}(\mathbb{G}_m)}$ and $(-)_{L_{pe}(\mathbb{A}^1)}$ commute with all colimits (thanks to Yoneda's lemma and the fact the evaluation map commutes with small colimits), this diagram is equivalent to

\begin{equation}
\xymatrix{
U(colim_I \, F_{i})\ar[r] \ar[d]& colim_I \, ((F_{i})_{L_{pe}(\mathbb{A}^1)} \coprod_{ F_{i}} (F_{i})_{L_{pe}(\mathbb{A}^1)}) \ar[d]\\
\ar[r]0&  colim_I \, ((F_{i})_{L_{pe}(\mathbb{G}_m)})
}
\end{equation}

\noindent and since in the stable context colimits commute with pullbacks we find a canonical equivalence

\begin{equation}
U(colim_I \, F_{i})\simeq colim_I \, U(F_{i})
\end{equation}

\end{proof}
\end{lemma}

\begin{lemma}
\label{nc2lemma2ultimo}
The two maps $U= Id\circ U\to U^2$ and $U= U\circ Id\to U^2$ induced by the natural transformation $Id\to U$, are homotopic. 
\begin{proof}
We are reduced to show that for any $F$ the natural maps $U(\alpha_F), \alpha_{U(F)}: U(F)\to U^2(F)$ are homotopic. Recall that by definition $\alpha_{U(F)}$ is determined by the universal property of pullbacks, as being the essentially unique map that makes the diagram

\begin{equation}
\label{nc2vista122}
\xymatrix{
U(F)\ar@{-->}[dr]^{\alpha_{U(F)}}\ar@/^1pc/[drr]\ar@/_1pc/[ddr]&&\\
&U^2(F)\ar[d]\ar[r]& U(F)_{L_{pe}(\mathbb{A}^1)}\coprod_{U(F)} U(F)_{L_{pe}(\mathbb{A}^1)}\ar[d]\\
&0\ar[r]&U(F)_{L_{pe}(\mathbb{G}_m)}
}
\end{equation}

\noindent commute. In this case, as $U$ commutes with colimits by the Lemma \ref{nc2lemma1ultimo} and as $Fun(\dg^{ft}, \widehat{Sp})$ is stable, $U$ also preserves pullbacks. In this case, and as we have equivalences  $U(F_{\mathbb{G}_m})\simeq U(F)_{{\mathbb{G}_m}}$ and $U(F_{L_{pe}(\mathbb{A}^1)})\simeq U(F)_{L_{pe}(\mathbb{A}^1)}$, the diagram in (\ref{nc2vista122}) is in fact equivalent to the image of the diagram

\begin{equation}
\label{nc2vista123}
\xymatrix{
F\ar@{-->}[dr]^{\alpha_{F}}\ar@/^1pc/[drr]\ar@/_1pc/[ddr]&&\\
&U(F)\ar[r]\ar[d]& F_{L_{pe}(\mathbb{A}^1)}\coprod_{F}F_{L_{pe}(\mathbb{A}^1)}\ar[d]\\
&0\ar[r]&F_{L_{pe}(\mathbb{G}_m)}
}
\end{equation}

\noindent under $U$, where the inner commutative square is a pullback so that $U(\alpha_F)$ is necessarily homotopic to $\alpha_{U(F)}$.
\end{proof}
\end{lemma}

These lemmas have the following consequences:

\begin{prop}
\label{nc2lemma3ultimo}
The natural transformation $(-)^B\circ Id\to (-)^B\circ U$ is an equivalence.
\begin{proof}
This amounts to check that for any $F$ the natural map $F^B\to U(F)^B$ is an equivalence. By construction, this is the map induced at the colimit level by the morphism of diagrams

\begin{equation}
\xymatrix{
...& ...\\
U^2(F)\ar[r]^{U^2(\alpha_{F})}\ar[u]^{\alpha_{U^2(F)}}& U(F)\ar[u]_{\alpha_{U^3(F)}}\\
U(F)\ar[r]^{U(\alpha_{F})}\ar[u]^{\alpha_{U(F)}}& U(F)\ar[u]_{\alpha_{U^2(F)}}\\
F\ar[u]^{\alpha_F}\ar[r]^{\alpha_F}& U(F)\ar[u]_{\alpha_{U(F)}}
}
\end{equation}

By iterating the Lemma \ref{nc2lemma2ultimo} we find that for any $k\geq 0$ the maps $U^k(\alpha_F)$ and $\alpha_{U^k(F)}$ are homotopic so that, by cofinality, the map $F^B\to U(F)^B$ induced between the colimit of each column is an equivalence.
\end{proof}
\end{prop}

\begin{prop}
\label{nc2lemma4ultimo}
The natural map $(-)^B\circ U\to U\circ (-)^B$ is an equivalence.
\begin{proof}
It is enough to show that for any $F$ the natural map $U(F)^B\to U(F^B)$ is an equivalence. As $F^B$ can be obtained as a colimit for the sequence (\ref{nc2chinesmacaco}) and as $U$ commutes with colimits, $U(F^B)$ is the colimit of 

\begin{equation}
\xymatrix{
U(F)\ar[rr]^{U(\alpha_F)} && U^2(F)\ar[rr]^{U(\alpha_{U(F)})} &&U^3(F)\ar[rr]^{U(\alpha_{U^2(F)})}&& ....
}
\end{equation}

\noindent and again, by using the Lemma \ref{nc2lemma2ultimo} together with cofinality, we deduce that this colimit is equivalent to $U(F)^B$.

\end{proof}
\end{prop}

We can now put these results together and show that

\begin{cor}
\label{FBbass}
For any object $F\in Fun(\dg^{ft}, \widehat{\Sp})$ the object $F^B$ satisfies the Bass-Thomason-Trobaugh exact sequences for any $n\in \mathbb{Z}$.
\begin{proof}
By combining the Propositions \ref{nc2lemma3ultimo} and \ref{nc2lemma4ultimo} we deduce that the canonical map $F^B\to U(F^B)$ is an equivalence. Therefore, we have a pullback-pushout square

\begin{equation}
\label{nc2hkilkklaka}
\xymatrix{
 (F^B)_{L_{pe}(\mathbb{A}^1)}\coprod_{F^B} (F^B)_{L_{pe}(\mathbb{A}^1)}\ar[d]\ar[r] & 0\ar[d]\\
(F^B)_{L_{pe}(\mathbb{G}_m)}\ar[r]& \Sigma F^B
}
\end{equation}

\noindent and using exactly the same arguments as in \ref{bassexactsequencefornisnevichlocal} we find that for any $T_{\X} \in \dg^{ft}$, the associated long exact sequence breaks up into short exact sequences

\begin{equation}
0\to\pi_n(F^B(L_{pe}(\mathbb{A}^1)\otimes T_{\X} )\coprod_{F^B(T_{\X} )} F^B(L_{pe}(\mathbb{A}^1)\otimes T_{\X} ))\to \pi_n(F^B(L_{pe}(\mathbb{G}_m)\otimes T_{\X} ))\to \pi_{n-1}(F^B( T_{\X} ))\to 0
\end{equation}

\noindent and again by the same arguments we are able to extract the familiar exact sequences of Bass-Thomason-Trobaugh, for all $n\in \mathbb{Z}$.
\end{proof}
\end{cor}

\begin{remark}
\label{FBretract}
As the canonical map  $F^B\to U(F^B)$ is an equivalence it follows from the Construction \ref{Bconstruction} that when we construct the diagram (\ref{nc2hkilla}) with $F^B$

\begin{equation}
\label{nc2hkilkkla}
\xymatrix{
&\Omega (F^B_{L_{pe}(\mathbb{G}_m)})\ar[r]\ar[d]&0\ar[d]\\
F^B\ar@{-->}[ur]^{\sigma_{F^B}}\ar[r]&U(F^B)\ar[r]\ar[d]& (F^B)_{L_{pe}(\mathbb{A}^1)}\coprod_{F^B} (F^B)_{L_{pe}(\mathbb{A}^1)}\ar[d]\\
&0\ar[r]&(F^B)_{L_{pe}(\mathbb{G}_m)}
}
\end{equation}

\noindent the section $\sigma_{F^B}$ makes $F^B$ a retract of $\Omega (F^B_{L_{pe}(\mathbb{G}_m)})$. In particular, by iteratively applying the construction $\Omega(-)_{_{L_{pe}(\mathbb{G}_m)}}$ we find (because $\Sigma^{\infty}_{+}\circ j_{nc}$ is monoidal) that for any $n\geq 1$, the composition

\begin{equation}
F^B\to \Omega (F^B_{L_{pe}(\mathbb{G}_m)}) \to ... \to \Omega^n (F^B_{L_{pe}(\mathbb{G}_m)^{\otimes n}})\to ... \to \Omega (F^B_{L_{pe}(\mathbb{G}_m)})\to F^B
\end{equation}

\noindent is the identity map so that, for any $n\geq 1$, $F^B$ is a rectract of $ \Omega^n (F^B_{L_{pe}(\mathbb{G}_m)^{\otimes n}})$. Equivalently, for any $n\geq 1$, the suspension $\Sigma^n F^B$ is a retract of $(F^B)_{L_{pe}(\mathbb{G}_m)^{\otimes n}}$.
\end{remark}

We will now show that the construction $(-)^B$ defines a localization:

\begin{prop}
\label{Bislocalization}
The functor $(-)^B: Fun(\dg^{ft}, \widehat{\Sp})\to Fun(\dg^{ft}, \widehat{\Sp})$  of the Construction \ref{Bconstruction} is an accessible localization functor.
\end{prop}

This result follows from the Lemmas \ref{nc2lemma1ultimo} and \ref{nc2lemma2ultimo} together with the following general result:

\begin{lemma}
\label{nc229ultimo}
Let $\C$ be a presentable $(\infty,1)$-category and let $U:\C\to \C$ be an colimit preserving endofunctor of $\C$, together with a natural transformation $f:Id_{\C}\to U$ such that the two obvious maps $U\circ Id_{\C}\to U^2$ and $Id_{\C}\circ U\to U^2$ are homotopic. Let

\begin{equation}
\label{nc229pinga}
\xymatrix{
Id \ar[dr]_{i_0}\ar[r]^f& \ar[d]_{i_1}U=Id_{\C}\circ U\ar[r]& \ar[dl]_{i_2}U^2= Id_{\C}\circ U^2\ar[r]& ....\\
& T&&
}
\end{equation}

\noindent be a colimit cone for the horizontal sequence (indexed by $\mathbb{N}$), necessarily in $Fun^L(\C,\C)$.  Then, the functor $T:\C\to \C$ provides a reflexive localization of $\C$. Moreover, since $T$ commutes with small colimits the localization is accessible.
\begin{proof}

The proof requires some preliminairs. To start with we observe that the arguments in the proof of the Propositions \ref{nc2lemma3ultimo} and  \ref{nc2lemma4ultimo} apply mutatis-mutandis to this general situation so that 
we have natural equivalences $T\circ Id_{\C}\simeq T\circ U$ and $T\circ U\simeq U\circ T$. These two facts combined force the canonical maps 

\begin{equation}
\xymatrix{
T \ar[r]& U\circ T\ar[r]& U^2\circ T\ar[r]& ...
}
\end{equation}

\noindent to be equivalences.

Let us now explain the main proof. For this purpose we will use the description of a reflexive localization functor given in \cite[5.2.7.4-(3)]{lurie-htt}. Namely, for a functor $T:\C\to \C$ from an $(\infty,1)$-category $\C$ to itself to provide a reflexive localization of $\C$  (which we recall, means that $T$ factors as $\C\to \C_0\subseteq \C$ with $\C_0$ a full subcategory of $\C$, $\C_{0}\subseteq \C$ the inclusion and $\C\to \C_0$  a left adjoint to the inclusion $\C_0\subseteq \C$) it is enough to  have $T$ equipped with a natural transformation $\alpha: Id_{\C}\to T$ such that for every object $X\in \C$, the morphisms $\alpha_{T(X)}$ and $T(\alpha_X)$ are equivalences .

In our case, we let $\alpha$ be the canonical natural transformation $i_0: Id_{\C}\to T$ appearing in the colimit cone (\ref{nc229pinga}). We show that for any $X\in \C$, the maps $(i_0)_{T(X)}$ and $T((i_0)_X)$ are equivalences. The first follows immediately from our preliminairs: since all the maps in the sequence 

\begin{equation}
T(X)\to U(T(X))\to U^2(T(X)) \to...
\end{equation}

\noindent are equivalences and $(i_0)_{T(X)}$ is by definition the first structural map in the colimit cone of this sequence, it is also an equivalence. 

Let us now discuss  $T((i_0)_X)$. By construction of the functor $T$, this is the map $colim_{n\in \mathbb{N}}\, U^n((i_0)_X):T(X)\to T(T(X))$ induced by the universal property of colimits by means of the morphism of sequences

\begin{equation}
\label{nc2281}
\xymatrix{
X\ar[r]\ar[d]^{(i_0)_X}& U(X)\ar[r]\ar@{-->}[d]^{U((i_0)_X)}&U^2(X)\ar[r]\ar@{-->}[d]^{U^2((i_0)_X)}&...\\
T(X)\ar[r]& U(T(X))\ar[r]&U^2(T(X))\ar[r]&...
}
\end{equation}

We will prove that

\begin{enumerate}[(i)]
\item For any $X$ there is a canonical homotopy between the maps $U((i_0)_X)$ and $(i_0)_{U(X)}$. By induction we get canonical homotopies between $U^n((i_0)_X)$ and $(i_0)_{U^n(X)}$;
\item For any  diagram $I\to \C$ in $\C$ (denoted as $\{X_k\}_{k\in I}$), there is a canonical homotopy between the maps $colim_{k\in I}\, ((i_0)_{X_k})$ and $(i_0)_{colim_{k\in I}\, (X_k)}$.
\end{enumerate}

\noindent so that by combining these two results we get

\begin{equation}
T((i_0)_X)\simeq colim_{n\in \mathbb{N}}\,U^n((i_0)_X)\simeq colim_{n\in \mathbb{N}}\,((i_0)_{U^n(X)})\simeq (i_0)_{colim_{n\in \mathbb{N}}\, (U^n(X))}\simeq (i_0)_{T(X)}
\end{equation}

\noindent and since we already know that $(i_0)_{T(X)}$ is an equivalence, we deduce the same for $T((i_0)_X)$.\\ 

To prove $(i)$ we observe the existence of a canonical commutative triangle 

\begin{equation}
\xymatrix{
T(U(X))\ar@{-->}[r]& U(T(X))\\
U(X) \ar[ru]_{U((i_0)_X)}\ar[u]^{(i_0)_{U(X)}}&
}
\end{equation}

\noindent  provided by the universal property of the colimit defining $T$. As explained in the preliminairs this dotted map is an equivalence so that the commutativity of this diagram holds the desired homotopy.\\

Let us now prove $(ii)$. Let  $I^{\triangleright}\to \C$ be a colimit diagram in $\C$ (which, by abusing the notation we denote as $\{X_{k}, \phi_k: X_k\to colim_{k\in I}\,X_k\}_{k\in I}$). Since $i_0$ is a natural transformation we find for any $k\in $I a commutative diagram

\begin{equation}
\label{nc2291}
\xymatrix{
X_k\ar[r]^{\phi_k}\ar[d]_{(i_0)_{(X_k)}}& colim_{k\in I}\,X_k\ar[d]^{(i_0)_{(colim_{k\in I}\,X_k)}}\\
T(X_k)\ar[r]^(.4){T(\phi_k)}& T(colim_{k \in I}\,X_k)
}
\end{equation}

\noindent and the universal property of colimits allows us to factor the lower horizontal arrows as

\begin{equation}
\xymatrix{
T(X_k)\ar[r] & colim_{k\in I}\,T(X_k)\ar@{-->}[r]^{\theta}& T(colim_{k\in I}\,X_k)
}
\end{equation}

\noindent where the dotted map $\theta$ is essentially unique. More importantly, since by construction $T$ commutes with colimits, $\theta$ is an equivalence. 

At the same time, the map $colim_{k\in I}\,(i_0)_{(X_k)}$ is by definition the essentially unique map $colim_{k\in I}\,X_k\to colim_{k\in I}\,T(X_k)$, induced by the universal property  colimits, that makes the diagrams

\begin{equation}
\xymatrix{
X_k\ar[r]^{\phi_k}\ar[d]_{(i_0)_{(X_k)}}& colim_{k\in I}\,X_k\ar@{-->}[d]^{colim_{k\in I}\,(i_0)_{(X_k)}}\\
T(X_k)\ar[r]&colim_{k\in I}T(X_k)
}
\end{equation}

\noindent commute. Finally, since $\theta$ is an equivalence, the commutativity of (\ref{nc2291}) implies the commutativity of 

\begin{equation}
\xymatrix{
X_k\ar[r]^{\phi_k}\ar[d]_{(i_0)_{(X_k)}}& colim_{k\in I}\,X_k\ar[d]^{\theta^{-1}\circ (i_0)_{(colim_{k\in I}\,X_k)}}\\
T(X_k)\ar[r]&colim_{k\in I}\,T(X_k)
}
\end{equation}

\noindent so that by the uniqueness property that defines $colim_{k\in I}(i_0){(X_k)}$, the diagram

\begin{equation}
\xymatrix{
colim_{k\in I}\,(X_k)\ar[d]_{colim_{k\in I}\,(i_0)_{(X_k)}} \ar[dr]^{(i_0)_{(colim_{k\in I}\,X_k)}}&\\
colim_{k\in I}\,T(X_k)\ar[r]^{\theta}&T(colim_{k\in I}\,X_k)
}
\end{equation}

\noindent must commute. This provides the desired homotopy and concludes $(ii)$ and the proof.

\end{proof}
\end{lemma}

\begin{remark}
\label{nc2nisnevichlocalimpliesBlocal}
It follows from the the Proposition \ref{Bislocalization} and from the Construction \ref{Bconstruction} that an object $F\in Fun(\dg^{ft},\widehat{\Sp})$ is local with respect to the localization $(-)^B$ if and only if the diagram

\begin{equation}
\label{nc2231}
\xymatrix{
F\ar[r]\ar[d]& F_{L_{pe}(\mathbb{A}^1)}\coprod_{F} F_{L_{pe}(\mathbb{A}^1)}\ar[d]\\
0\ar[r]&F_{L_{pe}(\mathbb{G}_m)}
}
\end{equation}

\noindent is a pullback-pushout square. In particular, the discussion in \ref{bassexactsequencefornisnevichlocal} implies that any Nisnevich local object $F$ is $(-)^B$-local.
\end{remark}

We now come to a series of technical steps in order prove each of the items in \ref{nc2proposicao1}. First thing, we give a precise sense to what it means for a functor $F$ with connective values to satisfy  all the  Bass exact sequences for $n\geq 1$.

\begin{defn}
\label{nc2232}
Let $F\in Fun(\dg^{ft}, \widehat{\Sp}_{\geq 0})$ and consider its associated diagram (\ref{nc2hkilla}) (constructed in $Fun(\dg^{ft}, \widehat{\Sp})$, where we identify $F$ with its inclusion). We say that $F$ satisfies all Bass exact sequences for $n\geq 1$ if the canonical induced map of connective functors $F\to \tau_{\geq 0}U(F)$ is an equivalence, or, in other words, since $\tau_{\geq 0}$ commutes with limits and because of the definition of $U(F)$, if the diagram (\ref{nc2231}) is a pullback in $Fun(\dg^{ft}, \widehat{\Sp}_{\geq 0})$.
\end{defn}

\begin{remark}
Let $F\in Fun(\dg^{ft}, \widehat{\Sp}_{\geq 0})$ and consider the pullback-pushout diagram in $Fun(\dg^{ft}, \widehat{\Sp})$

\begin{equation}
\xymatrix{
 \ar[d]\Omega(F_{L_{pe}(\mathbb{A}^1)}\coprod_{F} F_{L_{pe}(\mathbb{A}^1)})\ar[r]& \Omega(F_{L_{pe}(\mathbb{G}_m)})\ar[d]\\
0\ar[r]&U(F)
}
\end{equation}

Since, $\tau_{\geq 0}$ preserves pullbacks, we obtain a pullback diagram in $Fun(\dg^{ft}, \widehat{\Sp}_{\geq 0})$

\begin{equation}
\label{nc2251}
\xymatrix{
 \ar[d]\tau_{\geq 0}\Omega(F_{L_{pe}(\mathbb{A}^1)}\coprod_{F} F_{L_{pe}(\mathbb{A}^1)})\ar[r]& \tau_{\geq 0}\Omega(F_{L_{pe}(\mathbb{G}_m)})\ar[d]\\
0\ar[r]&\tau_{\geq 0}U(F)
}
\end{equation}

If $F$ satisfies the condition in the previous definition, then the zero truncation of the composition

\begin{equation}
 \xymatrix{F\ar@{-->}[r]^(.3){\sigma_F}&\Omega(F_{L_{pe}(\mathbb{G}_m)})\ar[r]& U(F)}
\end{equation}

\noindent makes $F$ a retract of $\tau_{\geq 0}\Omega(F_{L_{pe}(\mathbb{G}_m)})$. With this, and as before, once evaluated at $T_{\X}\in \dg^{ft}$, the long exact sequence associated to the pullback (\ref{nc2251}) splits up into short exact sequences

\begin{equation}
0\to\pi_n(F(L_{pe}(\mathbb{A}^1)\otimes T_{\X} )\coprod_{F(T_{\X} )} F(L_{pe}(\mathbb{A}^1)\otimes T_{\X} ))\to \pi_n(F(L_{pe}(\mathbb{G}_m)\otimes T_{\X} ))\to \pi_{n-1}(F( T_{\X} ))\to 0
\end{equation}

\noindent $\forall n\geq 1$, and again by the same arguments, we can extract the exact sequences of Bass-Thomason.
\end{remark}

\begin{lemma}
\label{nc2lemakiki1}
If $F\in Fun(\dg^{ft}, \widehat{\Sp})$ has connective values and satisties all the Bass exact sequences for $n\geq 1$ (in the sense of the Definition \ref{nc2232}),  then the canonical map $F\simeq \tau_{\geq 0} F\to \tau_{\geq 0}F^B$ is an equivalence. 
\begin{proof}

Assuming that $F$ satisfies the condition in the  Definition \ref{nc2232}, meaning the canonical map $F\to \tau_{\geq 0}U(F)$ is an equivalence, we will show that for any $k\geq 2$, the canonical map  $F\to \tau_{\geq 0}U^k(F)$ is an equivalence. Once we have this, the conclusion of the lemma will follow from the fact $\tau_{\geq 0}$ commutes with filtered colimits (because the $t$-structure in $\widehat{\Sp}$ is determined by the stable homotopy groups and these commute with filtered colimits), so that

\begin{equation}
\tau_{\geq 0}(F^B)\simeq \tau_{\geq 0}(colim_{i\in \mathbb{N}}U^i(F))\simeq colim_{i\in \mathbb{N}}\tau_{\geq 0}(U^i(F))\simeq colim_{i\in \mathbb{N}}F\simeq F
\end{equation}

So, let us prove the assertion for $k=2$. By definition, we have a pullback-pushout square in $Fun(\dg^{ft}, \widehat{\Sp})$

\begin{equation}
\xymatrix{
U^2(F)\ar[d]\ar[r]& U(F)_{L_{pe}(\mathbb{A}^1)}\coprod_{U(F)} U(F)_{L_{pe}(\mathbb{A}^1)}\ar[d]\\
0\ar[r]&U(F)_{L_{pe}(\mathbb{G}_m)}
}
\end{equation}

\noindent and since $\tau_{\geq 0}$ preserves pullbacks, we find

\begin{equation}
\tau_{\geq 0}U^2(F)\simeq \tau_{\geq 0}( U(F)_{L_{pe}(\mathbb{A}^1)}\coprod_{U(F)} U(F)_{L_{pe}(\mathbb{A}^1)})\times_{\tau_{\geq 0}(U(F)_{L_{pe}(\mathbb{G}_m)})}0
\end{equation}

We observe that

\begin{enumerate}[(i)]
\item $\tau_{\geq 0}(U(F)_{L_{pe}(\mathbb{G}_m)})\simeq F_{L_{pe}(\mathbb{G}_m)}$.\\

\item  $\tau_{\geq 0}( U(F)_{L_{pe}(\mathbb{A}^1)}\coprod_{U(F)} U(F)_{L_{pe}(\mathbb{A}^1)})\simeq F_{L_{pe}(\mathbb{A}^1)}\coprod_{F} F_{L_{pe}(\mathbb{A}^1)}$

\end{enumerate}

To deduce the first equivalence, we use the equivalence $\tau_{\geq 0}U(F)\simeq F$ together with the fact that $(-)_{L_{pe}(\mathbb{G}_m)}$ commutes with $\tau_{\geq 0}$. The second equivalence requires a more sophisticated discussion. Recall from the section \ref{bassexactsequencefornisnevichlocal} that for any $G\in Fun(\dg^{ft}, \widehat{\Sp})$ we are able to construct a pushout square in $Fun(\dg^{ft}, \widehat{\Sp})$

\begin{equation}
\label{nc225h1}
\xymatrix{
\ar[d] G\ar[r]^{i^G_1} &G\oplus G\ar[r]^{} & G_{L_{pe}(\mathbb{P}^1)}\ar[r] & G_{L_{pe}(\mathbb{A}^1)}\oplus G_{L_{pe}(\mathbb{A}^1)} \ar[d]\\
0 \ar[rrr]&&&G_{L_{pe}(\mathbb{A}^1)}\coprod_{G} G_{L_{pe}(\mathbb{A}^1)}
}
\end{equation}

\noindent such that the top horizontal composition admits a left inverse.  Applying this construction to $G=F$ and to $G=U(F)$, we construct a map between the associated pullback-pushout squares

\begin{equation}
\label{nc225h2}
\xymatrix{
&U(F)\ar[rr]\ar[dd]&&U(F)_{L_{pe}(\mathbb{A}^1)}\oplus U(F)_{L_{pe}(\mathbb{A}^1)} \ar[dd]\\
F\ar[ru]\ar[rr]\ar[dd]&&\ar[ru]\ar[dd]F_{L_{pe}(\mathbb{A}^1)}\oplus F_{L_{pe}(\mathbb{A}^1)} &\\
&0\ar[rr]&&U(F)_{L_{pe}(\mathbb{A}^1)}\coprod_{U(F)} U(F)_{L_{pe}(\mathbb{A}^1)}\\
0\ar[ru]\ar[rr]&&\ar[ru]F_{L_{pe}(\mathbb{A}^1)}\coprod_{F} F_{L_{pe}(\mathbb{A}^1)}&\\
}
\end{equation}

\noindent (obtained using the natural transformation $\underline{Hom}(-,F)\to \underline{Hom}(-,U(F))$ induced by canonical morphism $F\to U(F)$).

 Both the front and back faces are pullback-pushouts and both the top horizontal maps admite left-inverses. 

Finally, since $\tau_{\geq 0}U(F)\simeq F$ and because the top horizontal maps  admit left-inverses, the long exact sequences associated to each square breaks up into short exact sequences, and for each $n\geq 0$ and each $T_{\X}\in \dg^{ft}$ we find natural maps of short exact sequences

\begin{equation}
\label{nc225h13}
\xymatrix{
  \pi_n(U(F)(T_{\X} ))\ar[r]& \pi_n((U(F)_{L_{pe}(\mathbb{A}^1)}\oplus U(F)_{L_{pe}(\mathbb{A}^1)})(T_{\X}))\ar[r]&\pi_n((U(F)_{L_{pe}(\mathbb{A}^1)}\coprod_{U(F)} U(F)_{L_{pe}(\mathbb{A}^1)})(T_{\X}))\\
  \pi_n(F(T_{\X} ))\ar[r]\ar[u]^{\sim}&\pi_n((F_{L_{pe}(\mathbb{A}^1)}\oplus F_{L_{pe}(\mathbb{A}^1)} )(T_{\X}))\ar[r]\ar[u]^{\sim}& \pi_n((F_{L_{pe}(\mathbb{A}^1)}\coprod_{F} F_{L_{pe}(\mathbb{A}^1)})(T_{\X}))\ar[u]
}
\end{equation}

\noindent implying the equivalence in $(ii)$.\\

Finally, we deal with the case $k>2$. Applying the same strategy for $G=F$ and $G=U^k(F)$, we consider the analogue of the diagram (\ref{nc225h2}) induced by the canonical morphism $F\to U^k(F)$. By induction, we deduce that $\tau_{\geq 0}U^{k+1}(F)\simeq F$. This concludes the proof.

\end{proof}
\end{lemma}

\begin{prop}
\label{nc2lemakiki2}
Let $F$ be a connectively-Nisnevich local object. Then, it satisfies the Projective Bundle Theorem and all the Bass exact sequences for $n\geq 1$. In particular,  by the Lemma \ref{nc2lemakiki1} we have $F\simeq \tau_{\geq 0} F\simeq \tau_{\geq 0}F^B$. 
\begin{proof}

To start with, we prove that if is $F$ connectively-Nisnevich local then it satisfies the Projective bundle theorem. Indeed, we can use the arguments used in \ref{bassexactsequencefornisnevichlocal} together with the definition of being connectively-Nisnevich local to construct a pullback diagram in $Fun(\dg^{ft}, \widehat{\Sp}_{\geq 0})$

\begin{equation}
\xymatrix{
F\simeq F_{L_{pe}(k)}\ar[r]\ar[d]& F_{L_{pe}(\mathbb{P}^1)}\ar[d]\\
0\ar[r] & F\simeq F_{L_{pe}(k)}
}
\end{equation}

\noindent with splittings, which, as explained in the Remark \ref{nc2remarksplit}, provide a canonical equivalence $F_{L_{pe}(\mathbb{P}^1)}\simeq F\oplus F$ in $Fun(\dg^{ft}, \widehat{\Sp})$.
Secondly, and again by the definition of connectively-Nisnevich local, we can easily deduce that the canonical diagram

\begin{equation}
\label{nc223pika}
\xymatrix{
F_{L_{pe}(\mathbb{P}^1)}\ar[d]\ar[r] &  F_{L_{pe}(\mathbb{A}^1)}\oplus F_{L_{pe}(\mathbb{A}^1)} \ar[d]\\
0\ar[r]&F_{L_{pe}(\mathbb{G}_m)}
}
\end{equation}

\noindent associated to the covering of $\mathbb{P}^1$ by two affine lines (\ref{nc2classicalcoveringP1}) is a pullback in $Fun(\dg^{ft}, \widehat{\Sp}_{\geq 0})$.

With these two ingredients we prove that if is $F$ connectively-Nisnevich local then it satisfies all the Bass exact sequences for $n\geq 1$ in the sense of the Definition \ref{nc2232}, namely, we show that the canonical map $F\simeq \tau_{\geq 0}F \to \tau_{\geq 0}U(F)$ is an equivalence, or, in other words, that the diagram (\ref{nc2231}) is a pullback within connective functors.

Consider the pushout squares in $Fun(\dg^{ft}, \widehat{\Sp})$ described in (\ref{nc2h1}). More precisely, since $F$ satisfies the Projective bundle theorem, we are interested in the pullback-pushout square

\begin{equation}
\label{nc2233h1}
\xymatrix{
\ar[d]F\oplus F \simeq  F_{L_{pe}(\mathbb{P}^1)}\ar[d]\ar[r] & F_{L_{pe}(\mathbb{A}^1)}\oplus F_{L_{pe}(\mathbb{A}^1)} \ar[d]\\
 \ar[r] F\simeq F\coprod_{F\oplus F}F_{L_{pe}(\mathbb{P}^1)}\ar[r]&F_{L_{pe}(\mathbb{A}^1)}\coprod_{F} F_{L_{pe}(\mathbb{A}^1)}
}
\end{equation}

\noindent which, in particular,  is a pullback square in $Fun(\dg^{ft}, \widehat{\Sp}_{\geq 0})$ once truncated at level zero. Combining with the pullback square (\ref{nc223pika}) we find a series of pullback squares in $Fun(\dg^{ft}, \widehat{\Sp}_{\geq 0})$.

\begin{equation}
\label{nc2231823}
\xymatrix{
\Omega F\ar[d]\ar[r]&\ar[d]0&\\
                      \Omega(F_{L_{pe}(\mathbb{A}^1)}\coprod_{F} F_{L_{pe}(\mathbb{A}^1)})\ar[dd]\ar[r]&\ar[d]\Omega(F_{L_{pe}(\mathbb{G}_m)})  \ar[r]&0\ar[d]\\
                      &   F_{L_{pe}(\mathbb{P}^1)}\ar[d]\ar[r]&  F_{L_{pe}(\mathbb{A}^1)}\oplus F_{L_{pe}(\mathbb{A}^1)} \ar[d]\\
                 0\ar[r]     & F\ar[r] &F_{L_{pe}(\mathbb{A}^1)}\coprod_{F} F_{L_{pe}(\mathbb{A}^1)}
}
\end{equation}

Now comes the important ingredient: since the diagram (\ref{nc223pika}) is a pullback, we can still deduce (as before) the existence of a canonical map $\sigma_{F}$
such that the composition

\begin{equation}
\xymatrix{
F\ar@{-->}[r]^{\sigma_F}& \Omega F_{L_{pe}(\mathbb{G}_m)}\ar[r]& F_{L_{pe}(\mathbb{P}^1)}\ar[r]& F
}
\end{equation}

\noindent is the identity. We now explain how the existence of this section allows us to prove that the diagram (\ref{nc2231}) is a pullback. More precisely, by using $\sigma_F$ at each copy of $F$ in  (\ref{nc2231}) and applying the construction $\Omega(-)_{L_{pe}(\mathbb{G}_m)}$ we find the square (\ref{nc2231}) as a retract of the square

\begin{equation}
\label{nc2231jdj}
\xymatrix{
\Omega F_{L_{pe}(\mathbb{G}_m)}\ar[r]\ar[d]& (\Omega F_{L_{pe}(\mathbb{G}_m)})_{L_{pe}(\mathbb{A}^1)}\coprod_{\Omega(F)_{L_{pe}(\mathbb{G}_m)}} (\Omega F_{L_{pe}(\mathbb{G}_m)})_{L_{pe}(\mathbb{A}^1)}\ar[d]\\
0\ar[r]&\Omega (F_{L_{pe}(\mathbb{G}_m)})_{L_{pe}(\mathbb{G}_m)}
}
\end{equation}

\noindent but since both $\Omega$ and  $\underline{Hom}(\Sigma^{\infty}_{+}\circ j_{nc}(L_{pe}(\mathbb{G}_m)),-)$ commute with colimits, we can easily indentify this last square with the image of the top left pullback square in (\ref{nc2231823}) under $\underline{Hom}(\Sigma^{\infty}_{+}\circ j_{nc}(L_{pe}(\mathbb{G}_m)),-)$ and conclude that this is also a pullback square. We conclude the proof using the fact that the rectract of a pullback square is a pullback.

\end{proof}
\end{prop}

We now address the second item of the Proposition \ref{nc2proposicao1}, namely, 

\begin{prop}
\label{nc2proporkikiki}
Let $F\in Fun_{Nis_{\geq 0}}(\dg^{ft}, \widehat{\Sp}_{\geq 0})$. Then, the object $(i\circ \alpha(F))^{B}$ is Nisnevich local.
\end{prop}

The proof of this proposition is based on a very helpful criterium to decide if a given $F$ is Nisnevich local by studying its truncations $\tau_{\geq 0}\Sigma^n F$, namely:

\begin{lemma}
\label{nc2criteriumnisnevich}
Let $F$ be any object in $Fun(\dg^{ft}, \widehat{\Sp})$. Then, if for any $n\geq 0$ the truncations $\tau_{\geq 0}\Sigma^n F$ are connectively-Nisnevich local, the object $F$ itself is Nisnevich local.
\end{lemma}

This lemma follows from a somewhat more general situation, which we isolate in the following remark:

\begin{remark}
\label{nc2pullbacksandtruncatedpullbacks}
Let $\C$ be a stable $(\infty,1)$-category with a right-complete $t$-structure $(\C_{\geq 0}, \C_{\leq 0})$ and let  $\tau_{\geq n}$ and $\tau_{\leq n}$ denote the associated truncation functors (see \cite[Section 1.2.1]{lurie-ha} for the complete details or \refnci{Tstructures} for a fast review of the subject ).  We observe that a commutative square 

\begin{equation}
\label{nc2pikachucomsida}
\xymatrix{
A\ar[d]\ar[r]& B\ar[d]\\
C\ar[r]& D
}
\end{equation}

\noindent in $\C$ is a pullback (therefore pushout) if and only if for any $n\leq 0$ the truncated squares

\begin{equation}
\xymatrix{
\tau_{\geq n}A\ar[d]\ar[r]& \tau_{\geq n}B\ar[d]\\
\tau_{\geq n} C\ar[r]& \tau_{\geq n}D
}
\end{equation}

\noindent are pullbacks in $\C_{\geq n}$. Indeed, if we let $H$ denote the pullback of the square in $\C$, we want to show that the canonical map $A\to H$ in $\C$  induced by the universal property of the pullback, is an equivalence. But, since the truncation functors $\tau_{\geq n}$ are right adjoints to the inclusions $\C_{\geq n}\subseteq \C$,  $\tau_{\geq n}H$ is a pullback for the square in $\C_{\geq n}$ and therefore the induced maps $\tau_{\geq n}A\to \tau_{\geq n} H$ are equivalences for all $n\leq 0$. To conclude, we are reduce to show that if a map $f:X\to Y$ in $\C$ induces equivalences $\tau_{\geq n}X\simeq \tau_{\geq n}Y$ for all $n\leq 0$ then the map $f$ itself is an equivalence. To see this, and because $\C$ is stable if suffices to check that the fiber $fib(f)$ is equivalent to zero. This fiber fits in pullback-pushout square

\begin{equation}
\xymatrix{
fib(f)\ar[r]\ar[d]& \ar[d] X\\
0\ar[r]& Y
}
\end{equation}

\noindent and since $\tau_{\geq n}$ commutes with pullbacks and the maps $\tau_{\geq n}X\to \tau_{\geq n}Y$ are equivalences, we find that for any $n\leq 0$ we have $\tau_{\geq n}fib(f)\simeq 0$. Finally, we use the canonical pullback-pushout squares in $\C$

\begin{equation}
\xymatrix{
\tau_{\geq n}fib(f)\ar[r]\ar[d]& fib(f)\ar[d]\\
0\ar[r]& \tau_{\leq n-1} fib(f)
}
\end{equation}

\noindent to deduce that for all $n\leq -1$ the map $fib(f)\to \tau_{\leq n} fib(f)$ is an equivalence. In particular $fib(f)$ belongs to the intersection $\cap_n \C_{\leq n}$ so that, since the $t$-structure is assumed to be right-complete,  we have $fib(f)\simeq 0$.

In particular, since the truncations $\tau_{\geq n}$ can be obtained as the compositions $\Omega^n\circ \tau_{\geq 0}\circ \Sigma^{n}$ and since $\Omega$ commutes with limits, the previous discussion implies that for the square (\ref{nc2pikachucomsida}) to be a pullback in $\C$ it suffices to have for each $n\geq 0$, the induced square

\begin{equation}
\xymatrix{
\tau_{\geq 0}\Sigma^n A\ar[d]\ar[r]& \tau_{\geq 0}\Sigma^n B\ar[d]\\
\tau_{\geq 0}\Sigma^n C\ar[r]& \tau_{\geq 0}\Sigma^n D
}
\end{equation}

\noindent a pullback in $\C_{\geq 0}$.
\end{remark}

\textit{Proof of the Lemma \ref{nc2criteriumnisnevich}:}
Just apply the Remark \ref{nc2pullbacksandtruncatedpullbacks} to the commutative squares of spectra

\begin{equation}
\xymatrix{
F(T_{\X})\ar[r] \ar[d] & F(T_{\UU})\ar[d]\\
F(T_{\V})\ar[r] & F(T_{\W})
}
\end{equation}

\noindent induced by the Nisnevich squares of noncommutative spaces. The discussion therein works because the $t$-structure in $\widehat{\Sp}$ is known to be right-complete (see \cite[1.4.3.6]{lurie-ha}).

\hfill $\qed$

\textit{Proof of the Proposition \ref{nc2proporkikiki}:}
As explained in the Remark \ref{FBretract}, for any $n\geq 1$, the suspension $\Sigma^n F^B$ is a retract of $(F^B)_{L_{pe}(\mathbb{G}_m)^{\otimes n}}$. In particular, $\tau_{\geq 0} \Sigma^n F^B$ is a retract of $\tau_{\geq 0}((F^B)_{L_{pe}(\mathbb{G}_m)^{\otimes n}})$ which is a mere notation for  $\tau_{\geq 0}\underline{Hom}(\Sigma^{\infty}_{+}\circ j_{nc}(L_{pe}(\mathbb{G}_m)^{\otimes n}), F^B)$  so that

\begin{eqnarray}
\tau_{\geq 0}((F^B)_{L_{pe}(\mathbb{G}_m)^{\otimes n}})\simeq \underline{Hom}(\Sigma^{\infty}_{+}\circ j_{nc}(L_{pe}(\mathbb{G}_m)^{\otimes n}), \tau_{\geq 0}F^B)\simeq \underline{Hom}(\Sigma^{\infty}_{+}\circ j_{nc}(L_{pe}(\mathbb{G}_m)^{\otimes n}), F)
\end{eqnarray}

\noindent where the first equivalence follows because  the $t$-structure in $Fun(\dg^{ft}, \widehat{\Sp})$ is determined objectwise by the $t$-structure in $\Sp$ and the second  follows from the Proposition \ref{nc2lemakiki2}. In particular, since $F$ is connectively-Nisnevich local, $\underline{Hom}(\Sigma^{\infty}_{+}\circ j_{nc}(L_{pe}(\mathbb{G}_m)^{\otimes n}), F)$ is also connectively-Nisnevich local so that $\tau_{\geq 0} \Sigma^n F^B$ is the retract of a connectively-Nisnevich local and therefore, it is itself local \footnote{In general, the retract of a local object in a reflexive localization is local. This is, ultimately, because the retract of an equivalence is an equivalence.}. We conclude using the Lemma \ref{nc2criteriumnisnevich}, observing that for $n=0$ the condition follows by the hypothesis that $F$ is connectively-Nisnevich local.

\hfill $\qed$

Finally, 

\begin{cor}
\label{nc2corolariokikiki}
Let $F$ be any object in $Fun(\dg^{ft}, \widehat{\Sp})$. Then, there is a canonical equivalence 

\begin{equation}
(i\circ \alpha (F))^{B}\simeq l_{Nis}^{nc}((i\circ \alpha (F))
\end{equation}
\begin{proof}
This follows from the Proposition \ref{Bislocalization}, the Remark \ref{nc2nisnevichlocalimpliesBlocal} and the Proposition \ref{nc2proporkikiki}, using the universal properties of the two localizations.
\end{proof}
\end{cor}

\textit{Proof of the Proposition \ref{nc2proposicao1}:}
The three items correspond, respectively to the Propositions \ref{nc2lemakiki2}, \ref{nc2proporkikiki} and to the Corollary \ref{nc2corolariokikiki}.
The conclusion now follows from the universal property of the two localizations.

\hfill $\qed$

\subsection{Proof of the Theorem \ref{12}: Comparing the commutative and the noncommutative $\mathbb{A}^1$-localizations}
\label{comparisonA1}

We start by asking the reader to bring back to his attention the diagrams (\ref{diagramaleft}) and (\ref{diagramaright}) and to recall that after the Theorem \ref{teorema1}, together with Yoneda's lemma, $\M_2(l_{Nis}^{nc}(K^c))$  is the Bass-Thomason-Trobaugh $K$-theory of schemes. Recall also that, by definition\footnote{Either we take it as a definition or as a consequence of the explicit formula given in this section.}, Weibel's homotopy invariant $K$-theory of \cite{weibel-homotopyinvariantktheory} is the "commutative" localization $l_{\mathbb{A}^1}(\M_2(l_{Nis}^{nc}(K^c))$. With these ingredientes the conclusion of \ref{12} will follow if we prove that the commutative and noncommutative versions of the $\mathbb{A}^1$-localizations make the diagram

\begin{equation}
\label{nc2jerupiga}
\xymatrix{
\ar[d]^{l_{\mathbb{A}^1}} Fun_{Nis}(\aff^{op}, \widehat{\Sp})& \ar[d]^{l_{\mathbb{A}^1}^{nc}} \ar[l]_{\M_2} Fun_{Nis}(\dg^{ft}, \widehat{\Sp})\\
 Fun_{Nis, \mathbb{A}^1}(\aff^{op}, \widehat{\Sp})& \ar[l]_{\M_3} Fun_{Nis, L_{pe}(\mathbb{A}^1)}(\dg^{ft}, \widehat{\Sp})
}
\end{equation}

\noindent commute. In fact, we will be able to prove something slightly more general. We begin by recalling a well-known explicit formula for the $\mathbb{A}^1$-localization of presheaves of spectra. Let $\Delta_{\mathbb{A}^1}$ be the cosimplicial affine scheme given by

\begin{equation}
\Delta^{n}_{\mathbb{A}^1}:= Spec(k[t_0,...,t_n]/(t_0+... + t_n -1))
\end{equation}

Notice that at each level we have isomorphisms $\Delta^n_{\mathbb{A}^1}\simeq (\mathbb{A}^1_k)^n$. After \cite{cisinski-descentpar}, the endofunctor of $\C=Fun(\aff^{op}, \widehat{\Sp})$ defined by the formula

\begin{equation}
F\mapsto colim_{n\in \Delta^{op}}\underline{Hom}(\Delta^{n}_{\mathbb{A}^1},F)
\end{equation}

\noindent with $\underline{Hom}$ the internal-hom for presheaves of spectra, is an explicit model for the $\mathbb{A}^1$-localization in the commutative world. To see that this indeed gives something $\mathbb{A}^1$-local we use the $\mathbb{A}^1$-homotopy $m:\mathbb{A}^1\times \mathbb{A}^1\to \mathbb{A}^1$ between the identity of $\mathbb{A}^1$ and the constant map at zero. The map $m$ is given by the usual multiplication. 
It follows from this explicit description that the $\mathbb{A}^1$-localization preserves Nisnevich local objects (this is because in a stable context, sifted colimits commute with pullbacks and the Nisnevich local condition is determined by certain squares being pullbacks). \\

The important point now is that this mechanism applies mutadis-mutandis in the noncommutative world. Indeed, by taking the composition

\begin{equation}
\xymatrix{\Delta^{nc}_{\mathbb{A}^1}:\Delta \ar[r]^{\Delta_{\mathbb{A}^1}}& \aff\ar[r]^{L_{pe}}& \nck}
\end{equation}

\noindent we obtain a cosimplicial noncommutative space and as $L_{pe}$ is monoidal we get $(\Delta^{nc, n}_{\mathbb{A}^1})\simeq L_{pe}(\mathbb{A}^1)^{\otimes_n}$. Moreover, we can use exactly the same arguments to prove that the endofunctor of $\C=Fun(\nck^{op}, \widehat{\Sp})$ defined by the formula

\begin{equation}
F\mapsto  colim_{n\in \Delta^{op}}\underline{Hom}(\Delta^{nc, n}_{\mathbb{A}^1},F)
\end{equation}

\noindent is an explicit model for the noncommutative $\mathbb{A}^1$-localization functor on spectral presheaves and also by the same arguments, we conclude that Nisnevich local objects are preserved under this localization.\\

With this we can now reduce the proof that the diagram \ref{nc2jerupiga} commutes to the proof that the following diagram commutes

\begin{equation}
\xymatrix{
\ar[d]^{l_{\mathbb{A}^1}} Fun(\aff^{op}, \widehat{\Sp})& \ar[d]^{l_{\mathbb{A}^1}^{nc}} \ar[l]_{\M_1} Fun(\dg^{ft}, \widehat{\Sp})\\
 Fun_{\mathbb{A}^1}(\aff^{op}, \widehat{\Sp})& \ar[l]_{\M'} Fun_{ L_{pe}(\mathbb{A}^1)}(\dg^{ft}, \widehat{\Sp})
}
\end{equation}

\noindent where the lower part corresponds to the reflexive $\mathbb{A}^1$-localizations and $\M'$ is the right adjoint of this context obtained by the same formal arguments as $\M_2$ and $\M_3$. The commutativity of this diagram is measured by the existence of a canonical natural transformation of functors $l_{\mathbb{A}^1}\circ \M_1 \to \M'\circ l_{\mathbb{A}^1}^{nc}$ induced by the fact that $\M'$ sends $L_{pe}(\mathbb{A}^1)$-local objects to $\mathbb{A}^1$-local objects, together with the universal property of $l_{\mathbb{A}^1}$.  The diagram commutes if and only if this natural transformation is an equivalence of functors. In particular, since the diagram of right adjoints commutes

\begin{equation}
\label{Mwithinclusions}
\xymatrix{
 Fun(\aff^{op}, \widehat{\Sp})& \ar[l]_{\M_1} Fun(\dg^{ft}, \widehat{\Sp})\\
 Fun_{\mathbb{A}^1}(\aff^{op}, \widehat{\Sp})\ar@{^{(}->}[u]^{\alpha}& \ar@{^{(}->}[u]^{\beta} \ar[l]_{\M'} Fun_{ L_{pe}(\mathbb{A}^1)}(\dg^{ft}, \widehat{\Sp})
}
\end{equation}

\noindent and the vertical maps are fully-faithful, it will be enough to show that the induced natural transformation $\alpha\circ l_{\mathbb{A}^1}\circ \M_1 \to \alpha\circ\M'\circ l_{\mathbb{A}^1}^{nc}$ is an equivalence. But now, using our explicit descriptions for the $\mathbb{A}^1$-localization functors we know that for each  $F\in Fun(\dg^{ft}, \widehat{Sp})$ we have

\begin{eqnarray}
\alpha\circ l_{\mathbb{A}^1}(\M_1(F))\simeq colim_{n\in \Delta^{op}}\, \underline{Hom}( \Sigma^{\infty}_{+}\circ j(\mathbb{A}^1)^{\otimes^n}, \M_1(F))\simeq\\ 
\simeq colim_{n\in  \Delta^{op}}\, \M_1(\underline{Hom}( \Sigma^{\infty}_{+}\circ j_{nc}(L_{pe}(\mathbb{A}^1))^{\otimes^n}, F)\simeq \\
\simeq \M_1( colim_{n\in  \Delta^{op}} \,\underline{Hom}( \Sigma^{\infty}_{+}\circ j_{nc}(L_{pe}(\mathbb{A}^1))^{\otimes^n}, F)\simeq  \\
\M_1(\beta\circ l_{\mathbb{A}^1}^{nc}(F))\simeq \alpha\circ \M'\circ l_{\mathbb{A}^1}^{nc}(F)
\end{eqnarray}

\noindent where the first and penultimate equivalences follow from the explicit formulas for the $\mathbb{A}^1$-localizations, the middle equivalences follow, respectively, from the Remarks \ref{Mcompatibleinternalhom} and \ref{Mpreservescolimits} and the last equivalence follows from the commutativity of the diagram (\ref{Mwithinclusions}).\\

In particular, when applied to $F=l_{Nis}(K^c)$ we conclude the proof of the Theorem \ref{12}.

\subsection{Proof of the Theorem \ref{teorema2}: The $\mathbb{A}^1$-localization of non-connective $K$-theory is a unit for the monoidal structure in $\stnck$}
\label{proofteorema2}

We start by gathering some necessary preliminary remarks. To start with, and as explained in the Remark \refnci{usingpresheavesofspectra2} we have two different equivalent ways to construct $\stnck$: one by using presheaves of spaces, forcing Nisnevich descent and $\mathbb{A}^1$-invariance and a second one by using presheaves of spectra and forcing again the Nisnevich and $\mathbb{A}^1$-localizations. These two approaches are related by means of a commutative diagram of monoidal functors

\begin{equation}
\xymatrix{
&\nck\ar[dr]\ar[dl]&\\
Fun(\dg^{ft},\widehat{\Spaces})\ar[rr]^{\Sigma_{+}^{\infty}}\ar[d]^{l^{nc}_{0,Nis}}&&Fun(\dg^{ft},\widehat{\Sp})\ar[d]^{l^{nc}_{Nis}}\\
Fun_{Nis}(\dg^{ft},\widehat{\Spaces})\ar[rr]^{\Sigma_{+, Nis}^{\infty}}\ar[d]^{l^{nc}_{0,\mathbb{A}^1}}&&Fun_{Nis}(\dg^{ft},\widehat{\Sp})\ar[d]^{l^{nc}_{\mathbb{A}^1}}\\
\stnck:=Fun_{Nis,  L_{pe}(\mathbb{A}^1)}(\dg^{ft},\widehat{\Spaces})\ar[rr]_{\sim}^{\Sigma_{+, Nis, \mathbb{A}^1}^{\infty}}&&Fun_{Nis,  L_{pe}(\mathbb{A}^1)}(\dg^{ft},\widehat{\Sp})
}
\end{equation}

\noindent induced by the universal properties envolved and the last induced $\Sigma_{+, Nis, \mathbb{A}^1}^{\infty}$ is an equivalence because of the results in the Proposition \refnci{alreadystable}. To be completely precise we have to check that the class of maps with respect to which we localize the theory of presheaves of spaces is sent to the class of maps with respect to which we localize spectral presheaves. Following the description of the last given in the Remark \refnci{usingpresheavesofspectra2} it is enough to see that for any representable object $j(\X)$ we have $\Sigma^{\infty}_{+}j(\X)\simeq \delta_{j(\X)}(S)$ where the $S$ is the sphere spectrum. This is because $Map^{Sp}(-)$ is an internal-hom in $\widehat{\Sp}$ and the sphere spectrum is a unit for the monoidal structure. \\

In this section we will be considering the associated commutative diagram of right adjoints

\begin{equation}
\xymatrix{
Fun(\dg^{ft},\widehat{\Spaces})&&\ar[ll]^{\Omega^{\infty}}Fun(\dg^{ft},\widehat{\Sp})\\
Fun_{Nis}(\dg^{ft},\widehat{\Spaces})\ar@{^{(}->}[u]&&\ar[ll]^{\Omega^{\infty}_{Nis}}Fun_{Nis}(\dg^{ft},\widehat{\Sp})\ar@{^{(}->}[u]\\
Fun_{Nis,   L_{pe}(\mathbb{A}^1)}(\dg^{ft},\widehat{\Spaces})\ar@{^{(}->}[u]&&\ar[ll]_{\sim}^{\Omega^{\infty}_{Nis, \mathbb{A}^1}}Fun_{Nis,  L_{pe}(\mathbb{A}^1)}(\dg^{ft},\widehat{\Sp})\ar@{^{(}->}[u]
}
\end{equation}

\noindent where again the last map is an equivalence. We will now explain how to use this diagram to reduce the proof that $l^{nc}_{\mathbb{A}^1}(K^S)$ is unit for the monoidal structure in $Fun_{Nis,  L_{pe}(\mathbb{A}^1)}(\dg^{ft}, \widehat{\Sp})$ to the proof that $l_{0,\mathbb{A}^1}^{nc}(\Omega^{\infty}(K^c))$ is a unit for the monoidal structure in $Fun_{Nis,   L_{pe}(\mathbb{A}^1)}(\dg^{ft},\widehat{\Spaces})$. This will require some preliminaries. First we recall that thanks to the Prop. \ref{nc2proposicao1} we have an equivalence

\begin{equation}
\xymatrix{Fun_{Nis \geq 0}(\dg^{ft}, \widehat{\Sp}_{\geq 0})&\ar[l]_{\sim}^{\overline{\tau_{\geq 0}}}Fun_{Nis}(\dg^{ft}, \widehat{\Sp})}
\end{equation}

This equivalence provides a compatibility for the $\mathbb{A}^1$-localizations, in the sense that the diagram

\begin{equation}
\label{nc2DIAD6}
\xymatrix{
Fun_{Nis \geq 0}(\dg^{ft}, \widehat{\Sp}_{\geq 0})\ar[d]^{l_{\geq 0, \mathbb{A}^1}^{nc}}&\ar[l]_{\sim}^{\overline{\tau_{\geq 0}}}Fun_{Nis}(\dg^{ft}, \widehat{\Sp})\ar[d]^{l_{\mathbb{A}^1}^{nc}}\\
Fun_{Nis \geq 0,  L_{pe}(\mathbb{A}^1)}(\dg^{ft}, \widehat{\Sp}_{\geq 0})&\ar[l]_{\sim}^{\overline{\tau_{\geq 0}}}Fun_{Nis,   L_{pe}(\mathbb{A}^1)}(\dg^{ft}, \widehat{\Sp})
}
\end{equation}

\noindent commutes. Here $l_{\geq 0, \mathbb{A}^1}^{nc}$ is the (noncommutative) $\mathbb{A}^1$-localization functor for connectively-Nisnevich local presheaves.\\

The second preliminary result is a consequence of the equivalence between $\widehat{\Sp}_{\geq 0}$ and the $(\infty,1)$-category of commutative algebra objects $CAlg(\widehat{\Spaces})$ (see \cite[5.1.3.7]{lurie-ha}) and the equivalence of this last one with $Fun^{Segal-grplike}(N(\Fin),\widehat{\Spaces})$  - the full subcategory of the $(\infty,1)$-category $Fun(N(\Fin), \widehat{\Spaces})$ spanned by those functors satisfying the Segal condition and which are grouplike (see \cite[2.4.2.5]{lurie-ha}). This will be explained later in this paper in a proper way.

We can easily check that this equivalence induces equivalences

\begin{equation}
\label{nc2DIAD2}
Fun^{Segal-grplike}(N(\Fin),Fun_{Nis}(\dg^{ft},\widehat{\Spaces}))\simeq Fun_{Nis \geq 0}(\dg^{ft}, \widehat{\Sp}_{\geq 0})
\end{equation}

\noindent and 

\begin{equation}
\label{nc2DIAD3}
Fun^{Segal-grplike}(N(\Fin),Fun_{Nis,  L_{pe}(\mathbb{A}^1)}(\dg^{ft},\widehat{\Spaces}))\simeq Fun_{Nis \geq 0,  L_{pe}(\mathbb{A}^1)}(\dg^{ft}, \widehat{\Sp}_{\geq 0})
\end{equation}

\noindent and we claim that the $\mathbb{A}^1$-localization functor $l_{\geq 0, \mathbb{A}^1}^{nc}$ can be identified along this equivalence with the functor induced by the levelwise application of the $\mathbb{A}^1$-localization functor for spaces $l_{0, \mathbb{A}^1}^{nc}$. To confirm that this is indeed the case we observe first that the composition with $l_{0, \mathbb{A}^1}^{nc}$ produces a left-adjoint to the inclusion

\begin{equation}
Fun(N(\Fin), Fun_{Nis,  L_{pe}(\mathbb{A}^1)}(\dg^{ft}, \widehat{\Spaces}))\subseteq Fun(N(\Fin), Fun_{Nis}(\dg^{ft}, \widehat{\Spaces}))
\end{equation}

\noindent so that it suffices to check that this left-adjoint preserves Segal-grouplike objects. To prove this we will need an explicit descriptin for the  $\mathbb{A}^1$-localization for Nisnevich local objects functors $\dg^{ft}\to \widehat{\Spaces}$, in the same spirit of the description given in the previous section of this paper. Such an explicit model was already given for the commutative case in the foundational paper \cite{voevodsky-morel} where the authors prove (see the Lemma 1-3.20 and the Lemma 2 2.6) that the endofunctor of $\C=Fun(\aff^{op}, \widehat{\Spaces})$ given by the formula

\begin{equation}
F\mapsto colim_{i\in \mathbb{N}} U^i(F)
\end{equation}

\noindent with $U(F):= colim_{n\in \Delta^{op}}\underline{Hom}(\Delta_{ \mathbb{A}^1}^n, F)$, is an explicit model for the  $\mathbb{A}^1$-localization functor for presheaves of spaces. Furthermore, this formula preserves Nisnevich local objects.\footnote{In the original formulation of this result the authors use a different description of $U(F)$ that follows from the fact the the geometric realization of a simplicial space is homotopy equivalent to the diagonal of the underlying bisimplicial set.}

One can also check now that the same proof works in the noncommutative world replacing again $\Delta_{ \mathbb{A}^1}$ by its composition with $L_{pe}$. 

This description can now be used to prove that the composition with $l_{0, \mathbb{A}^1}^{nc}$ preserves the Segal-grouplike condition. Indeed, this follows immediately from this explicit description together with the fact that products in $Fun_{Nis}(\dg^{ft}, \Spaces)$ are computed objectwise in spaces and the fact that in spaces both sifted and filtered colimits commute with finite products (see \cite[5.5.8.11,5.5.8.12]{lurie-htt} for the sifted case). The grouplike condition follows from this and the functoriality of $l_{0, \mathbb{A}^1}^{nc}$. As a summary of this discussion, we concluded the existence of a commutative diagram

\begin{equation}
\label{nc2DIAD5}
\xymatrix{
Fun^{Segal-grplike}(N(\Fin),Fun_{Nis}(\dg^{ft},\widehat{\Spaces}))\ar[d]^{(l_{0, \mathbb{A}^1}^{nc}\circ -)}& \ar[l]_(0.35){\sim} Fun_{Nis \geq 0}(\dg^{ft}, \widehat{\Sp}_{\geq 0})\ar[d]^{l_{\geq 0, \mathbb{A}^1}^{nc}}\\
Fun^{Segal-grplike}(N(\Fin),Fun_{Nis,  L_{pe}(\mathbb{A}^1)}(\dg^{ft},\widehat{\Spaces}))&\ar[l]_(0.35){\sim} Fun_{Nis \geq 0,  L_{pe}(\mathbb{A}^1)}(\dg^{ft}, \widehat{\Sp}_{\geq 0})
}
\end{equation}

Finally, combining the commutativity of this diagram with the diagram (\ref{nc2DIAD6}) we obtain the commutativity of the diagram

\begin{equation}
\xymatrix{
Fun_{Nis}(\dg^{ft},\widehat{\Spaces})\ar[d]^{l_{0,\mathbb{A}^1}^{nc}}&&\ar[ll]^{\Omega^{\infty}_{Nis}}Fun_{Nis}(\dg^{ft},\widehat{\Sp})\ar[d]^{l_{\mathbb{A}^1}^{nc}}\\
Fun_{Nis,   L_{pe}(\mathbb{A}^1)}(\dg^{ft},\widehat{\Spaces})&&\ar[ll]_{\sim}^{\Omega^{\infty}_{Nis, \mathbb{A}^1}}Fun_{Nis,  L_{pe}(\mathbb{A}^1)}(\dg^{ft},\widehat{\Sp})
}
\end{equation}

This follows because $\Omega_{Nis}^{\infty}$ can now be identified with the evaluation at $\onefin\in N(\Fin)$ by means of the equivalence (\ref{nc2DIAD3}) and because as we conclude above, the $\mathbb{A}^1$-localization formula for connective spectra is determined levelwise by the formula for spaces.

The following lemma is the last step in our preliminaries:

\begin{lemma}
Let $F$ be a connectively-Nisnevich local object in $Fun(\dg^{ft},\widehat{\Sp})$. Then, $\Omega^{\infty}(F)$ is Nisnevich local and the canonical map $\Omega^{\infty}(F)\simeq l^{nc}_{0, Nis}(\Omega^{\infty}(F))\to \Omega_{Nis}(l^{nc}_{Nis}(F))$ is an equivalence in $Fun_{Nis}(\dg^{ft}, \widehat{\Spaces})$.
\begin{proof}
The proof relies on two observations. The first is that if $F$ is connectively-Nisnevich local, the looping $\Omega^{\infty}(F)$ is Nisnevich local as a functor $\dg^{ft}\to \widehat{\Spaces}$. This is because the composition $\xymatrix{\widehat{\Sp}_{\geq 0}\ar@{^{(}->}[r]& \widehat{\Sp}\ar[r]^{\Omega^{\infty}}& \widehat{\Spaces}}$ preserves limits (one possible way to see this is to use the equivalence between connective spectra and grouplike commutative algebras in $\widehat{\Spaces}$ for the cartesian product \cite[Theorem 5.1.3.16 and Remark 5.1.3.17]{lurie-ha} and the fact that this equivalence identifies the looping functor $\Omega^{\infty}$ with the forgetful functor which we know as a left adjoint and therefore commutates with limits. The conclusion now follows from the definition of connectively-Nisnevich local. The second observation is that the looping functor $\Omega^{\infty}$ only captures the connective part of a spectrum. This follows from the very definition of the canonical $t$-structure in $\widehat{\Sp}$ (see \cite[1.4.3.4]{lurie-ha}) In particular, since $F$ is connectively-Nisnevich local, our Proposition \ref{nc2proposicao1} implies that the canonical morphism $F\to l^{nc}_{Nis}(F)$ is an equivalence in the connective part so that its image under $\Omega^{\infty}$ is an equivalence. Putting together these two observations we have equivalences fitting in a commutative diagram

\begin{equation}
\xymatrix{
\Omega^{\infty}(F)\ar[d]^{\sim}\ar[dr]^{\sim}&\\
l^{nc}_{0, Nis}(\Omega^{\infty}(F))\ar[r]^{\delta}& \Omega^{\infty}(l^{nc}_{Nis}(F))
}
\end{equation}

\noindent so that the canonical map $\delta$ induced by the universal property of the localization is also an equivalence.
\end{proof}
\end{lemma}

Finally, we uncover the formulas

\begin{equation}
\Omega_{Nis, \mathbb{A}^1}^{\infty}(l_{ \mathbb{A}^1}^{nc}(K^S))\simeq l_{0,  \mathbb{A}^1}^{nc}(\Omega_{Nis}^{\infty}(K^S))\simeq l_{0, \mathbb{A}^1}(\Omega^{\infty}(K^c))
\end{equation}

\noindent where the first equivalence follows from the preceeding discussion and the last one follows from the previous lemma.\\

The first task is done. Now we explain the equivalence between  $l_{0, \mathbb{A}^1}(\Omega^{\infty}(K^c))$ and the unit for the monoidal structure in $Fun_{Nis,  L_{pe}(\mathbb{A}^1)}(\dg^{ft}, \widehat{\Spaces})$.

Our starting point is the formula (\ref{nc2formulaktheory}) describing the $K$-theory space of an idempotent complete dg-category $T$ by means of a colimit of mapping spaces. Since colimits and limits of functors are determined objetwise, the functor $\Omega^{\infty}K^c$ can itself be written as $\Omega\, colim_{[n]\in \Delta^{op}}\, Seq$ where $Seq$ is the object in the $(\infty,1)$-category $Fun(\Delta^{op}, Fun(\dg^{idem},\widehat{\Spaces}))$ resulting from the last stage of the Construction \ref{Sconstruction1}.

\begin{remark}
\label{nc2sephiroth2}
More precisely, at the end the Construction \ref{Sconstruction1} we obtained a functor 

\begin{equation}
\label{nc2mapaultimodia}
N(Cat_{Ch(k)})\to Fun(N(\Delta^{op}), N(\widehat{\Delta}_{big}))\to Fun(N(\Delta^{op}),\widehat{\Spaces})
\end{equation}

\noindent where the second map is induced by the localization functor $N(\widehat{\Delta}_{big})\to \widehat{\Spaces}$ with $\widehat{\Delta}_{big}$ the very big category big of simplicial sets equipped with the standard model structure. By the description of each space at level $n$ as a mapping space we conclude that this composition sends Morita equivalences of dg-categories to equivalences and therefore by the universal property the localization extends to a unique functor $\dg^{idem}\to Fun(N(\Delta)^{op}, \Spaces)$ which, using the equivalence between   $Fun(\dg^{idem},Fun(N(\Delta^{op}), \widehat{\Spaces})$ and $Fun(N(\Delta^{op}), Fun(\dg^{idem}, \widehat{\Spaces})$ gives what we call $Seq$.
\end{remark}

The value of $Seq$ at zero is the constant functor with value $\ast$ and its value at $n\geq 1$ is $Map_{\dg^{idem}}(\widehat{([n-1]_k)}_c, -)$. The boundary and degeneracy maps are obtained from the $S$-construction as explained in the Construction \ref{Sconstruction1}. We observe now that the dg-categories $\widehat{([n-1]_k)}_c$, for any $n\geq 0$, are of finite type so that each level of the simplicial object $Seq$ is in the full subcategory of $\omega$-continuous functors. Moreover, we can think of the dg-categories $\widehat{([n]_k)}_c$ as non-commutative spaces $I_n$ so that by means of the Yoneda's map $j_{nc}: \nck\hookrightarrow Fun(\dg^{ft}, \Spaces)$ we can identify $Seq_n$ with the representable $Map_{\nck}(-, I_{n-1})$. In particular, since the Yoneda's map is fully-faithfull, the simplicial object Seq is the image through $j_{nc}$ of a uniquely determined simplicial object $Seq_{nc}\in Fun(\Delta^{op}, \nck)$ whose value at level $n$ is the noncommutative space $I_{n-1}$. Finally, with these notations we can write $\Omega^{\infty}K^c$ as $\Omega\, colim_{[n]\in \Delta^{op}}\, j_{nc}\circ Seq_{nc}$ so that our main goal is to understand the localization $ l^{nc}_{0,\mathbb{A}^1}(\Omega\, colim_{[n]\in \Delta^{op}}\, j_{nc}\circ Seq_{nc})$. As the zero level of the simplicial object  $j_{nc}\circ Seq_{nc}$ is contractible, the realization $colim_{[n]\in \Delta^{op}}\, j_{nc}\circ Seq_{nc}$ is $1$-connective \footnote{Recall that a space is said to be $n$-connective if it is non-empty and all its homotopy groups for $i< n$ are zero.}. We have the following general fact:

\begin{lemma}
\label{Alocalizationpreservesloopings}
For any $F\in Fun(\dg^{ft}, \widehat{\Spaces})$ with values in $1$-connective spaces, the canonical map  $l^{nc}_{0,\mathbb{A}^1}(\Omega(F))\to  \Omega(l^{nc}_{0,\mathbb{A}^1}(F))$ is an equivalence.
\begin{proof}
This follows from the explicit description for $l_{0, \mathbb{A}^1}^{nc}$ given above as a filtered colimit of a sifted colimits, together with the fact that $\Omega$ preserves filtered colimits (this follows because homotopy groups preserve filtered colimits) and sifted colimits of $1$-connective spaces (see the proof of \cite[1.4.3.9]{lurie-ha}) Moreover, we use the fact that $1$-connective spaces are stable under sifted colimits (see again the arguments in the proof of \cite[1.4.3.9]{lurie-ha}).
\end{proof}
\end{lemma}

In particular, we have

\begin{equation}
 l^{nc}_{0,\mathbb{A}^1}(\Omega\, colim_{[n]\in \Delta^{op}}\, j_{nc}\circ Seq_{nc})\simeq \Omega\, l^{nc}_{0,\mathbb{A}^1}( colim_{[n]\in \Delta^{op}}\, j_{nc}\circ Seq_{nc}) 
\end{equation}

Our main goal now is to understand the simplicial object $Seq$. Following Waldhausen \cite{waldhausen-ktheoryofspaces} we recall the existence of a weaker version of the $S$-construction that considers only those sequence of cofibrations that split. More precisely, and using the same terminology as in the Construction \ref{Sconstruction1} we denote by $\mathbb{R}\underline{Hom}^{split}(Ar[n]_k, \widehat{T}_c)$ the full sub dg-category of $\mathbb{R}\underline{Hom}(Ar[n]_k, \widehat{T}_c)$ spanned by those $Ar[n]$-indexed diagrams satisfying the conditions given in the Construction \ref{Sconstruction0} and where the top sequence is given by the canonical inclusions $E_1\to E_1\oplus E_2\to E_1\oplus E_2\oplus E_3\to ...\to E_1\oplus...\oplus E_n$ for some list of perfect modules $(E_1,..., E_n)$. These are called \emph{split cofibrations}. As in the standard $S$-construction, the  categories subjacent to $\mathbb{R}\underline{Hom}^{split}(Ar[n]_k, \widehat{T}_c)$ carries a notion of weak-equivalences $W_n^{Split}$ and assemble to form a simplical space $[n]\to N( \mathbb{R}\underline{Hom}^{split}(Ar[n]_k, \widehat{T}_c)^{W_n^{Split}})$.

As in the Construction \ref{Sconstruction1} we can now describe these spaces in a somewhat more simple form. As the dg-categories $1_k$ are cofibrant (see \cite{tabuada-quillen}) they are also locally-cofibrant and for any $n\geq 0$ the coproduct $\coprod_{i=1}^n1_k$ is an homotopy coproduct. Moreover, for any locally-cofibrant dg-category $T$ we have equivalences $\mathbb{R}\underline{Hom}(\coprod_{i=1}^n1_k,\widehat{T}_c)\simeq \prod_{i=1}^n\widehat{(1_k \otimes^\mathbb{L} T)}_{pspe}\simeq  \prod_{i=1}^n\widehat{(1_k \otimes T)}_{pspe}\simeq  \prod_{i=1}^n\widehat{T}_{c}$ In this case, for every $n\geq 0$ and for every dg-category $T$ there is an equivalence between the category subjacent to $\mathbb{R}\underline{Hom}^{split}(Ar[n]_k, \widehat{T}_c)$ and the category subjacent to  $\mathbb{R}\underline{Hom}(\coprod_{i=1}^n 1_k, \widehat{T}_c)$, defined by sending a sequence $E_1\to E_1\oplus E_2\to E_1\oplus E_2\oplus E_3\to ...\to E_1\oplus...\oplus E_n$ to the sucessive quotients  $(E_1,..., E_n)$. This correspondence is functorial and defines an equivalence because of the universal property of direct sums. Moreover, and again thanks to the cube lemma, this equivalence preserves the natural notins of weak-equivalences. Finally, and again due to the main theorem of \cite{Toen-homotopytheorydgcatsandderivedmoritaequivalences} we found the spaces $N( \mathbb{R}\underline{Hom}^{split}(Ar[n]_k, \widehat{T}_c)^{W_n^{Split}})$ and $Map_{\dg^{idem}}(\oplus_{i=1}^n\widehat{(1_k)}_c, \widehat{T}_c)$ to be equivalent so that by the same arguments as in the Remark \ref{nc2sephiroth2} we obtain a simplicial object $Split\in Fun(N(\Delta)^{op}, Fun(\dg^{idem}, \widehat{\Spaces})$, which, because the dg-categories $\oplus_{i=1}^n\widehat{(1_k)}_c$ are of finite type, lives in the full subcategory of $\omega$-continuous functors, therefore being an object in $Fun(N(\Delta^{op}), \mathcal{P}(\nck))$. Moreover, for each $n\geq 0$ $Split_n$ is representable by the noncommutative space associated to the dg-category $\oplus_{i=1}^n\widehat{(1_k)}_c$ so that by Yoneda the whole simplicial object $Split$ is of the form $j_{nc}\circ \Theta$ for a simplicial object $\Theta\in Fun(N\Delta^{op}), \nck)$ with level $n$ given by  $\oplus_{i=1}^n\widehat{(1_k)}_c$.

Finally, the inclusion of split cofibrations into all sequences of morphisms provides a strict map of simplicial objects in the model category $\widehat{\Delta}$ between $[n]\to N(\mathbb{R}\underline{Hom}^{split}(Ar[n]_k, \widehat{T}_c)^{W_n})$ and $[n]\to  N(S_n^{dg}(T)^{W_n})$ and we define $\lambda$

\begin{equation}
\lambda:j_{nc}(\Theta)\simeq Split\to Seq 
\end{equation}

\noindent to be the image of this map under the composition in (\ref{nc2mapaultimodia}). This is where the result of \cite{Anthony-thesis} becomes crutial:

\begin{prop}\cite[Prop. 4.6]{Anthony-thesis}
The map $\lambda$ is a levelwise noncommutative $\mathbb{A}^1$-equivalence in $Fun(\dg^{ft},\Spaces)$.
\begin{proof}
In \cite[Prop. 4.6]{Anthony-thesis} the author uses and inductive argument to prove that for any $n\geq 0$ the map $\lambda_n$ is  an $\mathbb{A}^1$-equivalence. For the induction step we can use exactly the same argument but the induction basis $n=2$ requires further adaptation to our case. Namely, we are required to construct a noncommutative $\mathbb{A}^1$-homotopy between the identity of the noncommutative space $I_{2-1}$ and the zero map. Such an homotopy corresponds to a co-homotopy in $\dg^{idem}$, namely, a map $H:\widehat{([1]_k)}_c\to \widehat{([1]_k)}_c \otimes^{\mathbb{L}}L_{pe}(\mathbb{A}^1)$ in $\dg^{idem}$ fitting in a commutative diagram

\begin{equation}
\xymatrix{
&&\widehat{([1]_k)}_c\\
\widehat{([1]_k)}_c\ar[r]^(0.4){H}\ar@/^1pc/[rru]^{Id} \ar@/_1pc/[drr]_0&\widehat{([1]_k)}_c\otimes^{\mathbb{L}}L_{pe}(\mathbb{A}^1)\ar[ur]^{ev_1}\ar[dr]_{ev_0}&\\
&&\widehat{([1]_k)}_c
}
\end{equation}

Recall that $L_{pe}(\mathbb{A}^1)$ is canonically equivalent to $\widehat{k[X]}_c$ - the idempotent completion of the dg-category with one object and $k[X]$ concentrated in degree zero as endomorphisms. In this case the term in the middle is equivalent to $\widehat{(([1]_k)\otimes k[X])}_c$. We define $H$ to be the map induced by the universal property of the idempotent completion $\widehat{(-)}_c: \dg\to \dg^{idem}$ by means of the composition

\begin{equation}
([1]_k)\to ([1]_k)\otimes k[X]\subseteq \widehat{(([1]_k)\otimes k[X])}_c
\end{equation}

\noindent where the first map is obtained from the strict dg-functor defined by the identity on the objects, by the inclusion $k\subseteq k[X]$ on the endomorphisms of $0$ and by the composition $k\subseteq k[X]\to k[X]$ on the complex of maps between $0$ and $1$ and on the endomorphisms of $1$, where the last map is the multiplication by the variable $X$. This makes the diagram above commute and provides the required homotopy.
\end{proof}
\end{prop}

Finally, the fact that any colimit of $\mathbb{A}^1$-equivalences is an $\mathbb{A}^1$-equivalence gives us the following corollary:

\begin{cor}
\label{nc2sephiroth1}
The map induced by $\lambda$ between the colimits $colim_{\Delta^{op}}\, j_{nc}\circ \Theta\to  colim_{\Delta^{op}}\, j_{nc}\circ Seq_{nc}$ is an $\mathbb{A}^1$-equivalence. Moreover, and since $l_{0,\mathbb{A}^1}^{nc}$ commutes with colimits and representable objects are Nisnevich local, we have equivalences

\begin{equation}
colim_{\Delta^{op}}\, l_{0,\mathbb{A}^1}^{nc}\circ  j_{nc}\circ \Theta\simeq l_{0,\mathbb{A}^1}^{nc}(colim_{\Delta^{op}}\, j_{nc}\circ \Theta)\simeq  l_{0,\mathbb{A}^1}^{nc}(colim_{\Delta^{op}}\, j_{nc}\circ Seq_{nc})\simeq colim_{\Delta^{op}}\, l_{0,\mathbb{A}^1}^{nc}\circ  j_{nc}\circ Seq_{nc} 
\end{equation}

\noindent  in $\stnck$.
\end{cor}

Our next move requires a small preliminary digression. To start with, recall that any $(\infty,1)$-category endowed with finite sums and an initial object or finite products and a final object, can be considered as a symmetric monoidal $(\infty,1)$-category with respect to these two operations, respectively denoted as $\C^{\coprod}$ and $\C^{\times}$ (see \cite[Sections 2.4.1 and 2.4.3]{lurie-ha}). Monoidal structures appearing from this mechanism are called, respectively, \emph{cartesian} and \emph{cocartesian}. In particular, if $\C$ has direct sums and a zero object, these monoidal structures coincide $\C^{\oplus}$ (this follows from the Proposition \cite[2.4.3.19]{lurie-ha}). In this particular situation the theory of algebras over a given $\infty$-operad $\Opmonoidal$ gets simplified: the $(\infty,1)$-category of $\Op$-algebras on $\C^{\oplus}$ is equivalent to a full subcategory of $Fun(\Opmonoidal,\C)$, spanned by a class of functors satisfying a certain Segal condition (see \cite[2.4.2.1, 2.4.2.5]{lurie-ha}). In the particular case of associative algebras, and since the category $\Delta^{op}$ is a "model" for the associative operad (see \cite[4.1.2.6, 4.1.2.10, 4.1.2.14]{lurie-ha} for the precise statement) an associative algebra in $\C^{\oplus}$ is just a simplicial object in $\C$ satisfying the Segal condition. 

\begin{equation}
\label{algebrasandsimplicialobjects}
Alg_{\mathcal{A}ss}(\C)\simeq Fun^{Segal}(N(\Delta^{op}), \C)
\end{equation}

We shall now come back to our situation and observe that

\begin{lemma}
The simplicial object $\Theta$ satisfies the Segal conditions. 
\begin{proof}
As the Yoneda's embedding preserves limits and is fully-faithfull it suffices to check that $Split$ satisfies the Segal conditions. But this is obvious from the definition of the simplicial structure given by the $S$-construction. At each level the map appearing in the Segal condition is the map sending a sequence of dg-modules $E_0\to E_0\oplus E_1\to ....E_0\oplus...\oplus E_{n-1}$ to the quotients $(E_0,..., E_{n-1})$.
\end{proof}
\end{lemma}

We now characterize the simplicial object $\Theta$ in a somewhat more canonical fashion. An important aspect of a cocartesian symmetric monoidal structure $\C^{\coprod}$ is that any object $X$ in $\C$ admits a unique algebra structure, determined by the codiagonal map $X\coprod X\to X$. More precisely (see \cite[2.4.3.16]{lurie-ha} for the general result), the forgetful map $Alg(\C)\to \C$ is an equivalence of $(\infty,1)$-categories \footnote{Recall that the associative operad is unital.}. By choosing an inverse to this equivalence and composing with the equivalence (\ref{algebrasandsimplicialobjects}) we obtain an $\infty$-functor

\begin{equation}
\label{simplicialobjectcodiagonal}
\C\to Alg_{\mathcal{A}ss}(\C)\simeq Fun^{Segal}(N(\Delta)^{op}, \C) 
\end{equation}

\noindent providing for any object in $\C$ a uniquely determined simplicial object, encoding the algebra structure induced by the codigonal\footnote{The fact that the multiplication can be identified with the codiagonal map follows from the simplicial identities and from the universal property defining the codiagonal.}. Because of the Segal condition this simplicial object is a zero object of $\C$ in degree zero, $X$ in degree one and more generally is $X^{\oplus_n}$ in degree $n$. We now apply this discussion to $\C=\nck$ (it has direct sums and a zero object because $\dg^{idem}$ has and the inclusion $\dg^{ft}\subseteq \dg^{idem}$ preserves them) and to $X=\widehat{(1_k)}_c$ (see the Notation \ref{notationunit}). Since the simplicial object $\Theta$ satisfies the Segal condition and its first level is equivalent to $X$, the equivalence (\ref{simplicialobjectcodiagonal}) tells us that it is necessarily the simplicial object codifying the unique associative algebra structure on $X$ given by the codiagonal.\\

With the Corollary \ref{nc2sephiroth1} we are now reduced to study the colimit of the simplicial object  $l_{\mathbb{A}^1}^{nc}\circ  j_{nc}\circ \Theta$ in $\stnck$. As the last is a stable $(\infty,1)$-category it has direct sums and therefore can be understood as the underlying $(\infty,1)$-category of a symmetric monoidal structure $\stnck^{\oplus}$ which is simultaneously cartesian and cocartesian. As the canonical composition $\nck\to \stnck$ preserves direct sums (this follows from $1)$ the fact the Yoneda functor preserves limits; $2)$ the fact representables are Nisnevich local; $3)$ the fact the $\mathbb{A}^1$-localization preserves finite products by the same arguments of the lemma \ref{Alocalizationpreservesloopings} and finally $4)$ the fact that $\stnck$ is stable.) it can be lifted in a essentially unique way to a monoidal functor $\nck^{\oplus}\to \stnck^{\oplus}$ (\cite[Cor. 2.4.1.8]{lurie-ha}). This monoidal map allows us to transport algebras and provides a commutative diagram

\begin{equation}
\xymatrix{
\ar@/_5pc/[dd]_(0.6){ev_{[1]}}Fun^{Segal}(N(\Delta)^{op}, \nck)\ar[r]\ar[d]^{\sim}& Fun^{Segal}(N(\Delta)^{op}, \stnck)\ar[d]^{\sim} \ar@/^5pc/[dd]^(0.6){ev_{[1]}}\\
Alg_{\mathcal{A}ss}(\nck)\ar[r]\ar[d]^{\sim}&\ar[d]^{\sim} Alg_{\mathcal{A}ss}(\stnck)\\
\nck\ar[r]& \stnck
}
\end{equation}

\noindent where the upper map is the composition with $\nck\to \stnck$. It follows from the description of $\Theta$ above and from the commutativity of this diagram that the simplicial object $l_{\mathbb{A}^1}^{nc}\circ  j_{nc}\circ \Theta$ in $\stnck$ corresponds to the unique commutative algebra structure on $1_{nc}:=l_{\mathbb{A}^1}^{nc}\circ  j_{nc}(L_{pe}(k))$ created by the codiagonal.\\

Our next task is to study the theory of associative algebras on a stable $(\infty,1)$-category equipped with its natural simultaneously cartesian and cocartesian monoidal structure induced by the existence of direct sums. We recall some terminology. If $\C^{\otimes}$ is a cartesian symmetric monoidal structure, an associative algebra on $\C$ is said to be \emph{grouplike} if the simplicial object which codifies it $A\in Fun^{Segal}(N(\Delta^{op}), \C)$ is a groupoid object in $\C$ in the sense of the definition \cite[6.1.2.7]{lurie-htt}. We let $Alg_{\Ass}^{grplike}(\C)$ denote the full subcategory of $Alg_{\mathcal{A}ss}(\C)$ denote the full subcategory spanned by the grouplike associative algebras.\\

Let now $\Delta^{op}_{+}$ be the standard augmentation of the category $\Delta^{op}$. Following \cite[6.1.2.11]{lurie-htt}, an object $U_+\in Fun(\Delta^{op}_{+},\C)$ is said to be a \emph{Cech nerve of the morphism $U_0\to U_{-1}$} if the restriction $U_+|_{N(\Delta^{op})}$ is a groupoid object and the commutative diagram

\begin{equation}
\xymatrix{
U_1\ar[r]\ar[d]& U_0\ar[d]\\
U_0\ar[r]& U_{-1}
}
\end{equation} 

\noindent is a pullback diagram in $\C$. Again by \cite[6.1.2.11]{lurie-htt}, a Cech nerve $U_+$ is determined by the map $U_0\to U_{-1}$ in a essentially unique way as the right-Kan extension along the inclusion $N(\Delta^{op}_{+, \leq 0})\subseteq N(\Delta^{op}_{+})$.\\

We have the following lemma:

\begin{lemma}
\label{nc2DIAD1}
Let $\C^{\otimes}$ be a cartesian symmetric monoidal $(\infty,1)$-category whose underlying $(\infty,1)$-category is stable. Then 
\begin{enumerate}
\item The inclusion $Alg^{grplike}(\C)\subseteq Alg(\C)$ is an equivalence;
\item For any object $X$ in $\C$ the simplicial object associated to $X$ by means of the composition (\ref{simplicialobjectcodiagonal}) is a Cech nerve of the canonical morphism $0\to \Sigma X$.
\end{enumerate}
\begin{proof}

The first assertion is true because in any stable $(\infty,1)$-category every morphism $f:X\to Y$ has an inverse $-f$ with respect to the additive structure. In particular, for any object $X\in \C$ there is  map $-Id_{X}$ providing an inverse for the algebra structure given by the codiagonal map $X\oplus X\to X$. More precisely, let $X$ be an object in $\C$ and let $U_X$ be the simplicial object associated to $X$ by means of the mechanism (\ref{simplicialobjectcodiagonal}). By construction this simplicial object satisfies the Segal condition and in particular we have $(U_X)_0\simeq 0$ and $(U_X)_1\simeq X$. We aim to prove that this simplicial object is a groupoid object. 
For that we observe that for a simplicial object $A$ to be a groupoid object it is equivalent to ask for $A$ to satisfy the Segal conditions and to ask for the induced map

\begin{equation}
\label{grpcondition}
\xymatrix{
A([2])\ar[rr]^{A(\partial_1)\times A(\partial_0)}&& A([1])\times A([1])
}
\end{equation}

\noindent to be an equivalence. Indeed, if $A$ is a groupoid object, by the description in \cite[6.1.2.6 - (4'')]{lurie-htt} it satisfies these two requirements automatically. The converse follows by applying the same arguments as in the proof of \cite[6.1.2.6 - 4)' implies 3)]{lurie-htt}, together with the observation that for the induction step to work we don't need the full condition in $4')$ but only the Segal condition. The induction basis is equivalent to the Segal conditions for $n=2$ together with the condition that (\ref{grpcondition}) is an equivalence.

In our case (\ref{grpcondition}) is the map $\nabla\times id_X: X\oplus X\to X\oplus X$ where $\nabla$ is the codigonal map $X\oplus X\to X$. Of course, since the identity of $X$ admits an inverse $(-Id_{X})$ the map $(\nabla\circ (Id_X\times (-Id_X)))\times Id_X$ is an explicit inverse for $\nabla\times id_X$.\\

Let us now prove $2)$. Again by construction, we know that the colimit of the truncation $(U_X)_{|_{N(\Delta_{\leq 1})^{op}}}$ is canonically equivalent to the suspension $\Sigma X$. Therefore $U_X$ admits a canonical augmentation $(U_X)^{+}: N(\Delta_{+}^{op})\to \C$ with  $(U_X)^{+}_{-1}\simeq \Sigma X$. It follows from $1)$ that $U$ is a groupoid object and since $\C$ is stable, the diagram

\begin{equation}
\xymatrix{
(U_X)_1\simeq X\ar[r]\ar[d]& (U_X)_0\simeq 0\ar[d]\\
(U_X)_0\simeq 0\ar[r]& (U_X)_{-1}\simeq \Sigma X
}
\end{equation} 
 
\noindent is a pullback so that $(U_X)^{+}$ is the Cech nerve of the canonical map $0\to \Sigma X$.
\end{proof}
\end{lemma}

In particular, we find that the simplicial object $l_{\mathbb{A}^1}^{nc}\circ  j_{nc}\circ \Theta$ is a Cech nerve of the canonical map $0\to \Sigma 1_{nc}$. Finally, recall that a morphism $A\to B$ is said to be an effective epimorphism if the colimit of its Cech nerve is $B$. The following lemma holds the final step

\begin{lemma}
Let $\C$ be a stable $(\infty,1)$-category. Then, for any object $X$ in $\C$, the canonical morphism $0\to X$ is an effective epimorphism.
\begin{proof}
Let $U:N(\Delta^{op})\to \C$ be a simplicial object in $\C$. Then the colimit of $U$ can be computed as the sequential colimit of the sucessive colimits of its truncations $U_{|_{N(\Delta^{op}_{\leq n})}}$. Using the descriptions of Cech nerves as right-Kan extensions (see above) we know that if $U^{+}$ is the Cech nerve of the map $0\to X$, its level $n$ is given by the $n$-fold tensor product of $0$ over $X$. As $\C$ is stable this $n$-dimensional limit cube will also be a colimit $n$-cube so that the colimit of the truncation at level $n$ will necessarily be $X$ (See the Proposition \cite[1.2.4.13]{lurie-ha}). Since this holds for every $n\geq 0$ the colimit of the Cech nerve is necessarily canonically equivalent to $X$.
\end{proof}
\end{lemma}

We are done. Since $\stnck$ is stable we have $colim_{\Delta^{op}}\, l_{0,\mathbb{A}^1}^{nc}\circ  j_{nc}\circ \Theta\simeq \Sigma 1_{nc}$ so that  $ l^{nc}_{0,\mathbb{A}^1}(\Omega^{\infty}(K^c))$ is equivalent to $\Omega\Sigma 1_{nc}\simeq 1_{nc}$.

\newpage

\appendix

\section{Appendix: On the Comparison with Cisinski-Tabuada's approach to Noncommutative Motives}
\label{tabuada-comparison}

This section is independent of the rest of the paper. Our goal here is to explain the relation between our approach to noncommutative motives and the approach already studied by G. Tabuada in \cite{tabuada-higherktheory,MR2986869} and Cisinski-Tabuada in \cite{tabuada-cisinski, MR2822869}. Both theories have the $(\infty,1)$-category $Fun_{\omega}(\dg^{idem}, \widehat{\Sp})$ as a common ground. To start with we observe that our version $\stnck$ can be identified with the full subcategory spanned by those functors $F$ sending Nisnevich squares of dg-categories to pullback-pushout squares in spectra and satisfying $\mathbb{A}^1$-invariance. Indeed, our original definition of $\stnck$ as a localization of $Fun(\dg^{ft}. \widehat{\Sp})$ can be transported along the equivalence
  
\begin{equation}
\label{nc2indequi}
Fun_{\omega}(\dg^{idem},\widehat{\Sp})\simeq  Fun(\dg^{ft}, \widehat{\Sp})
\end{equation}

\begin{remark}
We give a more precise description of this localization. Given a noncommutative smooth space $\X$ associated with a dg-category of finite type $T_{\X}$, the image of $\X$ in $Fun(\dg^{ft},\widehat{\Sp})$ under the spectral Yoneda's embedding is just the corepresentable $\Sigma^{\infty}_{+}Map_{\dg^{ft}}(T_{\X}, -)$. Moreover, since $\dg^{ft}$ is the full subcategory of compact objects in $\dg^{idem}$ (\refnci{dgideminddgft}) the image of this corepresentable under the equivalence (\ref{nc2indequi}) is the corepresentable $\Sigma^{\infty}_{+}Map_{\dg^{idem}}(T_{\X}, -)$. We consider $Fun_{\omega, Nis}(\dg^{idem}, \widehat{\Sp})$ the reflexive accessible localization of $Fun_{\omega}(\dg^{idem}, \widehat{\Sp})$ obtained by inverting the small set of all maps

\begin{equation}
\Sigma^{\infty}_{+}Map_{\dg^{idem}}(T_{\V}, -)\coprod_{\Sigma^{\infty}_{+}Map_{\dg^{idem}}(T_{\W}, -)} \Sigma^{\infty}_{+}Map_{\dg^{idem}}(T_{\UU}, -) \to \Sigma^{\infty}_{+}Map_{\dg^{idem}}(T_{\X}, -)
\end{equation}

\noindent induced by the Nisnevich squares 

\begin{equation}
\xymatrix{
T_{\X}\ar[d]\ar[r]& T_{\UU}\ar[d]\\
T_{\V}\ar[r]&T_{\W}
}
\end{equation}

\noindent of dg-categories as described in the discussion following the definition in \refnci{defncnisnevich}. By the theorem \cite[5.5.4.15]{lurie-htt} this is an accessible reflexive localization of $Fun_{\omega}(\dg^{idem}, \widehat{\Sp})$ and now by construction  the local  objects are those functors $F: \dg^{idem}\to \widehat{\Sp}$ commuting with filtered colimits and sending the classical Nisnevich squares of dg-categories of finite type to pullback-pushout squares. To conclude this discussion we remark that the existing left adjoint to the inclusion $Fun_{\omega, Nis}(\dg^{idem}, \widehat{\Sp})\subseteq Fun_{\omega}(\dg^{idem}, \widehat{\Sp})$ fits in a commutative diagram

\begin{equation}
\xymatrix{
Fun_{\omega}(\dg^{idem},\widehat{\Sp})\ar[d]\ar[r]^{\sim}& Fun(\dg^{ft}, \widehat{\Sp})\ar[d]^{l_{Nis}}\\
Fun_{\omega,Nis}(\dg^{idem},\widehat{\Sp})\ar[r]^{\sim}& Fun_{Nis}(\dg^{ft}, \widehat{\Sp})
}
\end{equation}

We can now proceed in analogous terms and localize with respect to $\mathbb{A}^1$ to obtain our new description of $\stnck$.\\
\end{remark}

Tabuada's approach focuses on the full subcategory $Fun_{\omega, Loc}(\dg^{idem}, \widehat{\Sp})$ spanned by those functors sending exact sequences of dg-categories to fiber/cofiber sequences in spectra. His main theorem is the existence of a stable presentable $(\infty,1)$-category which we denote here as $\M^{Tab}_{Loc}$, together with a functor $\dg^{idem}\to  \M^{Tab}_{Loc}$ preserving filtered colimits, sending exact sequences to fiber/cofiber sequences and universal in this sense. We can also easily see that $\M^{Tab}_{Loc}$ is a stable presentable symmetric monoidal $(\infty,1)$-category with the monoidal structure extending the monoidal structure in $\dg^{idem}$. This  result was originally formulated using the language of derivators (see \cite{maltsiniotis-derivators} for an introduction) but we can easily extend it to the setting of $(\infty,1)$-categories by applying the same construction and the general machinery developed by J. Lurie in \cite{lurie-ha, lurie-htt}. In particular we have an equivalence of $(\infty,1)$-categories 

\begin{equation}
Fun_{\omega, Loc}(\dg^{idem}, \widehat{\Sp})\simeq Fun^{L}(\M^{Tab}_{Loc}, \widehat{\Sp})
\end{equation}

As we can see this is a theorem about a specific class of objects inside $Fun_{\omega}(\dg^{idem}, \widehat{\Sp})$, namely, those that satisfy localization. The comparison with our approach starts with the observation that any object $F$ satisfying localization satisfies also our condition of Nisnevich descent so that we have an inclusion of full subcategories $Fun_{\omega, Loc}(\dg^{idem}, \widehat{\Sp})\subseteq Fun_{\omega, Nis}(\dg^{idem}, \widehat{\Sp})$. In particular, we can identify $Fun_{\omega, Loc, \mathbb{A}^1}(\dg^{idem}, \widehat{\Sp})$ with a full subcategory of $\stnc$. We summarize this in the following diagram

\begin{equation}
\xymatrix{
&Fun_{\omega}(\dg^{idem}, \widehat{\Sp})&\\
 Fun_{\omega, Nis}(\dg^{idem}, \widehat{\Sp})\ar@{^{(}->}[ur]&&\ar@{_{(}->}[ll] \ar@{_{(}->}[ul]Fun_{\omega, Loc}(\dg^{idem}, \widehat{\Sp})\simeq Fun^{L}(\M^{Tab}_{Loc}, \widehat{\Sp})\\
 Fun_{\omega, Nis, \mathbb{A}^1}(\dg^{idem}, \widehat{\Sp})=: \stnck\ar@{^{(}->}[u]&&\ar@{_{(}->}[ll]\ar@{^{(}->}[u]Fun_{\omega, Loc, \mathbb{A}^1}(\dg^{idem}, \widehat{\Sp})
}
\end{equation}

The second observation is that the construction  $\M^{Tab}_{Loc}$ of Tabuada admits analogues adapted to each of  the full subcategories in this diagram. More precisely  one can easily show the existence of new stable presentable symmetric monoidal  $(\infty,1)$-categories $\M^{Tab}_{Nis}$, $\M^{Tab}_{Nis, \mathbb{A}^1}$, $\M^{Tab}_{Loc, \mathbb{A}^1}$ all equipped with $\omega$-continuous monoidal functors from $\dg^{idem}$, universal with respect to each of the obvious respective properties. In particular we find an equivalence

\begin{equation}
Fun^{L}(\M^{Tab}_{Nis, \mathbb{A}^1}, \widehat{\Sp})\simeq \stnck
\end{equation}

\noindent exhibiting the duality between our approach and the corresponding Nisnevich-$\mathbb{A}^1$-version of Tabuada's construction (recall that the very big $(\infty,1)$-category of big stable presentable $(\infty,1)$-categories has a natural symmetric monoidal structure \cite[6.3.2.10, 6.3.2.18 and 6.3.1.17]{lurie-ha} where the big $(\infty,1)$-category of spectra $\widehat{\Sp}$ is a unit and $Fun^{L}(-,-)$ is the internal-hom).  Because Localizing descent forces Nisnevich descent, the universal properties involved provide a zig-zag of canonical colimit preserving monoidal functor $\M^{Tab}_{Nis,\mathbb{A}^1}\rightarrow  \M^{Tab}_{Loc, \mathbb{A}^1}\leftarrow \M^{Tab}_{Loc}$ relating the dual of our new theory to the localizing theory of Tabuada. As emphazised before, the main advantage (in fact, la \emph{raison-d'être}) of our approach to noncommutative motives is the easy comparison with the motivic stable homotopy theory of schemes. The duality here presented explains why the original approach of Cisinski-Tabuada is not directly comparable.

\bibliographystyle{abbrv}	
\bibliography{biblio}

\begin{thebibliography}{10}

\bibitem{1204.3607}
C.~Barwick.
\newblock On the algebraic k-theory of higher categories, 2012.

\bibitem{MR0249491}
H.~Bass.
\newblock {\em Algebraic {$K$}-theory}.
\newblock W. A. Benjamin, Inc., New York-Amsterdam, 1968.

\bibitem{beilinsonprojective}
A.~Beilinson.
\newblock {\em The derived category of coherent sheaves on $\mathbb{P}^n$}.

\bibitem{Anthony-thesis}
A.~Blanc.
\newblock {\em Topological K-theory and its Chern character for non-commutative
  spaces}.
\newblock arXiv:1211.7360.

\bibitem{tabuada-gepner}
A.~Blumberg, D.~Gepner, and G.~Tabuada.
\newblock {\em A Universal Characterization of Higher Algebraic K-Theory}.

\bibitem{MR0116022}
A.~Borel and J.-P. Serre.
\newblock Le th\'eor\`eme de {R}iemann-{R}och.
\newblock {\em Bull. Soc. Math. France}, 86:97--136, 1958.

\bibitem{MR513569}
A.~K. Bousfield and E.~M. Friedlander.
\newblock Homotopy theory of {$\Gamma $}-spaces, spectra, and bisimplicial
  sets.
\newblock In {\em Geometric applications of homotopy theory ({P}roc. {C}onf.,
  {E}vanston, {I}ll., 1977), {II}}, volume 658 of {\em Lecture Notes in Math.},
  pages 80--130. Springer, Berlin, 1978.

\bibitem{cisinski-descentpar}
D.-C. Cisinski.
\newblock Descente par éclatements en k-théorie invariante par homotopie.
\newblock In {\em Annals of Mathematics}, volume 177, pages 425--448. 2013.

\bibitem{MR2822869}
D.-C. Cisinski and G.~Tabuada.
\newblock Non-connective {$K$}-theory via universal invariants.
\newblock {\em Compos. Math.}, 147(4):1281--1320, 2011.

\bibitem{tabuada-cisinski}
D.-C. Cisinski and G.~Tabuada.
\newblock Symmetric monoidal structure on non-commutative motives.
\newblock {\em J. K-Theory}, 9(2):201--268, 2012.

\bibitem{MR1644323}
W.~Fulton.
\newblock {\em Intersection theory}, volume~2 of {\em Ergebnisse der Mathematik
  und ihrer Grenzgebiete. 3. Folge. A Series of Modern Surveys in Mathematics
  [Results in Mathematics and Related Areas. 3rd Series. A Series of Modern
  Surveys in Mathematics]}.
\newblock Springer-Verlag, Berlin, second edition, 1998.

\bibitem{gepner-algebraiccobordismalgebraicKtheory}
D.~Gepner and V.~Snaith.
\newblock On the motivic spectra representing algebraic cobordism and algebraic
  {$K$}-theory.
\newblock {\em Doc. Math.}, 14:359--396, 2009.

\bibitem{hovey-modelcategories}
M.~Hovey.
\newblock {\em Model categories}, volume~63 of {\em Mathematical Surveys and
  Monographs}.
\newblock American Mathematical Society, Providence, RI, 1999.

\bibitem{MR1695653}
M.~Hovey, B.~Shipley, and J.~Smith.
\newblock Symmetric spectra.
\newblock {\em J. Amer. Math. Soc.}, 13(1):149--208, 2000.

\bibitem{joyal-article}
A.~Joyal.
\newblock {\em Quasi-categories and Kan complexes}.
\newblock J. Pure Appl. Algebra 175 (2002), 2002, no. 1-3.

\bibitem{MR0233871}
M.~Karoubi.
\newblock Foncteurs d\'eriv\'es et {$K$}-th\'eorie. {C}at\'egories filtr\'ees.
\newblock {\em C. R. Acad. Sci. Paris S\'er. A-B}, 267:A328--A331, 1968.

\bibitem{keller-exact}
B.~Keller.
\newblock On the cyclic homology of exact categories.
\newblock {\em J. Pure Appl. Algebra}, 136(1):1--56, 1999.

\bibitem{kontsevich3}
M.~Kontsevich.
\newblock {\em Mixed noncommutative motives - Talk at the Workshop on
  Homological Mirror Symmetry. University of Miami. 2010. Notes available at
  www-math.mit.edu/auroux/frg/miami10-notes}.

\bibitem{kontsevich1}
M.~Kontsevich.
\newblock {\em Noncommutative motives. Talk at the Institute for Advanced Study
  on the occasion of the 61 st birthday of Pierre Deligne, October 2005. Video
  available at http://video.ias.edu/Geometry-and-Arithmetic}.

\bibitem{kontsevich2}
M.~Kontsevich.
\newblock {\em Triangulated categories and geometry - Course at the Ecole
  Normale Superieure, Paris, 1998. Notes available at
  www.math.uchicago.edu/mitya/langlands.html}.

\bibitem{lurie-htt}
J.~Lurie.
\newblock {\em Higher topos theory}, volume 170 of {\em Annals of Mathematics
  Studies}.
\newblock Princeton University Press, Princeton, NJ, 2009.

\bibitem{lurie-ha}
J.~Lurie.
\newblock {\em \htmladdnormallink{ Higher Algebra}{
  http://www.math.harvard.edu/~lurie/papers/higheralgebra.pdf}}.
\newblock August 2012.

\bibitem{maltsiniotis-derivators}
G.~Maltsiniotis.
\newblock {\em Introduction à la théorie des dérivateurs (d'après
  Grothendieck)}.

\bibitem{voevodsky-morel}
F.~Morel and V.~Voevodsky.
\newblock {${\bf A}^1$}-homotopy theory of schemes.
\newblock {\em Inst. Hautes \'Etudes Sci. Publ. Math.}, (90):45--143 (2001),
  1999.

\bibitem{1010.3944}
N.~Naumann, M.~Spitzweck, and P.~A. Østvær.
\newblock Existence and uniqueness of e-infinity structures on motivic k-theory
  spectra, 2010.

\bibitem{MR0338129}
D.~Quillen.
\newblock Higher algebraic {$K$}-theory. {I}.
\newblock In {\em Algebraic {$K$}-theory, {I}: {H}igher {$K$}-theories ({P}roc.
  {C}onf., {B}attelle {M}emorial {I}nst., {S}eattle, {W}ash., 1972)}, pages
  85--147. Lecture Notes in Math., Vol. 341. Springer, Berlin, 1973.

\bibitem{quillen}
D.~G. Quillen.
\newblock {\em Homotopical algebra}.
\newblock Lecture Notes in Mathematics, No. 43. Springer-Verlag, Berlin, 1967.

\bibitem{riou-spanierwhitehead}
J.~Riou.
\newblock {\em Spanier-Whitehead duality in algebraic geometry}.
\newblock 2004.

\bibitem{nc1}
M.~Robalo.
\newblock From commutative to noncommutative motives.
\newblock March 2013.

\bibitem{MR1930883}
M.~Schlichting.
\newblock A note on {$K$}-theory and triangulated categories.
\newblock {\em Invent. Math.}, 150(1):111--116, 2002.

\bibitem{schlichting-negative}
M.~Schlichting.
\newblock Negative {$K$}-theory of derived categories.
\newblock {\em Math. Z.}, 253(1):97--134, 2006.

\bibitem{MR0353298}
G.~Segal.
\newblock Categories and cohomology theories.
\newblock {\em Topology}, 13:293--312, 1974.

\bibitem{tabuada-invariantsadditifs}
G.~Tabuada.
\newblock Invariants additifs de {DG}-cat\'egories.
\newblock {\em Int. Math. Res. Not.}, (53):3309--3339, 2005.

\bibitem{tabuada-quillen}
G.~Tabuada.
\newblock Une structure de cat\'egorie de mod\`eles de {Q}uillen sur la
  cat\'egorie des dg-cat\'egories.
\newblock {\em C. R. Math. Acad. Sci. Paris}, 340(1):15--19, 2005.

\bibitem{tabuada-higherktheory}
G.~Tabuada.
\newblock Higher {$K$}-theory via universal invariants.
\newblock {\em Duke Math. J.}, 145(1):121--206, 2008.

\bibitem{MR2986869}
G.~Tabuada.
\newblock A guided tour through the garden of noncommutative motives.
\newblock In {\em Topics in noncommutative geometry}, volume~16 of {\em Clay
  Math. Proc.}, pages 259--276. Amer. Math. Soc., Providence, RI, 2012.

\bibitem{thomasonalgebraic}
R.~W. Thomason and T.~Trobaugh.
\newblock Higher algebraic {$K$}-theory of schemes and of derived categories.
\newblock In {\em The {G}rothendieck {F}estschrift, {V}ol.\ {III}}, volume~88
  of {\em Progr. Math.}, pages 247--435. Birkh\"auser Boston, Boston, MA, 1990.

\bibitem{Toen-homotopytheorydgcatsandderivedmoritaequivalences}
B.~To{\"e}n.
\newblock The homotopy theory of {$dg$}-categories and derived {M}orita theory.
\newblock {\em Invent. Math.}, 167(3):615--667, 2007.

\bibitem{toen-vaquie}
B.~To{\"e}n and M.~Vaqui{\'e}.
\newblock Moduli of objects in dg-categories.
\newblock {\em Ann. Sci. \'Ecole Norm. Sup. (4)}, 40(3):387--444, 2007.

\bibitem{toenvezzosi-remarkonKtheory}
B.~To{\"e}n and G.~Vezzosi.
\newblock A remark on {$K$}-theory and {$S$}-categories.
\newblock {\em Topology}, 43(4):765--791, 2004.

\bibitem{Voevodsky-icm}
V.~Voevodsky.
\newblock {$\bold A^1$}-homotopy theory.
\newblock In {\em Proceedings of the {I}nternational {C}ongress of
  {M}athematicians, {V}ol. {I} ({B}erlin, 1998)}, number Extra Vol. I, pages
  579--604 (electronic), 1998.

\bibitem{waldhausen-ktheoryofspaces}
F.~Waldhausen.
\newblock Algebraic {$K$}-theory of spaces.
\newblock In {\em Algebraic and geometric topology ({N}ew {B}runswick,
  {N}.{J}., 1983)}, volume 1126 of {\em Lecture Notes in Math.}, pages
  318--419. Springer, Berlin, 1985.

\bibitem{weibel-homotopyinvariantktheory}
C.~A. Weibel.
\newblock Homotopy algebraic {$K$}-theory.
\newblock In {\em Algebraic {$K$}-theory and algebraic number theory
  ({H}onolulu, {HI}, 1987)}, volume~83 of {\em Contemp. Math.}, pages 461--488.
  Amer. Math. Soc., Providence, RI, 1989.

\end{thebibliography}

\end{document}